\documentclass[11pt,a4paper]{article}
\usepackage[mac]{inputenc}
\usepackage[english]{babel}
\usepackage{amsmath}
\usepackage{amsfonts}
\usepackage{graphicx}
\usepackage{amssymb}
\usepackage{indentfirst}
\usepackage{dsfont}
\usepackage{hyperref}
\usepackage{color}
\usepackage{pgf,pgfnodes,pgfarrows,pgfshade}

\pgfdeclareimage[interpolate=true,height=5.5cm,width=7.5cm]{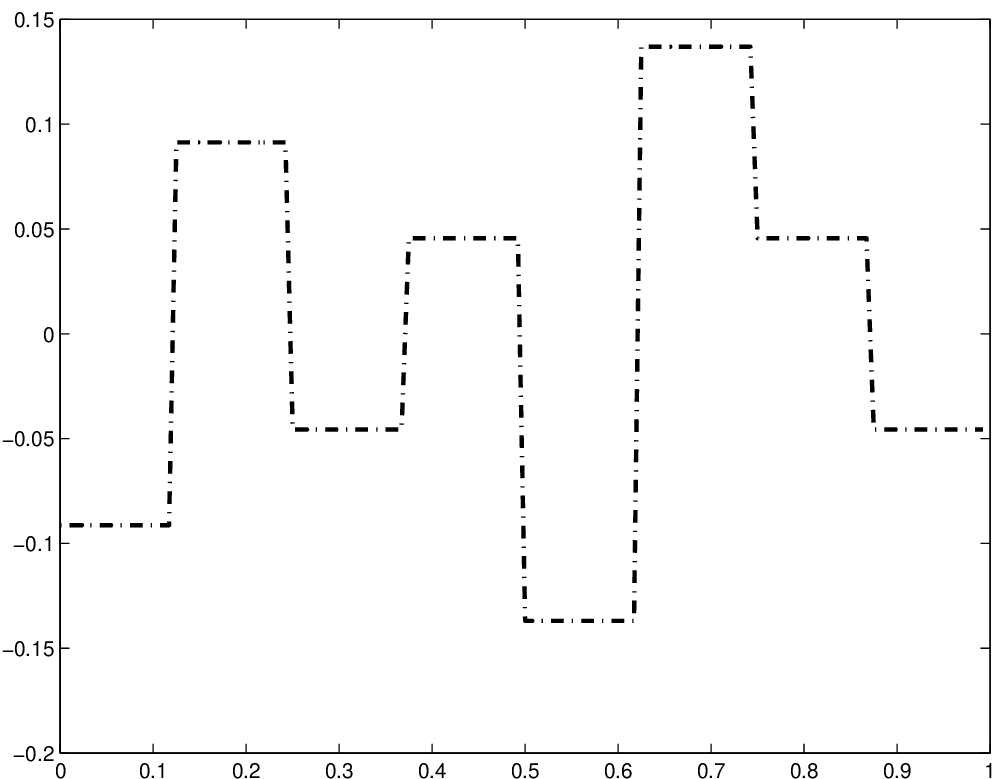}{f1}
\pgfdeclareimage[interpolate=true,height=5.5cm,width=7.5cm]{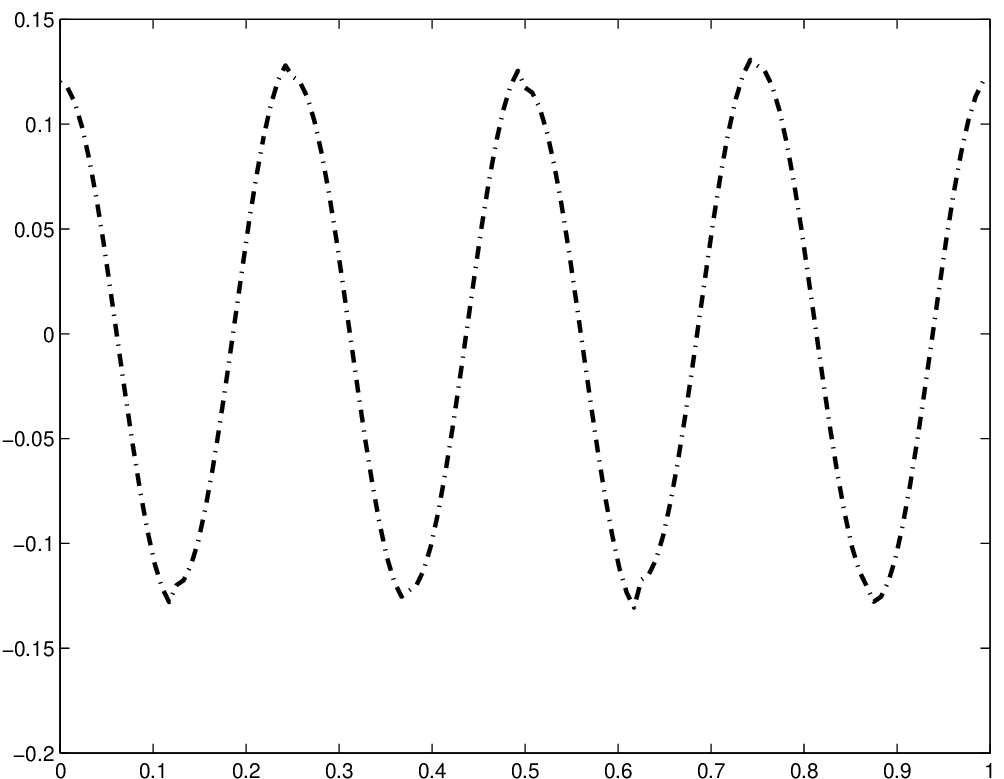}{f2}
\pgfdeclareimage[interpolate=true,height=5.5cm,width=7.5cm]{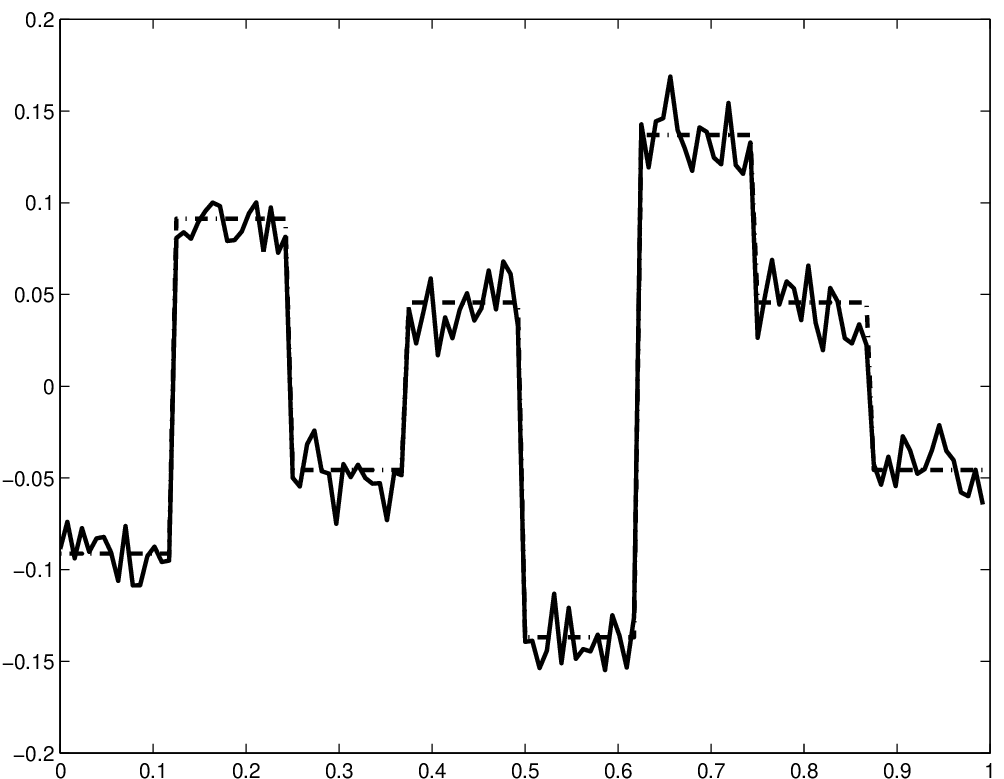}{hatf1_naive}
\pgfdeclareimage[interpolate=true,height=5.5cm,width=7.5cm]{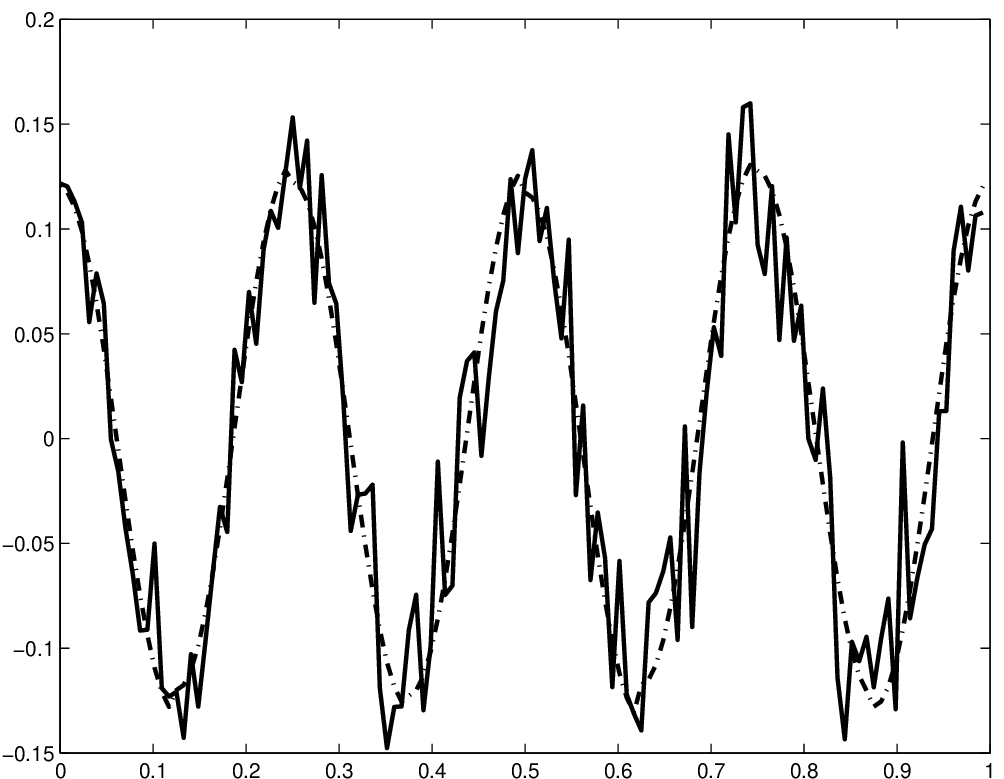}{hatf2_naive}
\pgfdeclareimage[interpolate=true,height=5.5cm,width=7.5cm]{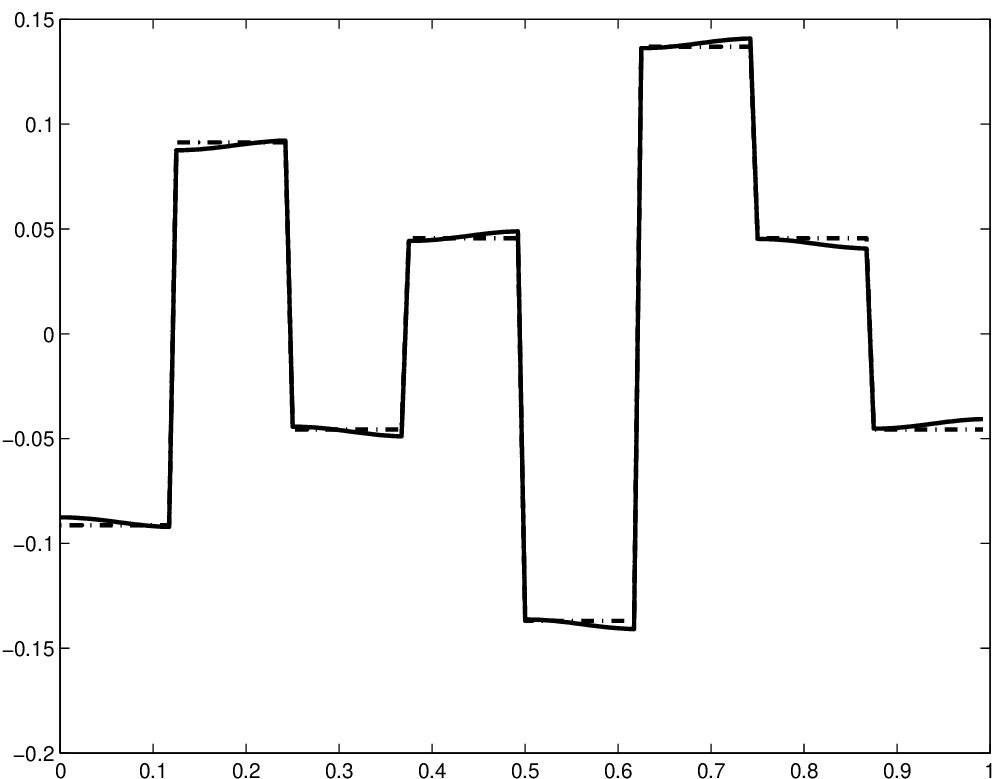}{hatf1_GLF}
\pgfdeclareimage[interpolate=true,height=5.5cm,width=7.5cm]{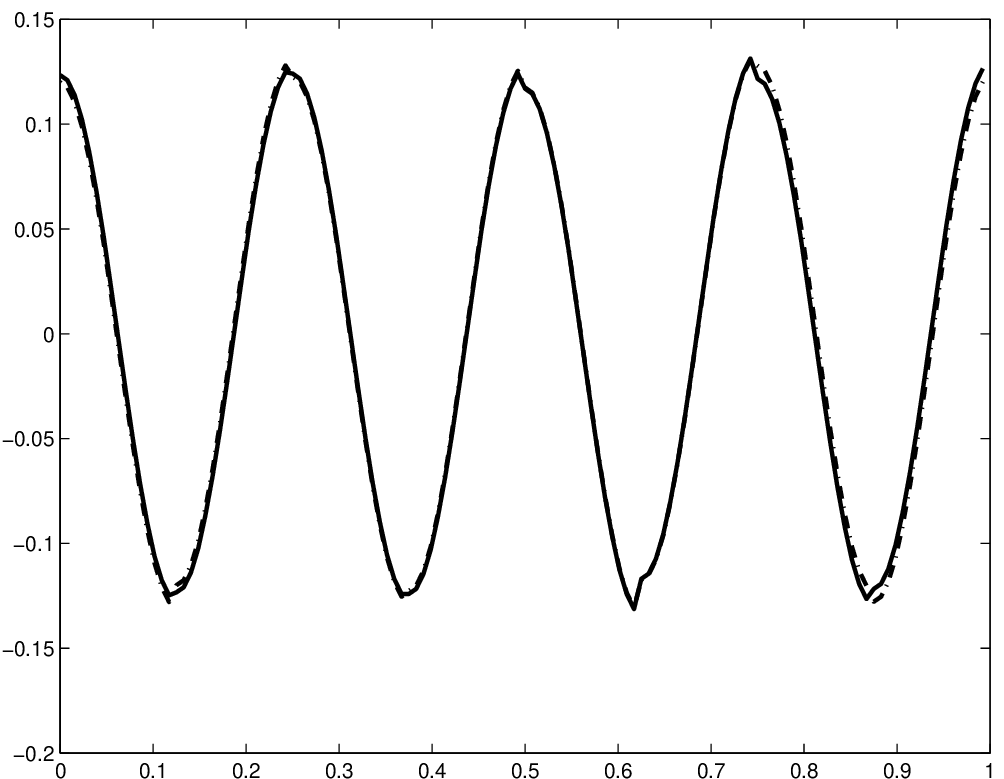}{hatf2_GLF}
\pgfdeclareimage[interpolate=true,height=5.5cm,width=7.5cm]{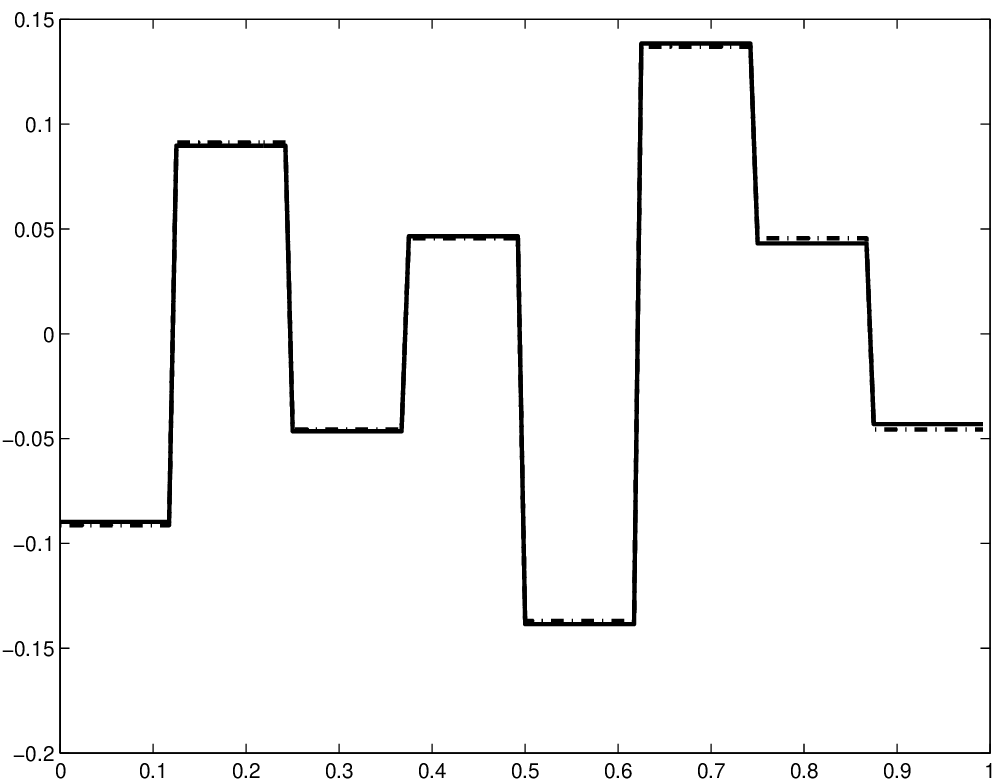}{hatf1_Haar}
\pgfdeclareimage[interpolate=true,height=5.5cm,width=7.5cm]{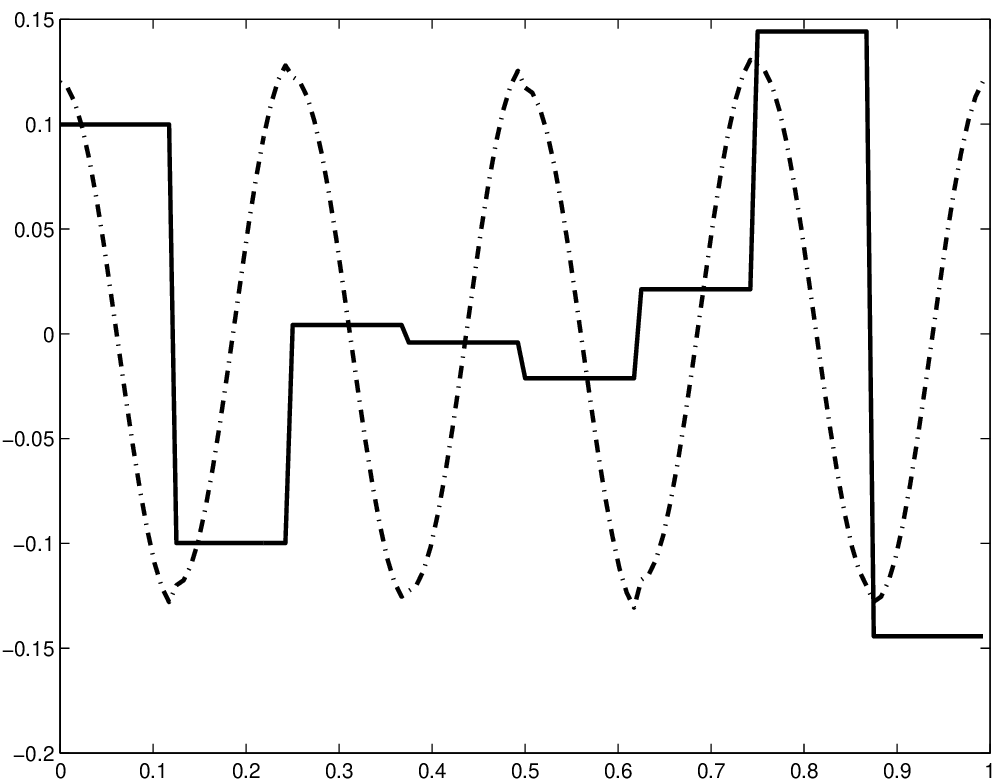}{hatf2_Haar}
\pgfdeclareimage[interpolate=true,height=5.5cm,width=7.5cm]{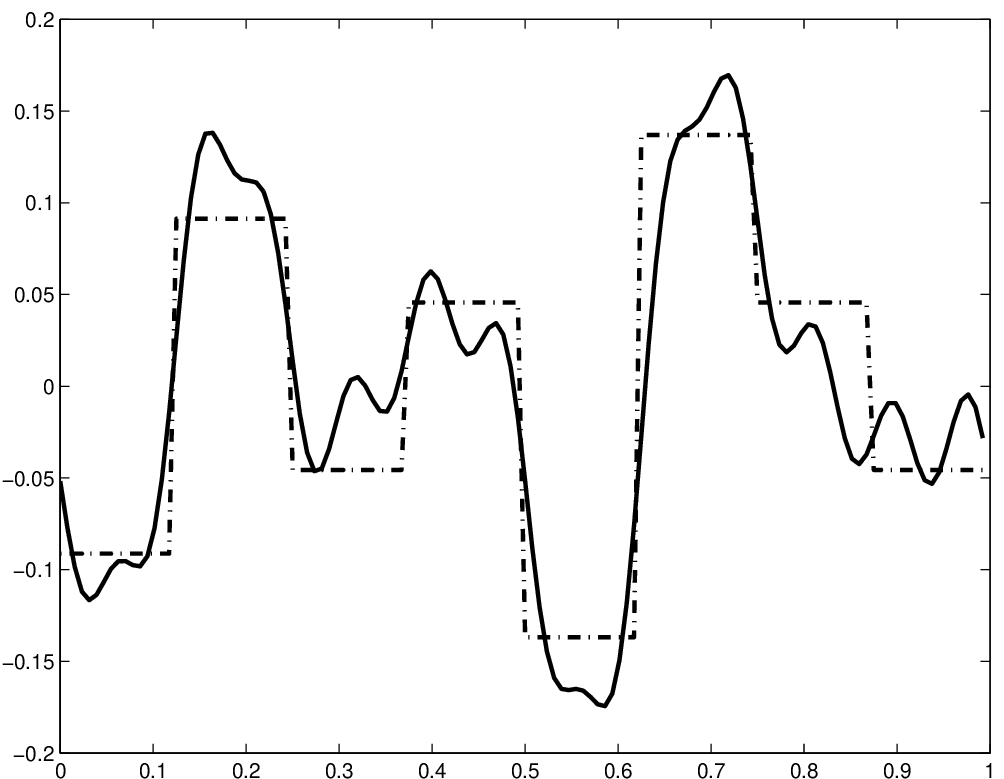}{hatf1_Fourier}
\pgfdeclareimage[interpolate=true,height=5.5cm,width=7.5cm]{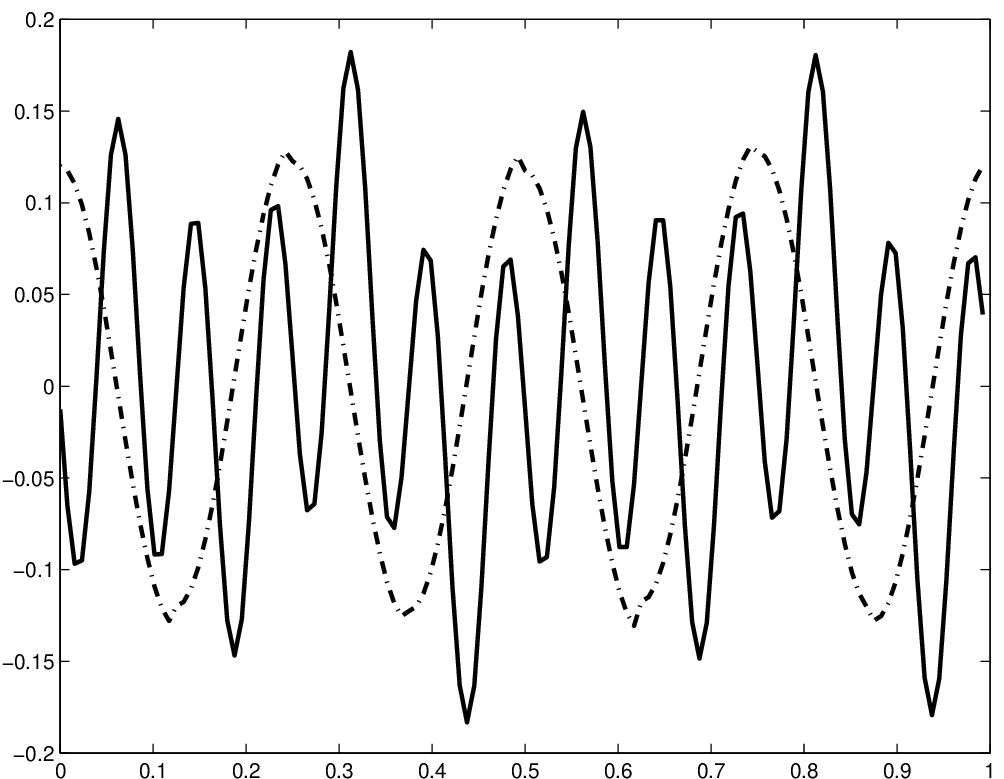}{hatf2_Fourier}
\pgfdeclareimage[interpolate=true,height=5.5cm,width=7.5cm]{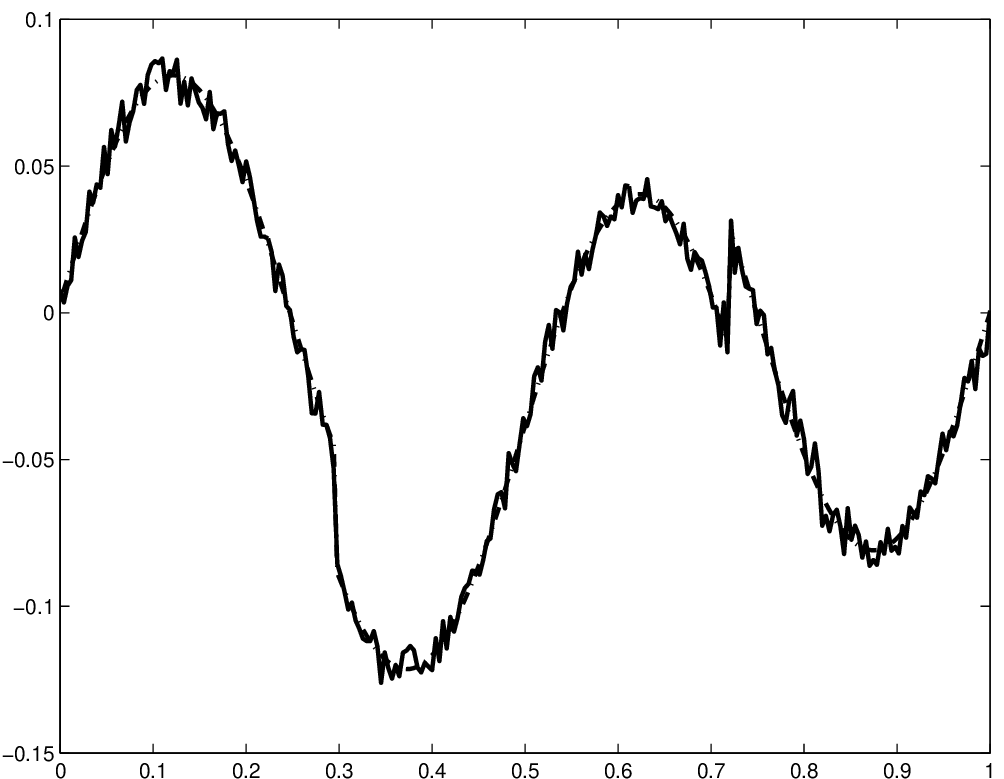}{Fig1H}
\pgfdeclareimage[interpolate=true,height=5.5cm,width=7.5cm]{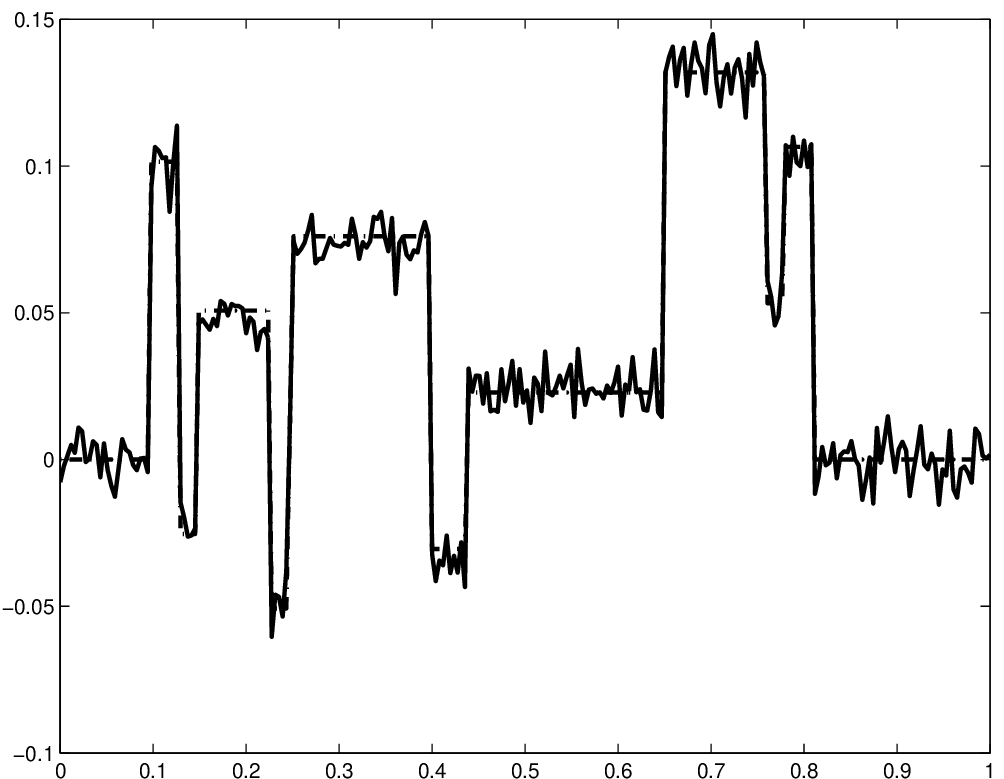}{Fig1B}
\pgfdeclareimage[interpolate=true,height=5.5cm,width=7.5cm]{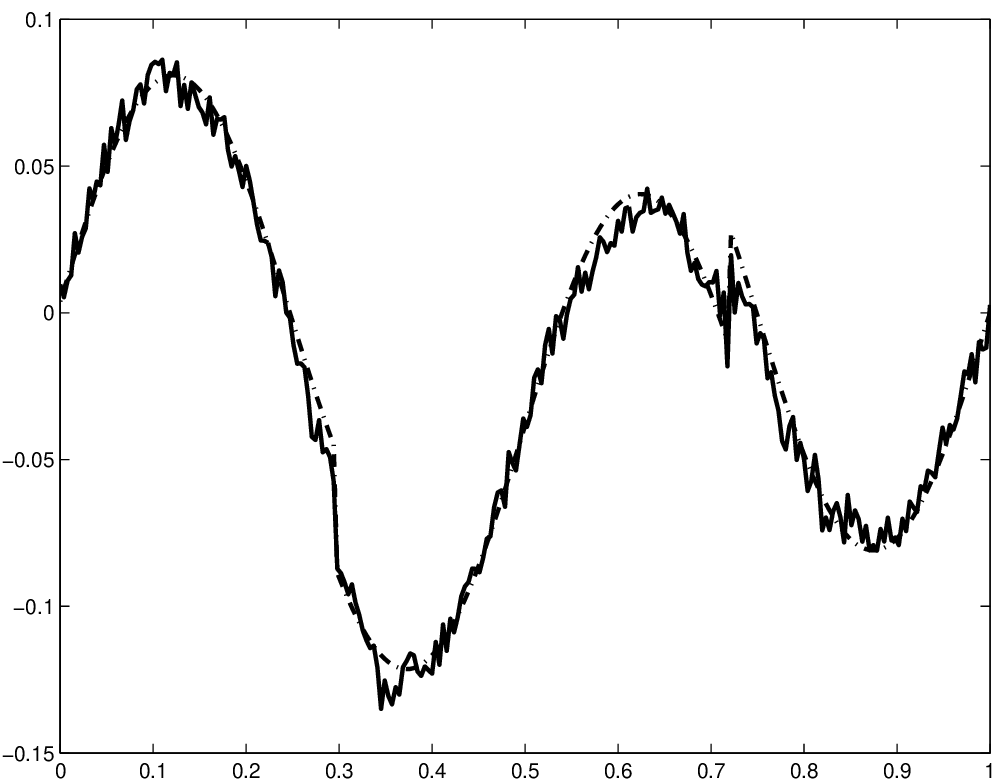}{Fig2H}
\pgfdeclareimage[interpolate=true,height=5.5cm,width=7.5cm]{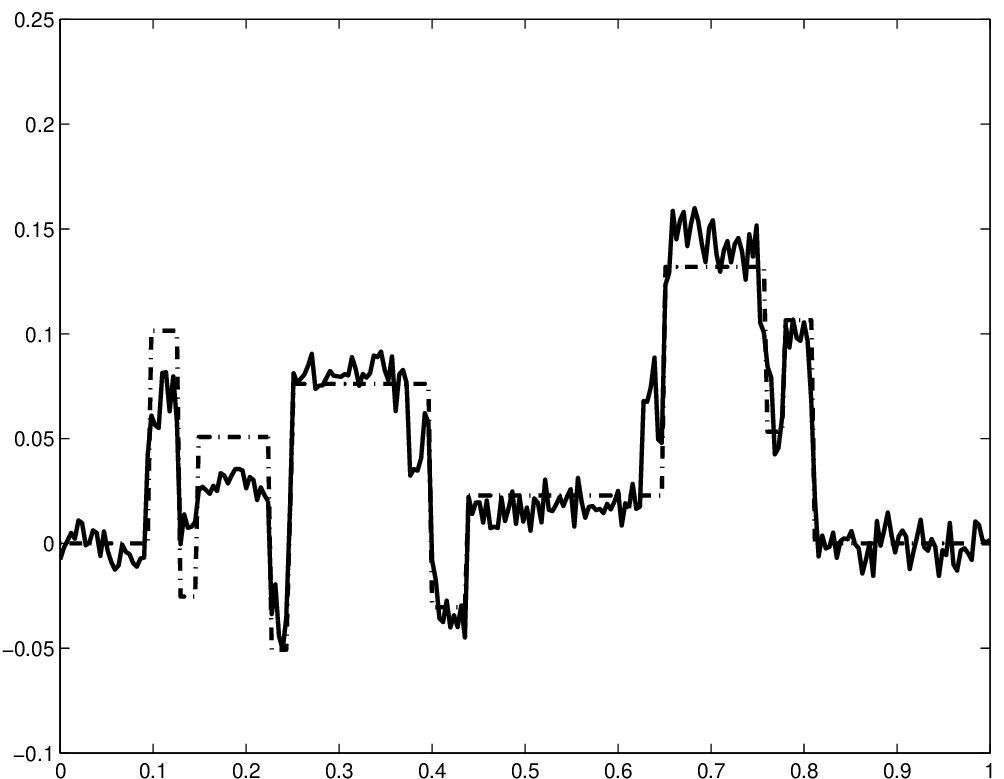}{Fig2B}
\pgfdeclareimage[interpolate=true,height=5.5cm,width=7.5cm]{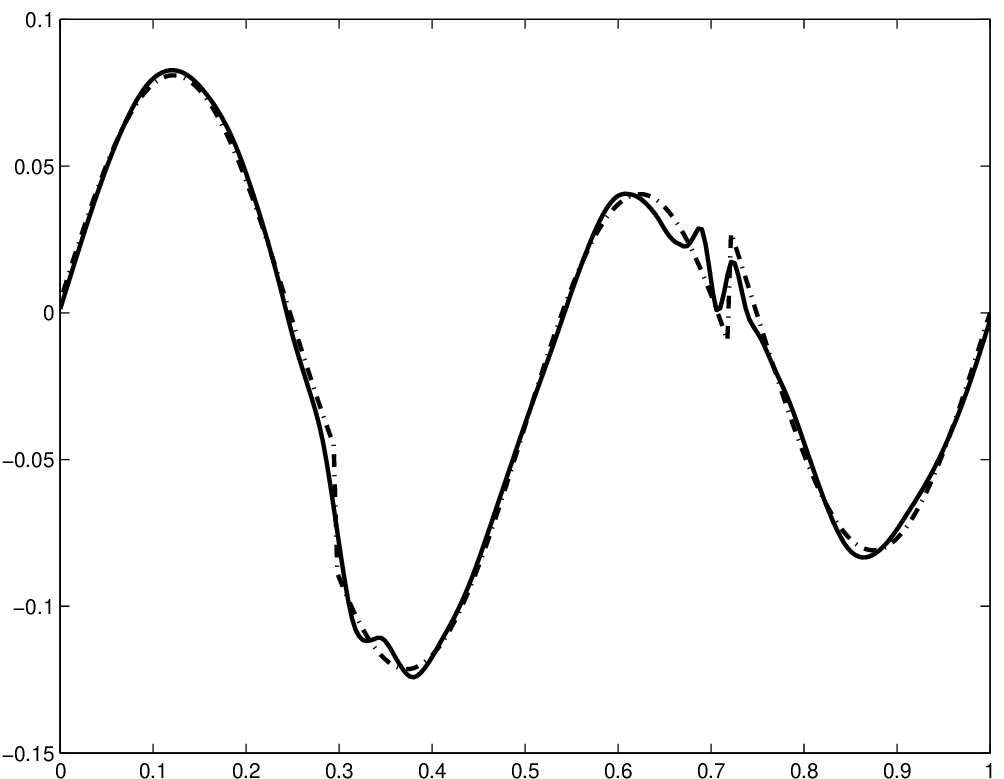}{Fig3H}
\pgfdeclareimage[interpolate=true,height=5.5cm,width=7.5cm]{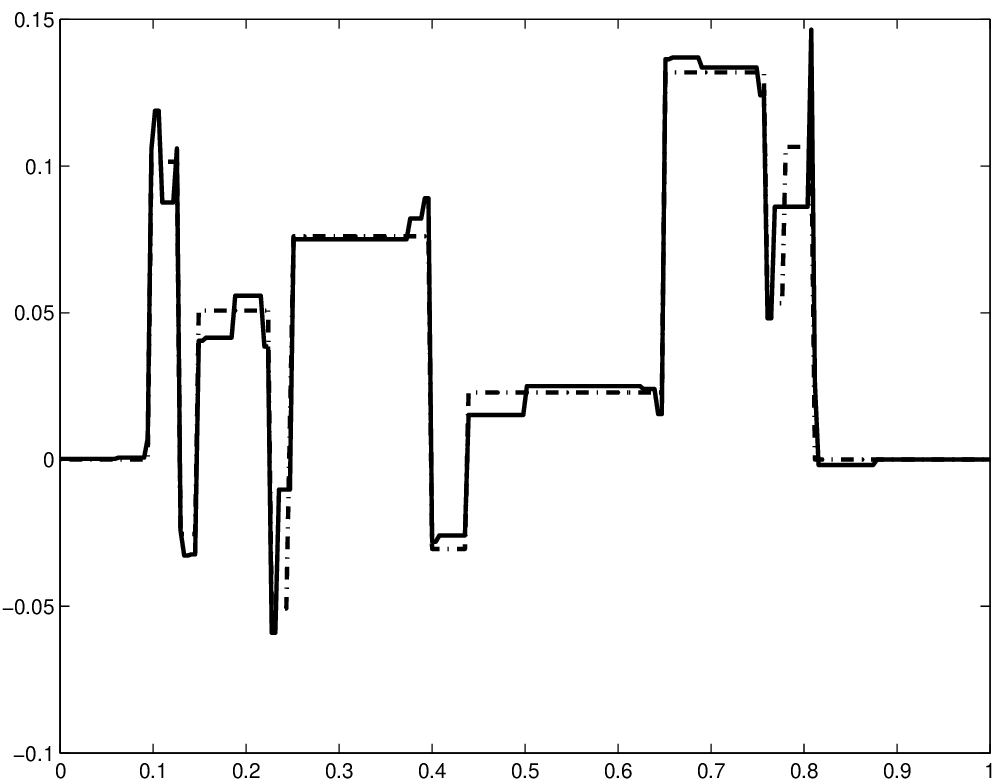}{Fig3B}
\pgfdeclareimage[interpolate=true,height=5.5cm,width=7.5cm]{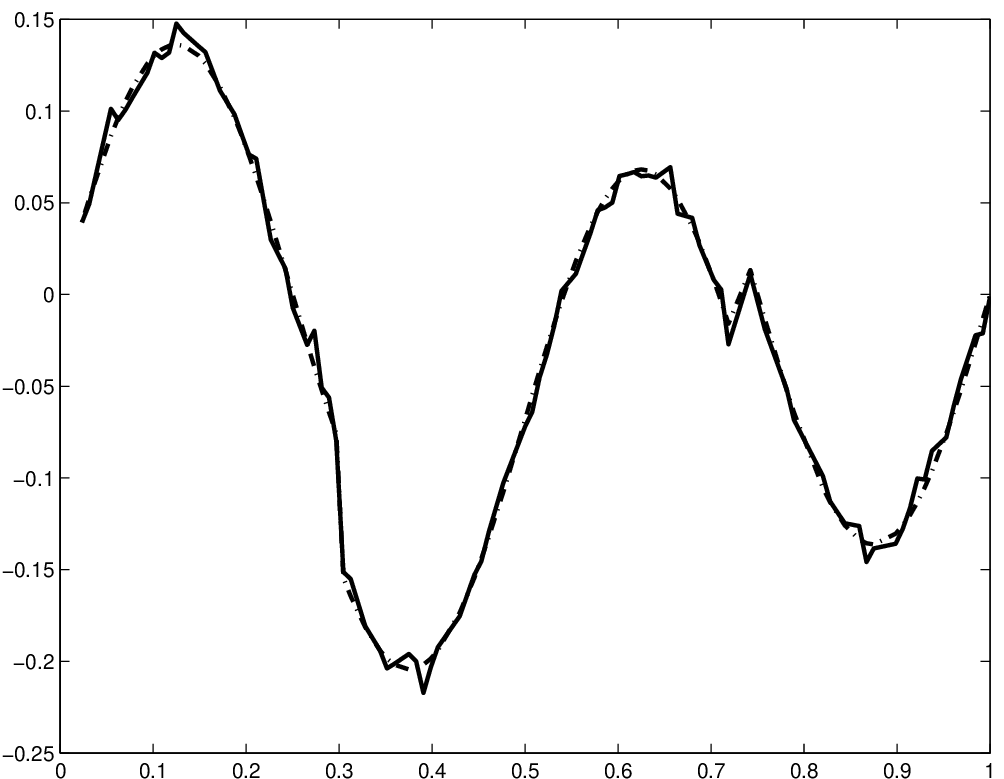}{Fig4Hnn}
\pgfdeclareimage[interpolate=true,height=5.5cm,width=7.5cm]{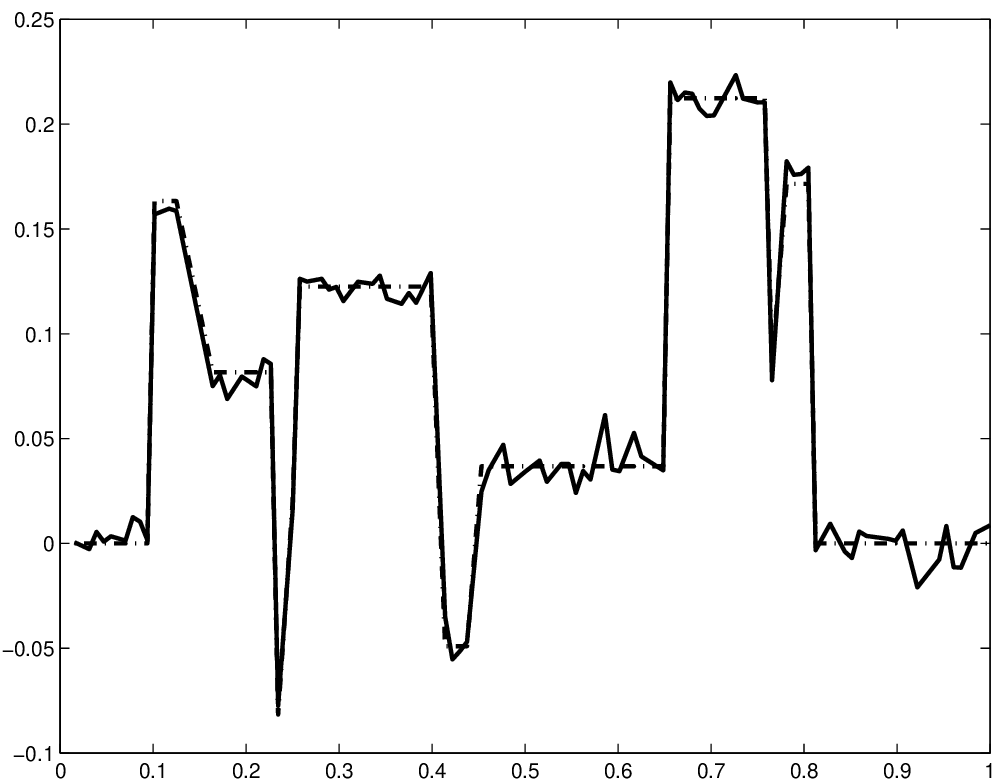}{Fig4Bn}
\pgfdeclareimage[interpolate=true,height=5.5cm,width=7.5cm]{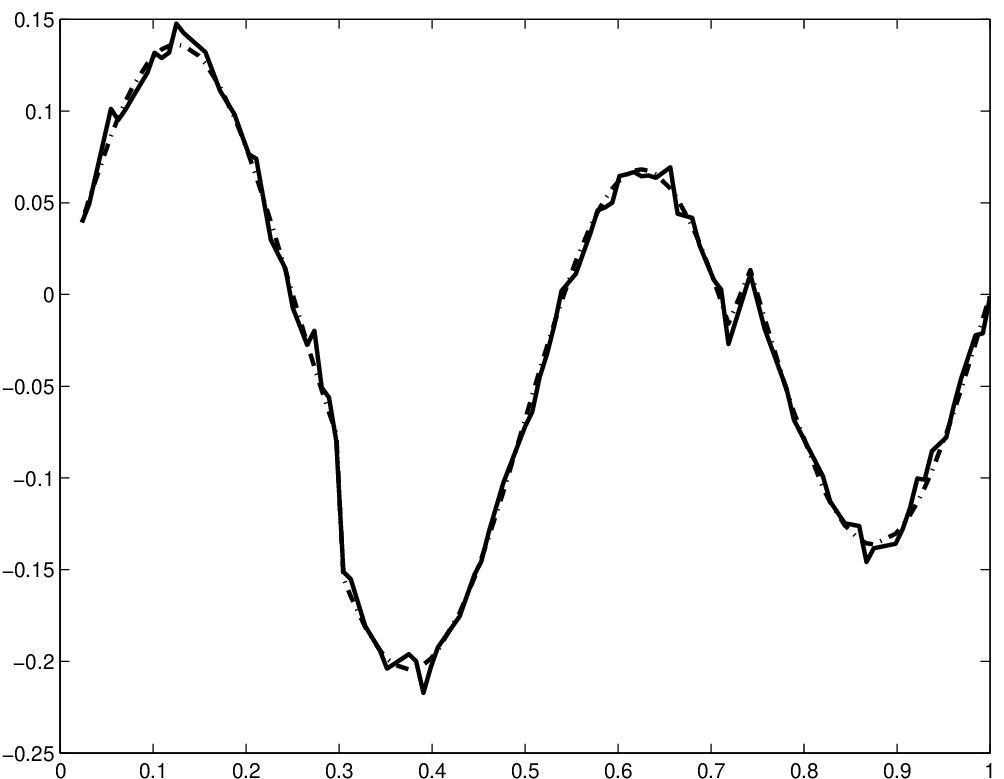}{Fig5Hnn}
\pgfdeclareimage[interpolate=true,height=5.5cm,width=7.5cm]{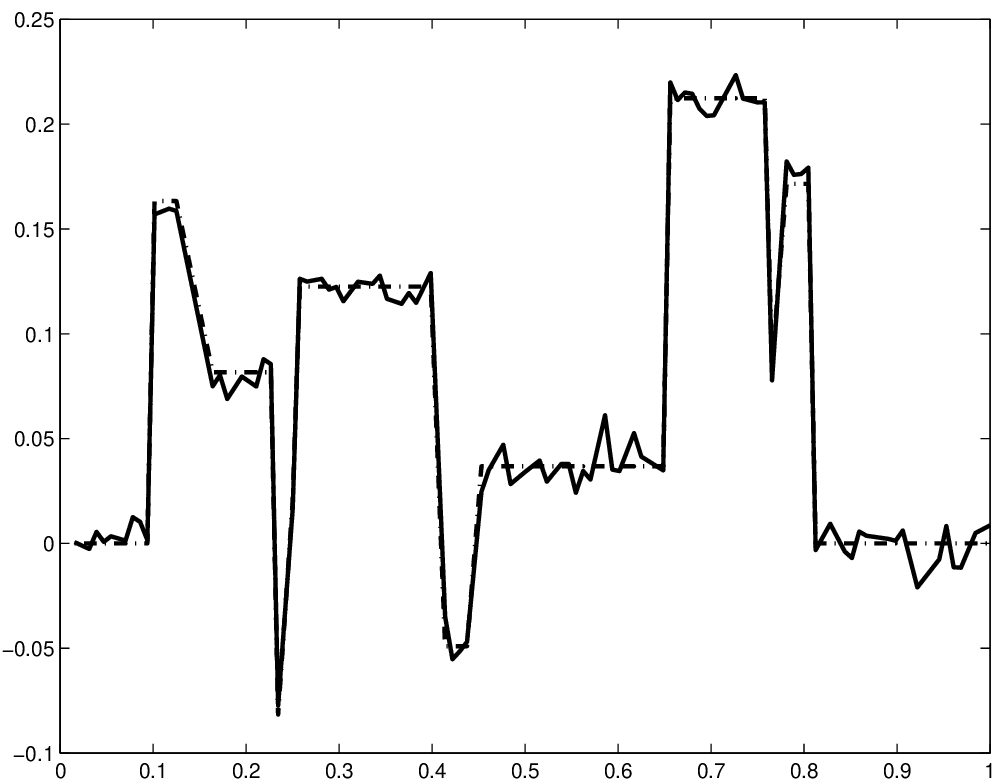}{Fig5Bn}
\pgfdeclareimage[interpolate=true,height=5.5cm,width=7.5cm]{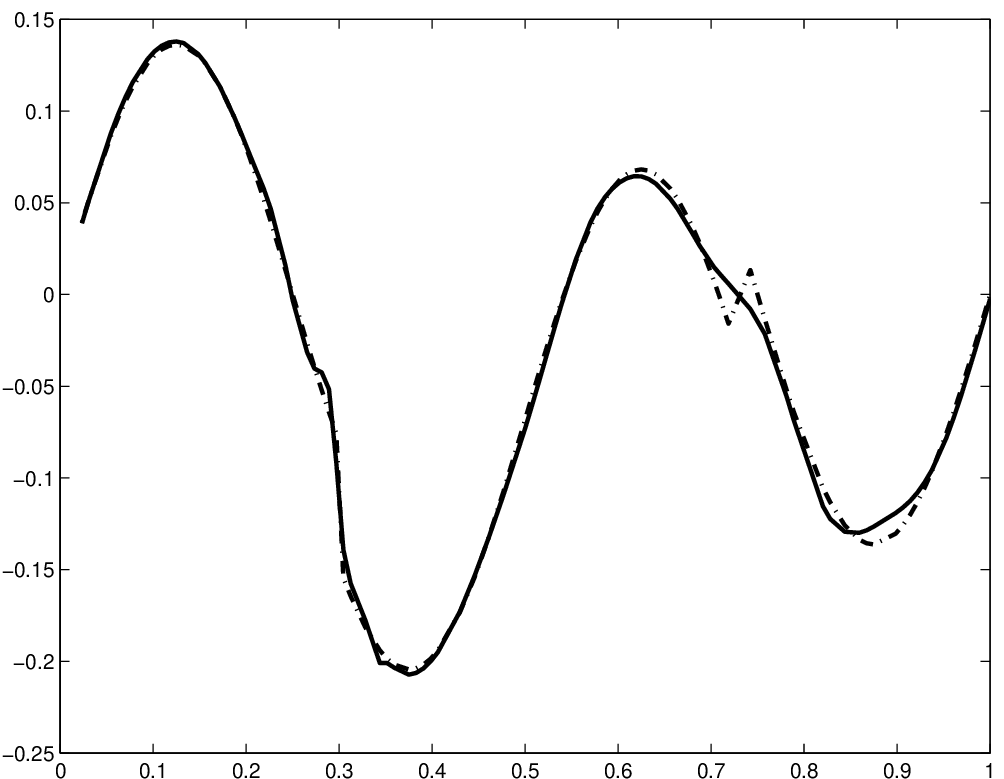}{Fig6Hnn}
\pgfdeclareimage[interpolate=true,height=5.5cm,width=7.5cm]{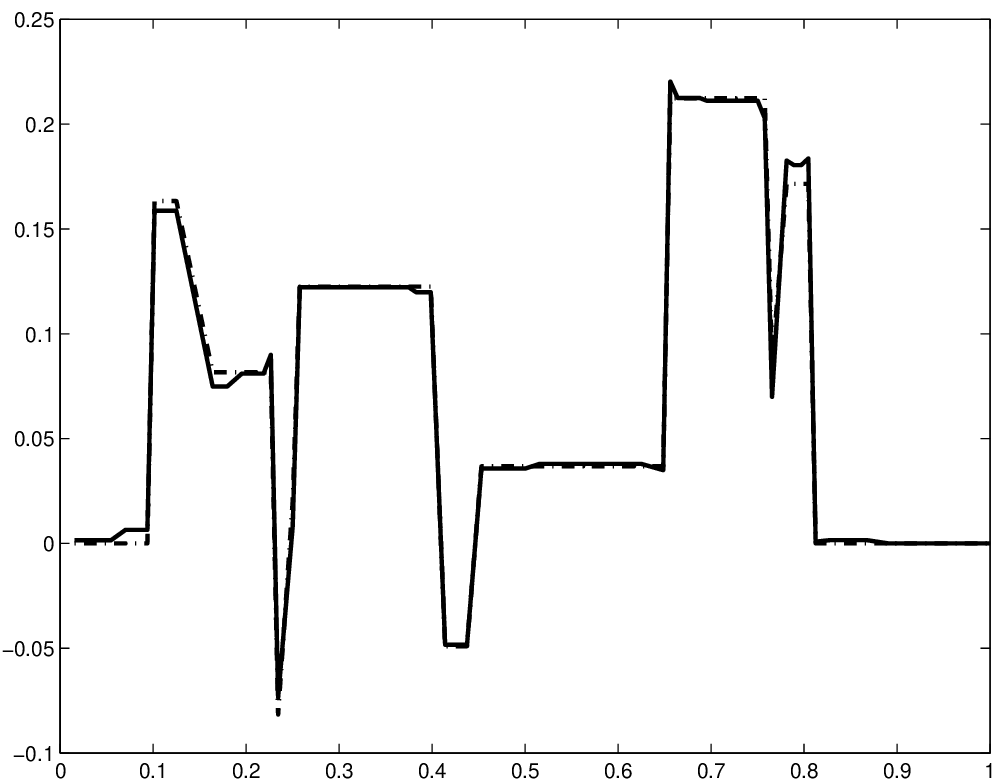}{Fig6Bn}
\pgfdeclareimage[interpolate=true,height=5.5cm,width=7.5cm]{Fig7H}{Fig7H}
\pgfdeclareimage[interpolate=true,height=5.5cm,width=7.5cm]{Fig7B}{Fig7B}
\pgfdeclareimage[interpolate=true,height=5.5cm,width=7.5cm]{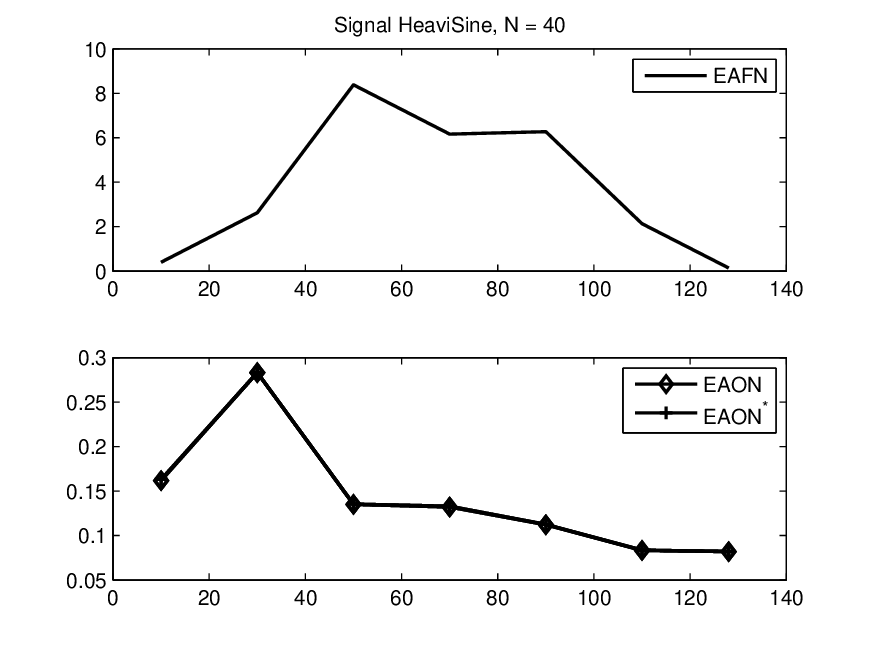}{FigHN40}
\pgfdeclareimage[interpolate=true,height=5.5cm,width=7.5cm]{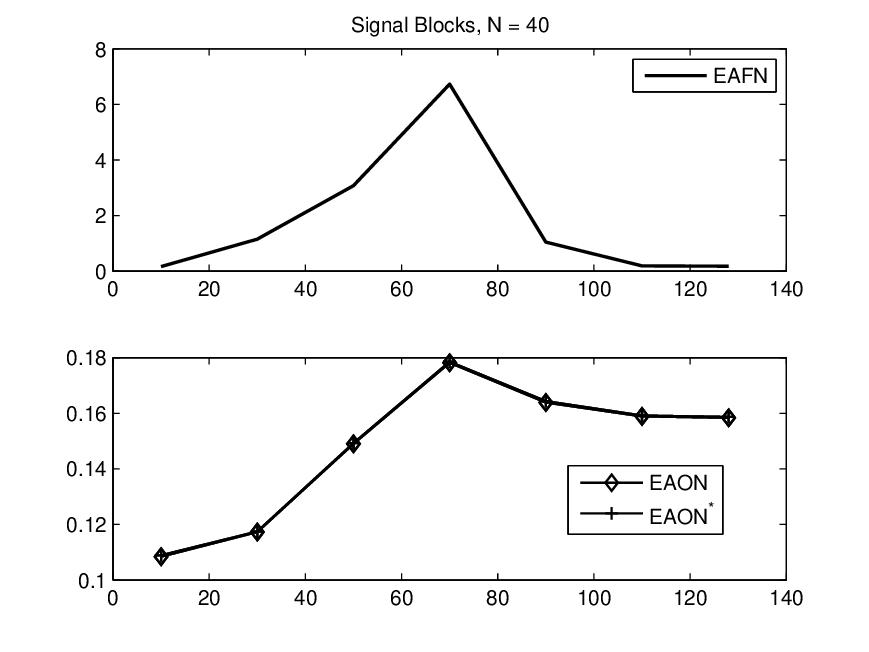}{FigBN40}
\newtheorem{theorem}{Theorem}
\newtheorem{lemma}{Lemma}
\newtheorem{proposition}{Proposition}
\newtheorem{definition}{Definition}
\newtheorem{corollary}{Corollary}
\newtheorem{assump}{Assumption}

\numberwithin{equation}{section}
\newcommand{\qed}{\nobreak \ifvmode \relax \else
      \ifdim\lastskip<1.5em \hskip-\lastskip
      \hskip1.5em plus0em minus0.5em \fi \nobreak
      \vrule height0.75em width0.5em depth0.25em\fi}

\def\R{\mathbb{R}}

\def\bX{\mathbf{X}}

\def\bSigma{\mathbf{\Sigma}}

\def\E{\mathbb{E}}

\DeclareMathOperator*{\argmin}{\mbox{argmin}}
\begin{document}

\author{J\'{e}r\'{e}mie Bigot$^{1,2}$, Rolando J. Biscay$^{4}$, Jean-Michel
Loubes$^{1}$ and Lilian Mu\~{n}iz-Alvarez$^{1,3}$ \\
\\
IMT, Universit\'{e} Paul Sabatier, Toulouse, France$^{1}$\\
Center for Mathematical Modelling, Universidad de Chile, Santiago, Chile$%
^{2} $ \\
Facultad de Matem\'{a}tica y Computaci\'{o}n, Universidad de La Habana, Cuba$%
^{3}$\\
DEUV-CIMFAV, Facultad de Ciencias, Universidad de Valparaiso, Chile$^{4}$ }
\title{Group Lasso estimation of high-dimensional covariance matrices}
\maketitle

\begin{abstract}
In this paper, we consider the Group Lasso estimator of the covariance
matrix of a stochastic process corrupted by an additive noise. We propose to
estimate the covariance matrix in a high-dimensional setting under the
assumption that the process has a sparse representation in a large
dictionary of basis functions. Using a matrix regression model, we propose a
new methodology for high-dimensional covariance matrix estimation based on
empirical contrast regularization by a group Lasso penalty. Using such a
penalty, the method selects a sparse set of basis functions in the
dictionary used to approximate the process, leading to an approximation of
the covariance matrix into a low dimensional space. Consistency of the
estimator is studied in Frobenius and operator norms and an application to
sparse PCA is proposed.
\end{abstract}

\noindent \textbf{Subject Class. MSC-2000} : 62G05, 62H25 \vskip .1in  \textwidth 16cm %
\evensidemargin 0cm \oddsidemargin 0cm \sloppy

\begin{small}
\noindent {\bf Acknowledgments:} this work was supported in part by Egide, under the Program of Eiffel
excellency Phd grants, as well as by the BDI CNRS grant. J. Bigot would like
to thank the Center for Mathematical Modeling and the CNRS for financial
support and excellent hospitality while visiting Santiago where part of this
work was carried out.
\end{small}

\section{Introduction}

Let $\mathbb{T}$ be some subset of $\mathbb{R}^{p}$, $p\in \mathbb{N}$, and
let $X=\left( X\left( t\right) \right) _{t\in \mathbb{T}}$ be a stochastic
process with values in $\mathbb{R}$. Assume that $X$ has zero mean $\mathbb{E%
}\left( X\left( t\right) \right) =0$ for all $t\in \mathbb{T}$, and finite
covariance $\sigma \left( s,t\right) =\mathbb{E}\left( X\left( s\right)
X\left( t\right) \right) $ for all $s,t\in \mathbb{T}$. Let $t_{1},\ldots
,t_{n}$ be fixed points in $\mathbb{T}$ (deterministic design), $%
X_{1},...,X_{N}$ independent copies of the process $X$, and suppose that we
observe the noisy processes 
\begin{equation}
\widetilde{X}_{i}\left( t_{j}\right) =X_{i}\left( t_{j}\right) +\mathcal{E}%
_{i}\left( t_{j}\right) \mbox{ for }i=1,...,N,\;j=1,...,n,
\label{eq:modnoise}
\end{equation}%
where $\mathcal{E}_{1},...,\mathcal{E}_{N}$ are independent copies of a
second order Gaussian process $\mathcal{E}$ with zero mean and independent
of $X$, which represent an additive source of noise in the measurements.
Based on the noisy observations (\ref{eq:modnoise}), an important problem in
statistics is to construct an estimator of the covariance matrix $\mathbf{%
\Sigma }=\mathbb{E}\left( \mathbf{X}\mathbf{X}^{\top }\right) $ of the
process $X$ at the design points, where $\mathbf{X}=\left( X\left(
t_{1}\right) ,...,X\left( t_{n}\right) \right) ^{\top }$. This problem is a
fundamental issue in many applications, ranging from geostatistics,
financial series or epidemiology for instance (see \cite{092462100}, \cite%
{MR0456314} or \cite{MR1239641,MR1741979} for general references and applications).
Estimating such a covariance matrix has also important applications in
dimension reduction by principal component analysis (PCA) or classification
by linear or quadratic discriminant analysis (LDA and QDA).

In \cite{BBLM10}, using $N$ independent copies of the process $X$, we have
proposed to construct an estimator of the covariance matrix $\mathbf{\Sigma }
$ by expanding the process $X$ into a dictionary of basis functions. The
method in \cite{BBLM10} is based on model selection techniques by empirical
contrast minimization in a suitable matrix regression model. This new
approach to covariance estimation is well adapted to the case of
low-dimensional covariance estimation when the number of replicates $N$ of
the process is larger than the number of observations points $n$. However,
many application areas are currently dealing with the problem of estimating
a covariance matrix when the number of observations at hand is small when
compared to the number of parameters to estimate. Examples include
biomedical imaging, proteomic/genomic data, signal processing in
neurosciences and many others. This issue corresponds to the problem of
covariance estimation for high-dimensional data. This problem is challenging
since, in a high-dimensional setting (when $n>>N$ or $n\sim N$), it is well
known that the sample covariance matrices 
\begin{equation*}
\mathbf{S}=\frac{1}{N}\sum\limits_{i=1}^{N}\mathbf{X}_{i}\mathbf{X}%
_{i}^{\top }\in \mathbb{R}^{n\times n},\mbox{ where }\mathbf{X}_{i}=\left(
X_{i}\left( t_{1}\right) ,...,X_{i}\left( t_{n}\right) \right) ^{\top
},i=1,\ldots ,N
\end{equation*}%
and 
\begin{equation*}
\widetilde{\mathbf{S}}=\frac{1}{N}\sum\limits_{i=1}^{N}\widetilde{\mathbf{X}}%
_{i}\widetilde{\mathbf{X}}_{i}^{\top }\in \mathbb{R}^{n\times n},%
\mbox{
where }\widetilde{\mathbf{X}}_{i}=\left( \widetilde{X}_{i}\left(
t_{1}\right) ,...,\widetilde{X}_{i}\left( t_{n}\right) \right) ^{\top
},i=1,\ldots ,N
\end{equation*}%
behave poorly, and are not consistent estimators of $\mathbf{\Sigma }$. For
example, suppose that the $\mathbf{X}_{i}$'s are independent and identically
distributed (i.i.d.) random vectors in $\mathbb{R}^{n}$ drawn from a
multivariate Gaussian distribution. Then, when $\frac{n}{N}\rightarrow c>0$
as $n,N\rightarrow +\infty $, neither the eigenvalues nor the eigenvectors
of the sample covariance matrix $\mathbf{S}$ are consistent estimators of
the eigenvalues and eigenvectors of $\mathbf{\Sigma }$ (see \cite{MR1863961}%
). This topic has thus recently received a lot of attention in the
statistical literature. To achieve consistency, recently developed methods
for high-dimensional covariance estimation impose sparsity restrictions on
the matrix $\mathbf{\Sigma }$. Such restrictions imply that the true (but
unknown) dimension of the model is much lower than the number $\frac{n(n+1)}{%
2}$ of parameters of an unconstrained covariance matrix. Under various
sparsity assumptions, different regularizing methods of the empirical
covariance matrix have been proposed. Estimators based on thresholding or
banding the entries of the empirical covariance matrix have been studied in 
\cite{MR2485008} and \cite{MR2387969}. Thresholding the components of the
empirical covariance matrix has also been proposed by \cite{MR2485011} and
the consistency of such estimates is studied using tools from random matrix
theory. \cite{Fan-etal08} impose sparsity on the covariance via a factor
model which is appropriate in financial applications. \cite{MR2415602} and 
\cite{MR2417391} propose regularization techniques with a Lasso penalty to
estimate the covariance matrix or its inverse. More general penalties have
been studied in \cite{MR2572459}. Another approach is to impose sparsity on
the eigenvectors of the covariance matrix which leads to sparse PCA. \cite%
{MR2252527} use a Lasso penalty to achieve sparse representation in PCA, 
\cite{MR2426043} study properties of sparse principal components by convex
programming, while \cite{citeulike:7339057} propose a PCA regularization by
expanding the empirical eigenvectors in a sparse basis and then apply a
thresholding step.

In this paper, we propose to estimate $\mathbf{\Sigma}$ in a
high-dimensional setting by using the assumption that the process $X$ has a
sparse representation in a large dictionary of basis functions. Using a
matrix regression model as in \cite{BBLM10}, we propose a new methodology
for high-dimensional covariance matrix estimation based on empirical
contrast regularization by a group Lasso penalty. Using such a penalty, the
method selects a sparse set of basis functions in the dictionary used to
approximate the process $X$. This leads to an approximation of the
covariance matrix $\mathbf{\Sigma}$ into a low dimensional space, and thus
to a new method of dimension reduction for high-dimensional data. Group
Lasso estimators have been studied in the standard linear model and in
multiple kernel learning to impose a group-sparsity structure on the
parameters to recover (see \cite{MR2426104}, \cite{MR2417268} and references
therein). However, to the best of our knowledge, it has not been used for
the estimation of covariance matrices using a functional approximation of
the process $X$.

The rest of the paper is organized as follows. In Section \ref{sec:model},
we describe a matrix regression model for covariance estimation, and we
define our estimator by group Lasso regularization. The consistency of such
a procedure is investigated in Section \ref{sec:oracle} using oracle
inequalities and a non-asymptotic point of view by holding fixed the number
of replicates $N$ and observation points $n$. Consistency of the estimator
is studied in Frobenius and operator norms. Various results existing in
matrix theory show that convergence in operator norm implies convergence of
the eigenvectors and eigenvalues (e.g.\ through the use of the $\sin(\theta)$
theorems in \ \cite{MR0264450}). Consistency in operator norm is thus well
suited for PCA applications. Numerical experiments are given in Section \ref%
{sec:simus}, and an application to sparse PCA is proposed. A technical
Appendix contains all the proofs.

\section{Model and definition of the estimator}

\label{sec:model}

To impose sparsity restrictions on the covariance matrix $\mathbf{\Sigma}$,
our approach is based on an approximation of the process in a finite
dictionary of (not necessarily orthogonal) basis functions $g_{m}: \mathbb{T}%
\rightarrow \mathbb{R}$ for $m=1,...,M$. Suppose that

\begin{equation}
X\left( t\right) \approx \sum\limits_{m=1}^{M}{}a_{m}g_{m}\left( t\right) ,
\label{model}
\end{equation}%
where $a_{m}$, $m=1,...,M$ are real valued random variables, and that for
each trajectory $X_{i}$ 
\begin{equation}
X_{i}\left( t_{j}\right) \approx \sum\limits_{m=1}^{M}{}a_{i,m}g_{m}\left(
t_{j}\right) .  \label{empmodel}
\end{equation}%
The notation $\approx $ means that the process $X$ can be well approximated
into the dictionary. A precise meaning of this will be discussed later on.
Then (\ref{empmodel}) can be written in matrix notation as:%
\begin{equation}
\mathbf{X}_{i}\approx \mathbf{G}\mathbf{a}_{i},\;i=1,...,N
\label{empmodelmatrixnotation}
\end{equation}%
where $\mathbf{G}$ is the $n\times M$ matrix with entries 
\begin{equation*}
\mathbf{G}_{jm}=g_{m}\left( t_{j}\right) \mbox{ for }1\leq j\leq n%
\mbox{
and }1\leq m\leq M,
\end{equation*}%
and $\mathbf{a}_{i}$ is the $M\times 1$ random vector of components $a_{i,m}$%
, with $1\leq m\leq M$.

Recall that we want to estimate the covariance matrix $\mathbf{\Sigma }=%
\mathbb{E}\left( \mathbf{X}\mathbf{X}^{\top }\right) $ from the noisy
observations (\ref{eq:modnoise}). Since $\mathbf{X}\approx \mathbf{G}\mathbf{%
a}$ with $\mathbf{a}=\left( a_{m}\right) _{1\leq m\leq M}$ with $a_{m}$ as
in (\ref{model}), it follows that 
\begin{equation*}
\mathbf{\Sigma }\approx \mathbb{E}\left( \mathbf{G}\mathbf{a}\left( \mathbf{G%
}\mathbf{a}\right) ^{\top }\right) =\mathbb{E}\left( \mathbf{G}\mathbf{a}%
\mathbf{a}^{\top }\mathbf{G}^{\top }\right) =\mathbf{G}\mathbf{\Psi }^{\ast }%
\mathbf{G}^{\top }\mbox{ with }\mathbf{\Psi }^{\ast }=\mathbb{E}\left( 
\mathbf{a}\mathbf{a}^{\top }\right) .
\end{equation*}

Given the noisy observations $\widetilde{\mathbf{X}}_{i}\ $as in (\ref%
{eq:modnoise}) with $i=1,...,N$, consider the following matrix regression
model 
\begin{equation}
\widetilde{\mathbf{X}}_{i}\widetilde{\mathbf{X}}_{i}^{\top }=\mathbf{\Sigma }%
+\mathbf{U}_{i}+\mathbf{W}_{i}\;i=1,\ldots ,N, \label{eq:regmatrix}
\end{equation}%
where $\mathbf{U}_{i}=\mathbf{X}_{i}\mathbf{X}_{i}^{\top }-\mathbf{\Sigma }$
are i.i.d centered matrix errors, and 
\begin{equation*}
\mathbf{W}_{i}=\mathcal{E}_{i}\mathcal{E}_{i}^{\top }\in \mathbb{R}^{n\times
n}\mbox{ where }\mathcal{E}_{i}=\left( \mathcal{E}_{i}\left( t_{1}\right)
,...,\mathcal{E}_{i}\left( t_{n}\right) \right) ^{\top },i=1,\ldots ,N.
\end{equation*}%
The size $M$ of the dictionary can be very large, but it is expected that
the process $X$ has a sparse expansion in this basis, meaning that, in
approximation (\ref{model}), many of the random coefficients $a_{m}$ are
close to zero. We are interested in obtaining an estimate of the covariance $%
\mathbf{\Sigma }$ in the form $\widehat{\mathbf{\Sigma }}=\mathbf{G}\widehat{%
\mathbf{\Psi }}\mathbf{G}^{\top }$ such that $\widehat{\mathbf{\Psi }}$ is a
symmetric $M\times M$ matrix with many zero rows (and so, by symmetry, many
corresponding zero columns). Note that setting the $k$-th row of $\widehat{%
\mathbf{\Psi }}$ to $\mathbf{0}\in \mathbb{R}^{M}$ means to remove the
function $g_{k}$ from the set of basis functions $\left( g_{m}\right)
_{1\leq m\leq M}$ in the function expansion associated to $\mathbf{G}$.

Let us now explain how to select a sparse set of rows/columns in the matrix $%
\widehat{\mathbf{\Psi }}$. For this, we use a group Lasso approach to
threshold some rows/columns of $\widehat{\mathbf{\Psi }}$ which corresponds
to removing some basis functions in the approximation of the process $X$.
For two $p\times p$ matrices $\mathbf{A},\mathbf{B}$ define the inner
product $\left\langle \mathbf{A},\mathbf{B}\right\rangle _{F}:=tr\left( 
\mathbf{A}^{\top }\mathbf{B}\right) $ and the associated Frobenius norm $%
\Vert \mathbf{A}\Vert _{F}^{2}:=tr\left( \mathbf{A}^{\top }\mathbf{A}\right) 
$. Let $\mathcal{S}_{M}$ denote the set of $M\times M$ symmetric matrices
with real entries. We define the group Lasso estimator of the covariance
matrix $\mathbf{\Sigma }$ by 
\begin{equation}
\widehat{\mathbf{\Sigma }}_{\lambda }=\mathbf{G}\widehat{\mathbf{\Psi }}%
_{\lambda }\mathbf{G}^{\top }\in \mathbb{R}^{n\times n},
\label{GroupLassoEstimator}
\end{equation}%
where $\widehat{\mathbf{\Psi }}_{\lambda }$ is the solution of the following
optimization problem:%
\begin{equation}
\widehat{\mathbf{\Psi }}_{\lambda }=\underset{\mathbf{\Psi \in }\mathcal{S}%
_{M}}{\argmin}\left\{ \frac{1}{N}\sum_{i=1}^{N}\left\Vert \widetilde{\mathbf{%
X}}_{i}\widetilde{\mathbf{X}}_{i}^{\top }-\mathbf{G\Psi G}^{\top
}\right\Vert _{F}^{2}+2\lambda \sum_{k=1}^{M}\gamma _{k}\sqrt{%
\sum_{m=1}^{M}\Psi _{mk}^{2}}\right\} ,  \label{MatrixGroupLassoEstimator0}
\end{equation}%
where $\mathbf{\Psi }=\left( \Psi _{mk}\right) _{1\leq m,k\leq M}\in \mathbb{%
R}^{M\times M}$, $\lambda $ is a positive number and $\gamma _{k}$ are some
weights whose values will be discuss later on. In (\ref%
{MatrixGroupLassoEstimator0}), the penalty term imposes to give preference
to solutions with components $\mathbf{\Psi }_{k}=\mathbf{0}$, where $\left( 
\mathbf{\Psi }_{k}\right) _{1\leq k\leq M}$ denotes the columns of $\mathbf{%
\Psi }$. Recall that $\widetilde{\mathbf{S}}=\frac{1}{N}\sum\limits_{i=1}^{N}%
\widetilde{\mathbf{X}}_{i}\widetilde{\mathbf{X}}_{i}^{\top }$ denotes the
sample covariance matrix from the noisy observations (\ref{eq:modnoise}). 
It can be checked that minimizing the criterion (\ref%
{MatrixGroupLassoEstimator0}) is equivalent to 
\begin{equation}
\widehat{\mathbf{\Psi }}_{\lambda }=\underset{\mathbf{\Psi \in }\mathcal{S}%
_{M}}{\argmin}\left\{ \left\Vert \widetilde{\mathbf{S}}-\mathbf{G\Psi G}%
^{\top }\right\Vert _{F}^{2}+2\lambda \sum_{k=1}^{M}\gamma _{k}\sqrt{%
\sum_{m=1}^{M}\Psi _{mk}^{2}}\right\} .  \label{MatrixGroupLassoEstimator}
\end{equation}%
Thus $\widehat{\mathbf{\Psi }}_{\lambda }\in \mathbb{R}^{M\times M}$ can be
interpreted as a group Lasso estimator of $\mathbf{\Sigma }$ in the
following matrix regression model%
\begin{equation}
\widetilde{\mathbf{S}}=\mathbf{\Sigma }+\mathbf{U}+\mathbf{W}\approx \mathbf{%
G}\mathbf{\Psi }^{\ast }\mathbf{G}^{\top }+\mathbf{U}+\mathbf{W},
\label{eq:regmodel}
\end{equation}%
where $\mathbf{U}\in \mathbb{R}^{n\times n}$ is a centered error matrix
given by $\mathbf{U}=\frac{1}{N}\sum_{i=1}^{N}\mathbf{U}_{i}$ and $\mathbf{W}%
=\frac{1}{N}\sum\limits_{i=1}^{N}\mathbf{W}_{i}$. {In the above regression
model \eqref{eq:regmodel}, there are two errors terms of a different nature.
The term $\mathbf{W}$ corresponds to the additive Gaussian errors $\mathcal{E%
}_{1},...,\mathcal{E}_{N}$ in model \eqref{eq:modnoise}, while the term $%
\mathbf{U}=\mathbf{S}-\mathbf{\Sigma }$ represents the difference between
the (unobserved) sample covariance matrix $\mathbf{S}$ and the matrix $%
\mathbf{\Sigma }$ that we want to estimate. }


This approach can be interpreted as a thresholding procedure of the entries
of an empirical matrix. To see this, consider the simple case where $M=n$
and the basis functions and observations points are chosen such that the
matrix $\mathbf{G}$ is orthogonal. Let $\mathbf{Y}=\mathbf{G}^{\top }%
\widetilde{\mathbf{S}}\mathbf{G}$ be a transformation of the empirical
covariance matrix $\widetilde{\mathbf{S}}$. In the orthogonal case, the
following proposition shows that the group Lasso estimator $\widehat{\mathbf{%
\Psi }}_{\lambda }$ defined by (\ref{MatrixGroupLassoEstimator}) consists in
thresholding the columns/rows of $\mathbf{Y}$ whose $\ell _{2}$-norm is too
small, and in multiplying the other columns/rows by weights between $0$ and $%
1$. Hence, the group Lasso estimate (\ref{MatrixGroupLassoEstimator}) can be
interpreted as covariance estimation by soft-thresholding  the
columns/rows of $\mathbf{Y}$.

\begin{proposition}
\label{prop:ortho} Suppose that $M=n$ and that $\mathbf{G}^{\top }\mathbf{G}=%
\mathbf{I}_{n}$ where $\mathbf{I}_{n}$ denotes the identity matrix of size $%
n\times n$. Let $\mathbf{Y}=\mathbf{G}^{\top }\widetilde{\mathbf{S}}\mathbf{G%
}$. Then, the group Lasso estimator $\widehat{\mathbf{\Psi }}_{\lambda }$
defined by (\ref{MatrixGroupLassoEstimator}) is the $n\times n$ symmetric
matrix whose entries are given by 
\begin{equation}
\left( \widehat{\mathbf{\mathbf{\Psi }}}_{\lambda }\right) _{mk}=\left\{ 
\begin{array}{ccc}
0 & \mbox{ if } & \sqrt{\sum_{j=1}^{M}\mathbf{Y}_{jk}^{2}}\leq \lambda
\gamma _{k}, \\ 
Y_{mk}\left( 1-\frac{\lambda \gamma _{k}}{\sqrt{\sum_{j=1}^{M}\mathbf{Y}%
_{mk}^{2}}}\right) & \mbox{ if } & \sqrt{\sum_{j=1}^{M}\mathbf{Y}_{jk}^{2}}%
>\lambda \gamma _{k},%
\end{array}%
\right.  \label{CGLEotho}
\end{equation}%
for $1\leq k,m\leq M$.
\end{proposition}

\section{Consistency of the group Lasso estimator}

\label{sec:oracle}

\subsection{Notations and main assumptions}

Let us begin by some definitions. For a symmetric $p\times p$ matrix $%
\mathbf{A}$ with real entries, $\rho _{min}(\mathbf{A})$ denotes the
smallest eigenvalue of $\mathbf{A}$, and $\rho _{max}(\mathbf{A})$ denotes
the largest eigenvalue of $\mathbf{A}$. For $\beta \in \mathbb{R}^{q}$, $%
\Vert \beta \Vert _{\ell _{2}}$ denotes the usual Euclidean norm of $\beta $%
. For $p\times q$ matrix $\mathbf{A}$ with real entries, $\Vert \mathbf{A}%
\Vert _{2}=\sup_{\beta \in \mathbb{R}^{q},\;\beta \neq 0}\frac{\Vert \mathbf{%
A}\beta \Vert _{\ell _{2}}}{\Vert \beta \Vert _{\ell _{2}}}$ denotes the
operator norm of $\mathbf{A}$. Recall that if $\mathbf{A}$ is a non negative
definite matrix with $p=q$ then $\Vert \mathbf{A}\Vert _{2}=\rho _{max}(%
\mathbf{A})$.


Let $\mathbf{\Psi} \in \mathcal{S}_{M}$ and $\beta$ a vector in $\mathbb{R}%
^M $. For a subset $J \subset \{1,\ldots,M\}$ of indices of cardinality $|J|$%
, then $\beta_{J}$ is the vector in $\mathbb{R}^{M}$ that has the same
coordinates as $\beta$ on $J$ and zeros coordinates on the complement $J^{c}$
of $J$. The $n \times |J|$ matrix obtained by removing the columns of $%
\mathbf{G}$ whose indices are not in $J$ is denoted by $\mathbf{G}_{J}$. The
sparsity of $\mathbf{\Psi}$ is defined as its number of non-zero columns
(and thus by symmetry non-zero rows) namely

\begin{definition}
For $\mathbf{\Psi} \in \mathcal{S}_{M}$, the sparsity of $\mathbf{\Psi}$ is 
\begin{equation*}
\mathcal{M}\left( \mathbf{\mathbf{\Psi}}\right) =\#\left\{ k:\mathbf{\Psi }%
_{k}\neq \mathbf{0}\right\}.
\end{equation*}
\end{definition}

Then, let us introduce the following quantities that control the minimal
eigenvalues of sub-matrices of small size extracted from the matrix $\mathbf{%
G}^{\top }\mathbf{G}$, and the correlations between the columns of $\mathbf{G%
}$:

\begin{definition}
\label{def:eig} Let $0<s\leq M$. Then, 
\begin{equation*}
\rho _{\min }(s):=\underset{%
\begin{array}{c}
J\subset \{1,\ldots ,M\} \\ 
|J|\leq s%
\end{array}%
}{\inf }\left( \frac{\beta _{J}^{\top }\mathbf{G}^{\top }\mathbf{G}\beta _{J}%
}{\Vert \beta _{J}\Vert _{\ell _{2}}^{2}}\right) =\underset{%
\begin{array}{c}
J\subset \{1,\ldots ,M\} \\ 
|J|\leq s%
\end{array}%
}{\inf }\rho _{\min }\left( \mathbf{G}_{J}^{\top }\mathbf{G}_{J}\right) .
\end{equation*}
\end{definition}

\begin{definition}
The mutual coherence $\theta (\mathbf{G})$ of the columns $\mathbf{G}_{k}$, $%
k=1,\ldots ,M$ of $\mathbf{G}$ is defined as 
\begin{equation*}
\theta (\mathbf{G}):=\max \left\{ \left\vert \mathbf{G}_{k^{\prime }}^{\top }%
\mathbf{G}_{k}\right\vert ,\;k\neq k^{\prime },\;1\leq k,k^{\prime }\leq
M\right\} ,
\end{equation*}%
and let 
\begin{equation*}
\mathbf{G}_{\max }^{2}:=\max \left\{ \Vert \mathbf{G}_{k}\Vert _{\ell
_{2}}^{2},\;1\leq k\leq M\right\} .
\end{equation*}
\end{definition}

To derive oracle inequalities showing the consistency of the group Lasso
estimator $\widehat{\mathbf{\Psi }}_{\lambda }$ the correlations between the
columns of $\mathbf{G}$ (measured by $\theta (\mathbf{G})$) should not be
too large when compared to the minimal eigenvalues of small matrices
extracted from $\mathbf{G}^{\top }\mathbf{G}$, which is formulated in the
following assumption:

\begin{assump}
\label{ass:eig} Let $c_{0}>0$ be some constant and $0<s\leq M$. Then 
\begin{equation*}
\theta (\mathbf{G})<\frac{\rho _{\min }(s)^{2}}{c_{0}\rho _{\max }(\mathbf{G}%
^{\top }\mathbf{G})s}.
\end{equation*}
\end{assump}

Assumption \ref{ass:eig} is inspired by recent results in \cite{BRT} on the
consistency of Lasso estimators in the standard nonparametric regression
model using a large dictionary of basis functions. In \cite{BRT}, a general
condition called \textit{restricted eigenvalue assumption} is introduced to
control the minimal eigenvalues of the Gram matrix associated to the
dictionary over sets of sparse vectors. In the setting of nonparametric
regression, a condition similar to Assumption \ref{ass:eig} is given in \cite%
{BRT} as an example for which the restricted eigenvalue assumption holds. 

Let us give some examples for which Assumption \ref{ass:eig} is satisfied.
If $M \leq n$ and the design points are chosen such that the columns of the
matrix $\mathbf{G}$ are orthonormal vectors in $\mathbb{R}^n$, then for any $%
0 < s \leq M$ one has that $\rho_{\min}(s) = 1$ and $\theta(\mathbf{G}) = 0$
and thus Assumption \ref{ass:eig} holds for any value of $c_{0}$ and $s$.

Now, suppose that the columns of $\mathbf{G}$ are normalized to one, i.e $%
\Vert \mathbf{G}_{k}\Vert _{\ell _{2}}=1$, $k=1,\ldots ,M$ implying that $%
\mathbf{G}_{\max }=1$. Let $\beta \in \mathbb{R}^{M}$. Then, for any $%
J\subset \{1,\ldots ,M\}$ with $\left\vert J\right\vert \leq s\leq \min
(n,M) $ 
\begin{equation*}
\beta _{J}^{\top }\mathbf{G}^{\top }\mathbf{G}\beta _{J}\geq \Vert \beta
_{J}\Vert _{\ell _{2}}^{2}-\theta (\mathbf{G})s\Vert \beta _{J}\Vert _{\ell
_{2}}^{2},
\end{equation*}%
which implies that 
\begin{equation*}
\rho _{\min }(s)\geq 1-\theta (\mathbf{G})s.
\end{equation*}%
Therefore, if $(1-\theta (\mathbf{G})(s-1))^{2}>c_{0}\theta (\mathbf{G})\rho
_{\max }(\mathbf{G}^{\top }\mathbf{G})s$, then Assumption \ref{ass:eig} is
satisfied.\newline

Let us now specify the law of the stochastic process $X$. For this, recall
that for a real-valued random variable $Z$, the $\psi _{\alpha }$ Orlicz
norm of $Z$ is 
\begin{equation*}
\Vert Z\Vert _{\psi _{\alpha }}:=\inf \left\{ C>0\;;\;\mathbb{E}\exp \left( 
\frac{|Z|^{\alpha }}{C^{\alpha }}\right) \leq 2\right\} .
\end{equation*}%
Such Orlicz norms are useful to characterize the tail behavior of random
variables. Indeed, if $\Vert Z\Vert _{\psi _{\alpha }}<+\infty $ then this
is equivalent to assuming that there exists two constants $%
K_{1},K_{2}>0 $ such that for all $x>0$ 
\begin{equation*}
\mathbb{P}\left( \left\vert Z\right\vert \geq x\right) \leq K_{1}\exp \left(
-\frac{x^{\alpha }}{K_{2}^{\alpha }}\right) ,
\end{equation*}%
(see e.g.\ \cite{MR2265341} for more
details on Orlicz norms of random variables) .
Therefore, if $\Vert Z\Vert _{\psi _{2}}<+\infty $ then $Z$ is said to have
a sub-Gaussian behavior and if $\Vert Z\Vert _{\psi _{1}}<+\infty $ then $Z$
is said to have a sub-Exponential behavior. In the next sections, oracle
inequalities for the group Lasso estimator will be derived under the
following assumption on $X$:

\begin{assump}
\label{ass:X} The random vector $\mathbf{X}=\left( X\left( t_{1}\right)
,...,X\left( t_{n}\right) \right) ^{\top }\in \mathbb{R}^{n}$ is such that

\begin{description}
\item[(A1)] There exists $\rho \left( \mathbf{\Sigma }\right) >0$ such that,
for all vector $\beta \in \mathbb{R}^{n}$ with $\Vert \beta \Vert _{\ell
_{2}}=1$, then $\left( \mathbb{E}|\mathbf{X}^{\top }\beta |^{4}\right)
^{1/4}<\rho \left( \mathbf{\Sigma }\right) $.

\item[(A2)] Set $Z=\Vert \mathbf{X}\Vert _{\ell _{2}}$. There exists $\alpha
\geq 1$ such that $\Vert Z\Vert _{\psi _{\alpha }}<+\infty $.
\end{description}
\end{assump}
Note that \textbf{(A1)} implies that $\Vert \mathbf{\Sigma }\Vert _{2}\leq
\rho \left( \mathbf{\Sigma }\right) ^{2}$. Indeed, one has that 
\begin{eqnarray*}
\Vert \mathbf{\Sigma }\Vert _{2}=\rho _{\max }(\mathbf{\Sigma })
&=&\sup_{\beta \in \mathbb{R}^{n},\;\Vert \beta \Vert _{\ell _{2}}=1}\beta
^{\top }\mathbf{\Sigma }\beta =\sup_{\beta \in \mathbb{R}^{n},\;\Vert \beta
\Vert _{\ell _{2}}=1}\mathbb{E}\beta ^{\top }\mathbf{X}\mathbf{X}^{\top
}\beta \\
&=&\sup_{\beta \in \mathbb{R}^{n},\;\Vert \beta \Vert _{\ell _{2}}=1}\mathbb{%
E}|\beta ^{\top }\mathbf{X}|^{2}\leq \sup_{\beta \in \mathbb{R}^{n},\;\Vert
\beta \Vert _{\ell _{2}}=1}\sqrt{\mathbb{E}|\beta ^{\top }\mathbf{X}|^{4}}%
\leq \rho ^{2}\left( \mathbf{\Sigma }\right) .
\end{eqnarray*}%
When $X$ is a Gaussian process, it follows that for any  $\beta \in \mathbb{R}^{n}$ with $\Vert \beta \Vert _{\ell
_{2}}=1$ then $ \left( \E | \bX^{\top} \beta |^{4} \right)^{1/4}= 3^{1/4} \left( \beta^{\top} \bSigma  \beta \right)^{1/2} $ since $ \bX^{\top} \beta \sim N(0, \beta^{\top} \bSigma  \beta) $. Therefore, under the assumption that $X$ is a Gaussian process, Assumption \textbf{(A1)} holds with $\rho \left( \mathbf{\Sigma }\right) = 3^{1/4} \Vert \mathbf{\Sigma }\Vert _{2}^{1/2}$.

Assumption \textbf{(A2)} requires that $\Vert Z\Vert _{\psi _{\alpha
}}<+\infty $, where $Z=\Vert \mathbf{X}\Vert _{\ell _{2}}$. The following
proposition provides some examples where such an assumption holds.

\begin{proposition}
\label{prop:orlicz} Let $Z=\Vert \mathbf{X}\Vert _{\ell _{2}}=\left(
\sum_{i=1}^{n}|X(t_{i})|^{2}\right) ^{1/2}$. Then

\begin{description}
\item[-] If $X$ is a Gaussian process 
\begin{equation*}
\| Z \|_{\psi_{2}} < \sqrt{8/3} \sqrt{tr( \mathbf{\Sigma})}.
\end{equation*}

\item[-] If the random process $X$ is such that $\| Z \|_{\psi_{2}} < +
\infty$, and there exists a constant $C_{1}$ such that $\| \mathbf{\Sigma}%
_{ii}^{-1/2} |X(t_{i})| \|_{\psi_{2}} \leq C_{1}$ for all $i =1,\ldots,n$,
then 
\begin{equation*}
\| Z \|_{\psi_{2}} < C_{1} \sqrt{tr( \mathbf{\Sigma})}.
\end{equation*}

\item[-] If $X$ is a bounded process, meaning that there exists a constant $%
R>0$ such that for all $t\in \mathbb{T}$, $|X(t)|\leq R$, then for any $%
\alpha \geq 1$, 
\begin{equation*}
\Vert Z\Vert _{\psi _{\alpha }}\leq \sqrt{n}R(\log 2)^{-1/\alpha }.
\end{equation*}
\end{description}
\end{proposition}

Assumption \ref{ass:X} will be used to control the deviation in operator
norm between the sample covariance matrix $\mathbf{S}$ and the true
covariance matrix $\mathbf{\Sigma }$ in the sense of the following
proposition whose proof follows from Theorem 2.1 in \cite{MR2265341}.

\begin{proposition}
\label{prop:mp} Let $X_{1},...,X_{N}$ be independent copies of the
stochastic process $X$, let $Z=\Vert \mathbf{X}\Vert _{\ell _{2}}$ and $%
\mathbf{X}_{i}=\left( X_{i}\left( t_{1}\right) ,...,X_{i}\left( t_{n}\right)
\right) ^{\top }$ for $i=1,\ldots ,N$. Recall that $\mathbf{S}=\frac{1}{N}%
\sum\limits_{i=1}^{N}\mathbf{X}_{i}\mathbf{X}_{i}^{\top }\mbox{ and }\mathbf{%
\Sigma }=\mathbb{E}\left( \mathbf{X}\mathbf{X}^{\top }\right) $. Suppose
that $X$ satisfies Assumption \ref{ass:X}. Let $d=\min (n,N)$. Then, there
exists a universal constant $\delta _{\ast }>0$ such that for all $x>0$ 
\begin{equation}
\mathbb{P}\left( \Big\|\mathbf{S}-\mathbf{\Sigma }\Big\|_{2}\geqslant \tau
_{d,N,n}x\right) \leqslant \exp \left( -(\delta _{\ast }^{-1}x)^{\frac{%
\alpha }{2+\alpha }}\right) ,  \label{eq:MP}
\end{equation}%
where $\tau_{N,n}=\max (A_{N,n}^{2},B_{N,n})$, with 
\begin{equation*}
A_{N,n}=\Vert Z\Vert _{\psi _{\alpha }}\frac{\sqrt{\log d}(\log N)^{1/\alpha
}}{\sqrt{N}}\mbox{ and }B_{N,n}=\frac{\rho ^{2}\left( \mathbf{\Sigma }%
\right) }{\sqrt{N}}+\Vert \mathbf{\Sigma }\Vert _{2}^{1/2}A_{N,n}.
\end{equation*}
\end{proposition}

Let us briefly comment Proposition \ref{prop:mp} in some specific cases. If $%
X$ is Gaussian, then Proposition \ref{prop:orlicz} implies that $A_{N,n}\leq
A_{N,n,1}$, where 
\begin{equation}
A_{N,n,1}=\sqrt{8/3}\sqrt{tr(\mathbf{\Sigma })}\frac{\sqrt{\log d}(\log
N)^{1/\alpha }}{\sqrt{N}}\leq \sqrt{8/3}\;\Vert \mathbf{\Sigma }\Vert
_{2}^{1/2}\sqrt{\frac{n}{N}}\sqrt{\log d}(\log N)^{1/\alpha },
\label{eq:boundA1}
\end{equation}%
and in this case inequality (\ref{eq:MP}) becomes 
\begin{equation}
\mathbb{P}\left( \Big\|\mathbf{S}-\mathbf{\Sigma }\Big\|_{2}\geqslant \max
\left( A_{N,n,1}^{2},B_{N,n,1}\right) x\right) \leqslant \exp \left(
-(\delta _{\ast }^{-1}x)^{\frac{\alpha }{2+\alpha }}\right)  \label{eq:MP1}
\end{equation}%
for all $x>0$, where $B_{N,n,1}=\frac{\rho ^{2}\left( \mathbf{\Sigma }%
\right) }{\sqrt{N}}+\Vert \mathbf{\Sigma }\Vert _{2}^{1/2}A_{N,n,1}$.

If $X$ is a bounded process by some constant $R>0$ , then using Proposition %
\ref{prop:orlicz} and by letting $\alpha \rightarrow +\infty $, Proposition %
\ref{prop:mp} implies that for all $x>0$, 
\begin{equation}
\mathbb{P}\left( \Big\|\mathbf{S}-\mathbf{\Sigma }\Big\|_{2}\geqslant \max
\left( A_{N,n,2}^{2},B_{N,n,2}\right) x\right) \leqslant \exp \left( -\delta
_{\ast }^{-1}x\right) ,  \label{eq:MP2}
\end{equation}%
where 
\begin{equation}
A_{N,n,2}=R\sqrt{\frac{n}{N}}\sqrt{\log d}\mbox{ and }B_{N,n,2}=\frac{\rho
^{2}\left( \mathbf{\Sigma }\right) }{\sqrt{N}}+\Vert \mathbf{\Sigma }\Vert
_{2}^{1/2}A_{N,n,2}.  \label{eq:boundA2}
\end{equation}

Contrary to the low-dimensional case ($n << N$), in a high-dimensional setting when $n>>N$ or when $n$ and $N$ are of the same
magnitude ($\frac{n}{N}\rightarrow c>0$ as $n,N\rightarrow +\infty $), inequalities (\ref{eq:MP1}) and (\ref{eq:MP2}) cannot be used to conclude that the norm $\Big\|\mathbf{S}-\mathbf{\Sigma }\Big\|_{2}$ concentrates around zero. Actually, it is well known that the sample covariance  $\mathbf{S}$ is a bad estimator of $\mathbf{\Sigma }$ in a high-dimensional setting, and that  without any further restriction on the structure of
the covariance matrix $\mathbf{\Sigma }$, then $\mathbf{S}$ cannot be a consistent estimator. However, we would like to point out that Proposition \ref{prop:mp} relates the quality of $\mathbf{S}$ to the ``true dimensionality'' of the vector $ \mathbf{X} =\left( X \left( t_{1}\right) ,...,X\left( t_{n}\right) \right) ^{\top} \in \R^{n}$ that is measured by the quantity $\Vert Z\Vert _{\psi _{\alpha }}$ with $ Z = \Vert \mathbf{X}\Vert _{\ell _{2}}$. Indeed, if $X$ is a low-dimensional Gaussian process such that $tr(\mathbf{\Sigma }) = 1$ then Proposition  \ref{prop:mp}  and  inequality \eqref{eq:boundA1} imply that
\begin{equation}
\mathbb{P}\left( \Big\|\mathbf{S}-\mathbf{\Sigma }\Big\|_{2}\geqslant \max
\left( A_{N}^{2},B_{N}\right) x\right) \leqslant \exp \left(
-(\delta _{\ast }^{-1}x)^{\frac{1 }{2}}\right)  \label{eq:MPlow}
\end{equation}
for all $x>0$, where $A_{N} = \sqrt{8/3} \frac{\sqrt{\log N}(\log
N)^{1/\alpha }}{\sqrt{N}}$ and  $B_{N}=\frac{\rho ^{2}\left( \mathbf{\Sigma }%
\right) }{\sqrt{N}}+\Vert \mathbf{\Sigma }\Vert _{2}^{1/2}A_{N}$. Hence, inequality \eqref{eq:MPlow} shows that, under an assumption of low-dimensionality of the process $X$,  the deviation in operator
norm between $\mathbf{S}$ and $\mathbf{\Sigma }$  depends on the ratio $\frac{1}{N}$ and not  on  $\frac{n}{N}$, and thus the quality of $\mathbf{S}$ as an estimator of $\mathbf{\Sigma }$ is much better in such settings.

More generally, another assumption of low-dimensionality for the process $X$ is to suppose that it has a sparse representation in a
dictionary of basis functions, which may also improve the quality of $\mathbf{S}$ as an
estimator of $\mathbf{\Sigma }$. To see this, consider the simplest case $%
X=X^{0}$, where the process $X^{0}$ has a sparse representation in the basis 
$(g_{m})_{1\leq m\leq M}$ given by 
\begin{equation}
X^{0}(t)=\sum_{m\in J^{\ast }}a_{m}g_{m}(t),\;t\in \mathbb{T},  \label{eq:X0}
\end{equation}%
where $J^{\ast }\subset \{1,\ldots ,M\}$ is a subset of indices of
cardinality $|J^{\ast }|=s_{\ast }$ and $a_{m}$, $m\in J^{\ast }$ are random
coefficients (possibly correlated). Under such an assumption, the following
proposition holds.

\begin{proposition}
\label{prop:MPsparse} Suppose that $X=X^{0}$ with $X^{0}$ defined by (\ref%
{eq:X0}) with $s_{\ast }\leq \min (n,M)$. Assume that $X$ satisfies
Assumption \ref{ass:X} and that the matrix $\mathbf{G}_{J^{\ast }}^{\top }%
\mathbf{G}_{J^{\ast }}$ is invertible, where $\mathbf{G}_{J^{\ast }}$
denotes the $n\times |J^{\ast }|$ matrix obtained by removing the columns of 
$\mathbf{G}$ whose indices are not in $J^{\ast }$. Then, there exists a
universal constant $\delta _{\ast }>0$ such that for all $x>0$, 
\begin{equation}
\mathbb{P}\left( \Big\|\mathbf{S}-\mathbf{\Sigma }\Big\|_{2}\geqslant \tilde{\tau}_{N,s_{\ast }} x\right) \leqslant \exp \left( -(\delta _{\ast
}^{-1}x)^{\frac{\alpha }{2+\alpha }}\right) ,  \label{eq:MPsparse}
\end{equation}%
where $\tilde{\tau}_{N,s_{\ast }}=\max (\tilde{A}_{N,s_{\ast }}^{2},\tilde{B}_{N,s_{\ast }})$, with 
\begin{equation*}
\tilde{A}_{N,s_{\ast }}=\rho _{\max }^{1/2}\left( \mathbf{G}%
_{J^{\ast }}^{\top }\mathbf{G}_{J^{\ast }}\right) \Vert \tilde{Z}\Vert
_{\psi _{\alpha }}\frac{\sqrt{\log d^{\ast }}(\log N)^{1/\alpha }}{\sqrt{N}},
\end{equation*}%
and%
\begin{equation*}
\tilde{B}_{N,s_{\ast }}=\left( \frac{\rho _{\max }\left( \mathbf{G}%
_{J^{\ast }}^{\top }\mathbf{G}_{J^{\ast }}\right) }{\rho _{\min }\left( 
\mathbf{G}_{J^{\ast }}^{\top }\mathbf{G}_{J^{\ast }}\right) }\right) \frac{%
\rho ^{2}\left( \mathbf{\Sigma }\right) }{\sqrt{N}}+\left( \frac{\rho _{\max
}\left( \mathbf{G}_{J^{\ast }}^{\top }\mathbf{G}_{J^{\ast }}\right) }{\rho
_{\min }\left( \mathbf{G}_{J^{\ast }}^{\top }\mathbf{G}_{J^{\ast }}\right) }%
\right) ^{1/2}\left\Vert \mathbf{\Sigma }\right\Vert _{2}^{1/2}\tilde{A}%
_{d^{\ast },N,s_{\ast }},
\end{equation*}%
with $d^{\ast }=\min (N,s_{\ast })$ and $\tilde{Z}=\Vert \mathbf{a}_{J^{\ast
}}\Vert _{\ell _{2}}$, where $\mathbf{a}_{J^{\ast }}=(\mathbf{G}_{J^{\ast
}}^{\top }\mathbf{G}_{J^{\ast }})^{-1}\mathbf{G}_{J^{\ast }}^{\top }\mathbf{X%
}\in \mathbb{R}^{s_{\ast }}$.
\end{proposition}

Using Proposition \ref{prop:orlicz} and Proposition \ref{prop:MPsparse} it
follows that

\begin{description}
\item[-] If $X=X^{0}$ is a Gaussian process then 
\begin{equation}
\tilde{A}_{N,s_{\ast }} \leq \sqrt{8/3}\;\left( \frac{\rho _{\max
}\left( \mathbf{G}_{J^{\ast }}^{\top }\mathbf{G}_{J^{\ast }}\right) }{\rho
_{\min }\left( \mathbf{G}_{J^{\ast }}^{\top }\mathbf{G}_{J^{\ast }}\right) }%
\right) ^{1/2}\Vert \mathbf{\Sigma }\Vert _{2}^{1/2}\sqrt{\frac{s_{\ast }}{N}%
}\sqrt{\log d^{\ast }}(\log N)^{1/\alpha }  \label{eq:boundA1sparse}
\end{equation}

\item[-] If $X=X^{0}$ is such that the random variables $a_{m}$ are bounded
by for some constant $R>0$, then 
\begin{equation}
\tilde{A}_{N,s_{\ast }} \leq R\Vert g\Vert _{\infty }\sqrt{\frac{%
s_{\ast }}{N}}\sqrt{\log d^{\ast }}  \label{eq:boundA2sparse}
\end{equation}%
with $\Vert g\Vert _{\infty }=\max_{1\leq m\leq M}\Vert g_{m}\Vert _{\infty
} $ where $\Vert g_{m}\Vert _{\infty }=\sup_{t\in \mathcal{T}}|g_{m}(t)|$.
\end{description}

Therefore, let us compare the bounds (\ref{eq:boundA1sparse}) and (\ref%
{eq:boundA2sparse}) with the inequalities (\ref{eq:boundA1}) and (\ref%
{eq:boundA2}). It follows that, in the case $X=X^{0}$, if the sparsity $%
s_{\ast }$ of $X$ in the dictionary is small compared to the number of time
points $n$ then the deviation between $\mathbf{S}$ and $\mathbf{\Sigma }$ is
much smaller than in the general case without any assumption on the
structure of $\mathbf{\Sigma }$. Obviously, the gain also depends on the
control of the ratio $\frac{\rho _{\max }\left( \mathbf{G}_{J^{\ast }}^{\top
}\mathbf{G}_{J^{\ast }}\right) }{\rho _{\min }\left( \mathbf{G}_{J^{\ast
}}^{\top }\mathbf{G}_{J^{\ast }}\right) }$. Note that in the case of an
orthonormal design ($M=n$ and $\mathbf{G}^{\top }\mathbf{G}=\mathbf{I}_{n}$)
then $\rho _{\max }\left( \mathbf{G}_{J^{\ast }}^{\top }\mathbf{G}_{J^{\ast
}}\right) =\rho _{\min }\left( \mathbf{G}_{J^{\ast }}^{\top }\mathbf{G}%
_{J^{\ast }}\right) =1$ for any $J^{\ast }$, and thus the gain in operator
norm between $\mathbf{S}$ and $\mathbf{\Sigma }$ clearly depends on the size
of $\frac{s_{\ast }}{N}$ compared to $\frac{n}{N}$. {Supposing that $X=X^{0}$
also implies that the operator norm of the error term $\mathbf{U}$ in the
matrix regression model \eqref{eq:regmodel} is controlled by the ratio $%
\frac{s_{\ast }}{N}$ instead of the ratio $\frac{n}{N}$ when no assumptions
are made on the structure of $\mathbf{\Sigma }$. This means that if $X$ has
a sparse representation in the dictionary then the error term $\mathbf{U}$
becomes smaller.}

\subsection{An oracle inequality for the Frobenius norm}

Consistency is first studied for the normalized Frobenius norm $\frac{1}{n}%
\left\Vert \mathbf{A}\right\Vert _{F}^{2}$ for an $n\times n$ matrix $%
\mathbf{A}$. The following theorem provides an oracle inequality for the
group Lasso estimator $\widehat{\mathbf{\Sigma }}_{\lambda }=\mathbf{G}%
\widehat{\mathbf{\Psi }}_{\lambda }\mathbf{G}^{\top }$.

\begin{theorem}
\label{theo:oracle1} Assume that $X$ satisfies Assumption \ref{ass:X}. Let $%
\epsilon >0$ and $1\leq s\leq \min (n,M)$. Suppose that Assumption \ref%
{ass:eig} holds with $c_{0}=3+4/\epsilon $, and that the covariance matrix  $\mathbf{\Sigma }_{noise}=\mathbb{E}\left( \mathbf{W}_{1}\right) $ of the noise is positive-definite. Consider the group Lasso
estimator $\widehat{\mathbf{\Sigma }}_{\lambda }$ defined by (\ref%
{GroupLassoEstimator}) with the choices 
\begin{equation*}
\gamma _{k}=2\Vert \mathbf{G}_{k}\Vert _{\ell _{2}}\sqrt{\rho _{\max }(%
\mathbf{G}\mathbf{G}^{\top })},
\end{equation*}%
and 
\begin{equation*}
\lambda =\Vert \mathbf{\Sigma }_{noise}\Vert _{2}\left( 1+\sqrt{\frac{n}{N}}+%
\sqrt{\frac{2\delta \log M}{N}}\right) ^{2}\mbox{ for some constant }\delta
>1.
\end{equation*}%
Then, with probability at least $1-M^{1-\delta }$ one has that 
\begin{eqnarray}
\frac{1}{n}\left\Vert \widehat{\mathbf{\Sigma }}_{\lambda }-\mathbf{\Sigma }%
\right\Vert _{F}^{2} &\leq &(1+\epsilon )\underset{%
\begin{array}{l}
\mathbf{\mathbf{\Psi }}\in \mathcal{S}_{M} \\ 
\mathcal{M}\left( \mathbf{\mathbf{\Psi }}\right) \leq s%
\end{array}%
}{\inf }\left( \frac{4}{n}\left\Vert \mathbf{G\Psi G}^{\top }-\mathbf{\Sigma 
}\right\Vert _{F}^{2}\right. +\frac{8}{n}\left\Vert \mathbf{S}-\mathbf{%
\Sigma }\right\Vert _{F}^{2}  \label{eq:oracle} \\
&&\left. +C(\epsilon )\frac{\mathbf{G}_{\max }^{2}\rho _{\max }(\mathbf{G}%
^{\top }\mathbf{G})}{\kappa _{s,c_{0}}^{2}} 
\Vert \mathbf{\Sigma }_{noise}\Vert _{2}^{2}\left( 1+\sqrt{\frac{n}{N}}+%
\sqrt{\frac{2\delta \log M}{N}}\right) ^{4}
\frac{\mathcal{M}(%
\mathbf{\Psi })}{n}\right) ,  \notag
\end{eqnarray}%
where $\kappa _{s,c_{0}}^{2}=\rho _{\min }(s)^{2}-c_{0}\theta (\mathbf{G}%
)\rho _{\max }(\mathbf{G}^{\top }\mathbf{G})s,$ and $C(\epsilon )=8\frac{%
\epsilon }{1+\epsilon }(1+2/\epsilon )^{2}$.
\end{theorem}

The first term $\frac{1}{n}\left\Vert \mathbf{G\Psi G}^{\top }-\mathbf{%
\Sigma }\right\Vert _{F}^{2}$ in inequality (\ref{eq:oracle}) is the bias of
the estimator $\widehat{\mathbf{\Sigma }}_{\lambda }$. It reflects the
quality of the approximation of $\mathbf{\Sigma }$ by the set of matrices of
the form $\mathbf{G\Psi G}^{\top }$, with $\mathbf{\Psi }\in \mathcal{S}_{M}$
and $\mathcal{M}\left( \mathbf{\Psi }\right) \leq s$. As an example, suppose
that $X=X^{0}$, where the process $X^{0}$ has a sparse representation in the
basis $(g_{m})_{1\leq m\leq M}$ given by 
\begin{equation*}
X^{0}(t)=\sum_{m\in J^{\ast }}a_{m}g_{m}(t),\;t\in \mathbb{T},
\end{equation*}%
where $J^{\ast }\subset \{1,\ldots ,M\}$ is a subset of indices of
cardinality $|J^{\ast }|=s_{\ast }\leq s$ and $a_{m},m\in J^{\ast }$ are
random coefficients. Then, in this case, since $s_{\ast }\leq s$ the bias
term in (\ref{eq:oracle}) is equal to zero.

The second term $\frac{1}{n}\left\Vert \mathbf{S}-\mathbf{\Sigma }%
\right\Vert _{F}^{2}$ in (\ref{eq:oracle}) is a variance term as the
empirical covariance matrix $\mathbf{S}$ is an unbiased estimator of $%
\mathbf{\Sigma }$. Using the inequality $\frac{1}{n}\left\Vert A\right\Vert
_{F}^{2}\leq \left\Vert A\right\Vert _{2}^{2}$ that holds for any $n\times n$
matrix $A$, it follows that $\frac{1}{n}\left\Vert \mathbf{S}-\mathbf{\Sigma 
}\right\Vert _{F}^{2}\leq \left\Vert \mathbf{S}-\mathbf{\Sigma }\right\Vert
_{2}^{2}$. Therefore, under the assumption that $X$ has a sparse
representation in the dictionary (e.g. when $X=X_{0}$ as above) then the
variance term $\frac{1}{n}\left\Vert \mathbf{S}-\mathbf{\Sigma }\right\Vert
_{F}^{2}$ is controlled by the ratio $\frac{s_{\ast }}{N}\leq \frac{s}{N}$
(see Proposition \ref{prop:MPsparse}) instead of the ratio $\frac{n}{N}$
without any assumption on the structure of $\mathbf{\Sigma }$.

The third term in (\ref{eq:oracle}) is also a variance term due to the noise
in the measurements (\ref{eq:modnoise}). If there exists a constant $c>0$
independent of $n$ and $N$ such that $\frac{n}{N}\leq c$ then the decay of
this third variance term is essentially controlled by the ratio $\frac{%
\mathcal{M}(\mathbf{\Psi })}{n}\leq \frac{s}{n}$. Therefore, if $\mathcal{M}%
\left( \mathbf{\Psi }\right) \leq s$ with sparsity $s$ much smaller than $n$
then the variance of the group Lasso estimator $\widehat{\mathbf{\Sigma }}%
_{\lambda }$ is smaller than the variance of $\widetilde{\mathbf{S}}$. This
shows some of the improvements achieved by regularization (\ref%
{MatrixGroupLassoEstimator}) of the empirical covariance matrix $\widetilde{%
\mathbf{S}}$ with a group Lasso penalty.

An important assumption of Theorem \ref{theo:oracle1} is that the covariance matrix of the noise $\mathbf{\Sigma }_{noise}=\mathbb{E}\left( \mathbf{W}_{1}\right) $ is positive definite. This restriction is clearly necessary as illustrated by the following example: suppose that the contaminating process $\mathcal{E}\left( t \right) = \zeta g_{1}(t)$ with $\zeta \sim N(0,\sigma_{1}^{2})$, implying that $\mathbf{\Sigma }_{noise} =\sigma_{1}^{2} \mathbf{g}_{1} \mathbf{g}_{1}^{\top}$ with $ \mathbf{g}_{1} = ( g_{1}(t_{1}),\ldots,g_{1}(t_{n}) )^{\top} $  has $n-1$ eigenvalues equal to zero. Now, suppose that $X(t) = a_{2} g_{2}(t)$ with $a_{2}  \sim N(0,\sigma_{2}^{2})$. If $\sigma_{1} > \sigma_{2}$ then the group LASSO regularization alone cannot get rid of the additive error term without eliminating first the right component $g_{2}$. Hence, in such settings,   group LASSO regularization does not yield to a consistent estimation of $\bSigma = \sigma_{2}^{2} \mathbf{g}_{2} \mathbf{g}_{2}^{\top}$ with $ \mathbf{g}_{2} = ( g_{2}(t_{1}),\ldots,g_{2}(t_{n}) )^{\top} $.

\subsection{An oracle inequality for the operator norm}

The \textquotedblleft normalized\textquotedblright\ Frobenius norm $\frac{1}{%
n}\left\Vert \widehat{\mathbf{\Sigma }}_{\lambda }-\mathbf{\Sigma }%
\right\Vert _{F}^{2}$, i.e the average of the eigenvalues of $\left( \widehat{%
\mathbf{\Sigma }}_{\lambda }-\mathbf{\Sigma }\right) ^{2}$, can be viewed as
a reasonable proxy for the operator norm $\left\Vert \widehat{\mathbf{\Sigma 
}}_{\lambda }-\mathbf{\Sigma }\right\Vert _{2}^{2}$ (maximum eigenvalue of $%
\left( \widehat{\mathbf{\Sigma }}_{\lambda }-\mathbf{\Sigma }\right) ^{2}$).
It is thus expected that the results of Theorem \ref{theo:oracle1} imply
that the group Lasso estimator $\widehat{\mathbf{\Sigma }}_{\lambda }$ is a
good estimator of $\mathbf{\Sigma }$ in operator norm. Let us recall that
controlling the operator norm enables to study the convergence of the
eigenvectors and eigenvalues of $\widehat{\mathbf{\Sigma }}_{\lambda }$ by
controlling of the angles between the eigenspaces of a population and a
sample covariance matrix through the use of the $\sin (\theta )$ theorems in
\ \cite{MR0264450}.

Now, let us consider the case where $X$ consists in noisy observations of
the process $X^{0}$ (\ref{eq:X0}) meaning that 
\begin{equation}
\widetilde{X}(t_{j})=X^{0}(t_{j})+\mathcal{E}\left( t_{j}\right)
,\;j=1,\ldots ,n,  \label{eq:sparseX}
\end{equation}%
where $\mathcal{E}$ is a second order Gaussian process $\mathcal{E}$ with
zero mean and independent of $X^{0}$. In this case, one has that 
\begin{equation*}
\mathbf{\Sigma }=\mathbf{G}\mathbf{\Psi }^{\ast }\mathbf{G}^{\top },%
\mbox{
where }\mathbf{\Psi }^{\ast }=\mathbb{E}\left( \mathbf{a}\mathbf{a}^{\top
}\right) ,
\end{equation*}%
where $\mathbf{a}$ is the random vector of $\mathbb{R}^{M}$ with $\mathbf{a}%
_{m}=a_{m}$ for $m\in J^{\ast }$ and $\mathbf{a}_{m}=0$ for $m\notin J^{\ast
}$. Therefore, using Theorem \ref{theo:oracle1} by replacing $s$ by $s^\ast=\left\vert J^{\ast
}\right\vert $, since $\mathbf{\mathbf{\Psi }}^{\ast }\in \left\{ 
\mathbf{\mathbf{\Psi }}\in \mathcal{S}_{M}:M\left( \mathbf{\mathbf{\Psi }}%
\right) \leq s_{\ast }\right\} $, one can derive the following corrollary:

\begin{corollary}
\label{cor:oracle1} Suppose that the observations $\widetilde{X}_{i}(t_{j})$
with $i=1,...,N$ and $j=1,\ldots ,n$ are i.i.d random variables from model (%
\ref{eq:sparseX}) and that the conditions of Theorem \ref{theo:oracle1} are
satisfied with $1\leq s=s_{\ast }\leq \min (n,M)$. Then, with probability at
least $1-M^{1-\delta }$ one has that 
\begin{equation}
\frac{1}{n}\left\Vert \widehat{\mathbf{\Sigma }}_{\lambda }-\mathbf{\Sigma }%
\right\Vert _{F}^{2}\leq C_{0}\left( n,M,N,s_{\ast },\mathbf{S,}\mathbf{\Psi 
}^{\ast },\mathbf{G},\mathbf{\Sigma }_{noise}\right) ,  \label{eq:oracle2}
\end{equation}%
where 
\begin{equation*}
C_{0}\left( n,M,N,s_{\ast },\mathbf{S,}\mathbf{\Psi }^{\ast },\mathbf{G},%
\mathbf{\Sigma }_{noise}\right) =(1+\epsilon )\left( \frac{8}{n}\left\Vert 
\mathbf{S}-\mathbf{G}\mathbf{\Psi }^{\ast }\mathbf{G}^{\top }\right\Vert
_{F}^{2}+C(\epsilon )\frac{\mathbf{G}_{\max }^{2}\rho _{\max }(\mathbf{G}%
^{\top }\mathbf{G})}{\kappa _{s_{\ast },c_{0}}^{2}}\lambda ^{2}\frac{s_{\ast
}}{n}\right) .
\end{equation*}
\end{corollary}

To simplify notations, write $\widehat{\mathbf{\mathbf{\Psi }}}=\widehat{%
\mathbf{\mathbf{\Psi }}}_{\lambda }$, with $\widehat{\mathbf{\mathbf{\Psi }}}%
_{\lambda }$ given by (\ref{MatrixGroupLassoEstimator}). Define $\hat{J}%
_{\lambda }\subset \{1,\ldots ,M\}$ as%
\begin{equation}
\hat{J}_{\lambda }\equiv \hat{J}:=\left\{ k:\frac{\delta _{k}}{\sqrt{n}}%
\left\Vert \widehat{\mathbf{\mathbf{\Psi }}}_{k}\right\Vert _{\ell
_{2}}>C_{1}\left( n,M,N,s_{\ast },\mathbf{S,}\mathbf{\Psi }^{\ast },\mathbf{G%
},\mathbf{\Sigma }_{noise}\right) \right\} ,\mbox{ with }\delta _{k}=\frac{%
\Vert \mathbf{G}_{k}\Vert _{\ell _{2}}}{\mathbf{G}_{\max }},  \label{eq:hatJ}
\end{equation}%
and $C_{1}\left( n,M,N,s_{\ast },\mathbf{S,}\mathbf{\Psi }^{\ast },\mathbf{G},%
\mathbf{\Sigma }_{noise}\right) = C_{1}$ with
\begin{equation}
C_{1} = \max \left( \gamma _{\max }^{-1} n^{-1/2}  \frac{1+\epsilon}{\lambda}  \left\Vert \mathbf{S}-\mathbf{G\mathbf{\Psi }^{\ast } G}^{\top }\right\Vert _{F}^{2} ; \frac{4\left( 1+\epsilon \right) \sqrt{%
s_{\ast }}}{\epsilon \kappa _{s_{\ast },c_{0}}}\sqrt{C_{0}\left(
n,M,N,s_{\ast },\mathbf{S,}\mathbf{\Psi }^{\ast },\mathbf{G},\mathbf{\Sigma }%
_{noise}\right) } \right). \label{eq:C1}
\end{equation}%
with $\gamma _{\max }=2\mathbf{G}_{\max }\sqrt{\rho _{\max }(\mathbf{G}^{\top }\mathbf{G})}$.
The set of indices $\hat{J}$ is an estimation of the set of active basis
functions $J^{\ast }$. {Note that such thresholding procedure \eqref{eq:hatJ}
does not lead immediately to a practical way to choose the set $\hat{J}$.
Indeed the constant $C_{1}$ in \eqref{eq:hatJ} depends on the a priori
unknown sparsity $s_{\ast }$ and on the amplitude of the noise in the matrix
regression model \eqref{eq:regmodel} measured by the quantities $\frac{8}{n}%
\left\Vert \mathbf{S}-\mathbf{G}\mathbf{\Psi }^{\ast }\mathbf{G}^{\top
}\right\Vert _{F}^{2}$ and $\Vert \mathbf{\Sigma }_{noise}\Vert _{2}^{2}$.
Nevertheless, in Section \ref{sec:simus} on numerical experiments we give a
simple procedure to automatically threshold the $\ell _{2}$-norm of the
columns of the matrix $\widehat{\Psi }_{\lambda }$ that are two small. }

{Note that to estimate $J^{\ast }$ we did not simply take $\hat{J}=\hat{J}%
_{0}:=\left\{ k:\left\Vert \widehat{\mathbf{\mathbf{\Psi }}}_{k}\right\Vert
_{\ell _{2}}\neq 0\right\} $, but rather apply a thresholding step to
discard the columns of $\widehat{\mathbf{\Psi }}$ whose $\ell _{2}$-norm are
too small. By doing so, we want to stress the fact that to obtain a
consistent procedure with respect to the operator norm it is not sufficient to simply take $%
\hat{J}=\hat{J}_{0}$. A similar thresholding step is proposed in \cite%
{MR2386087} and \cite{LPTG} in the standard linear model to select a sparse
set of active variables when using regularization by a Lasso or group-Lasso
penalty. In the paper (\cite{MR2386087}), the second thresholding step used
to estimate the true sparsity pattern depends on a unknown constant that is
related to the amplitude of the unknown coefficients to estimate. }

Then, the
following theorem holds.

\begin{theorem}
\label{theo:activeset} Under the assumptions of Corollary \ref{cor:oracle1},
for any solution of problem (\ref{MatrixGroupLassoEstimator}), we have that
with probability at least $1-M^{1-\delta }$, 
\begin{equation}
\underset{1\leq k\leq M}{\max }\frac{\delta _{k}}{\sqrt{n}}\left\Vert 
\widehat{\mathbf{\mathbf{\Psi }}}_{k}-\mathbf{\mathbf{\Psi }}_{k}^{\ast
}\right\Vert _{\ell _{2}}\leq C_{1}\left( n,M,N,s_{\ast },\mathbf{S,}\mathbf{%
\Psi }^{\ast },\mathbf{G},\mathbf{\Sigma }_{noise}\right) .  \label{Ineq2}
\end{equation}%
If in addition 
\begin{equation}
\underset{k\in J^{\ast }}{\min }\frac{\delta _{k}}{\sqrt{n}}\left\Vert 
\mathbf{\mathbf{\Psi }}_{k}^{\ast }\right\Vert _{\ell _{2}}>2C_{1}\left(
n,M,N,s_{\ast },\mathbf{S,}\mathbf{\Psi }^{\ast },\mathbf{G},\mathbf{\Sigma }%
_{noise}\right)  \label{Hyp1}
\end{equation}%
then with the same probability the set of indices $\hat{J}$, defined by (\ref%
{eq:hatJ}), estimates correctly the true set of active basis functions $%
J^{\ast }$, that is $\hat{J}=J^{\ast }$ with probability at least $%
1-M^{1-\delta }$.
\end{theorem}

{The results of Theorem \ref{theo:activeset} indicate that if the $\ell _{2}$%
-norm of the columns of $\mathbf{\Psi }_{k}^{\ast }$ for $k\in J^{\ast }$
are sufficiently large with respect to the level of noise in the matrix
regression model \eqref{eq:regmodel} and the sparsity $s_{\ast }$, then $%
\hat{J}$ is a consistent estimation of the active set of variables. Indeed,
if $\mathcal{M}\left( \mathbf{\Psi ^{\ast }}\right) =s_{\ast }$, then by
symmetry the columns of $\mathbf{\Psi ^{\ast }}$ such $\mathbf{\mathbf{\Psi }%
}_{k}^{\ast }\neq 0$ have exactly $s_{\ast }$ non-zero entries. Hence, the
condition \eqref{Hyp1} means that the $\ell _{2}$-norm of $\mathbf{\mathbf{%
\Psi }}_{k}^{\ast }\neq 0$ (normalized by $\frac{\delta _{k}}{\sqrt{n}}$)
has to be larger than $\frac{4\left( 1+\epsilon \right) }{\epsilon \kappa
_{s_{\ast },c_{0}}}\sqrt{s_{\ast }}\sqrt{C_{0}}$. A simple condition to
satisfy such an assumption is that the amplitude of the $s_{\ast }$
non-vanishing entries of $\mathbf{\mathbf{\Psi }}_{k}^{\ast }\neq 0$ are
larger than $\frac{\sqrt{n}}{\delta _{k}}\frac{4\left( 1+\epsilon \right) }{%
\epsilon \kappa _{s_{\ast },c_{0}}}\sqrt{C_{0}}$ which can be interpreted as
a kind of measure of the noise in model \eqref{eq:regmodel}. } This suggests
to take as a final estimator of $\mathbf{\Sigma }$ the following matrix: 
\begin{equation}
\widehat{\mathbf{\Sigma }}_{\hat{J}}=\mathbf{G}_{\hat{J}}\widehat{\mathbf{%
\Psi }}_{\hat{J}}\mathbf{G}_{\hat{J}}  \label{SigmaHatJhat}
\end{equation}%
where $\mathbf{G}_{\hat{J}}$ denotes the $n\times |\hat{J}|$ matrix obtained
by removing the columns of $\mathbf{G}$ whose indices are not in $\hat{J}$,
and 
\begin{equation*}
\widehat{\mathbf{\Psi }}_{\hat{J}}=\argmin_{\mathbf{\Psi }\in \mathcal{S}_{|%
\hat{J}|}}\left\{ \left\Vert \widetilde{\mathbf{S}}-\mathbf{G}_{\hat{J}}%
\mathbf{\Psi }\mathbf{G}_{\hat{J}}^{\top }\right\Vert _{F}^{2}\right\} ,
\end{equation*}%
where $\mathcal{S}_{|\hat{J}|}$ denotes the set of $|\hat{J}|\times |\hat{J}%
| $ symmetric matrices. Note that if $\mathbf{G}_{\hat{J}}^{\top }\mathbf{G}%
_{\hat{J}}$ is invertible, then 
\begin{equation*}
\widehat{\mathbf{\Psi }}_{\hat{J}}=\left( \mathbf{G}_{\hat{J}}^{\top }%
\mathbf{G}_{\hat{J}}\right) ^{-1}\mathbf{G}_{\hat{J}}^{\top }\widetilde{%
\mathbf{S}}\mathbf{G}_{\hat{J}}\left( \mathbf{G}_{\hat{J}}^{\top }\mathbf{G}%
_{\hat{J}}\right) ^{-1}.
\end{equation*}%
Let us recall that if the observations are i.i.d random variables from model
(\ref{eq:sparseX}) then 
\begin{equation*}
\mathbf{\Sigma }=\mathbf{G}\mathbf{\Psi }^{\ast }\mathbf{G}^{\top },
\end{equation*}%
where $\mathbf{\Psi }^{\ast }=\mathbb{E}\left( \mathbf{a}\mathbf{a}^{\top
}\right) $, and $\mathbf{a}$ is the random vector of $\mathbb{R}^{M}$ with $%
\mathbf{a}_{m}=a_{m}$ for $m\in J^{\ast }$ and $\mathbf{a}_{m}=0$ for $%
m\notin J^{\ast }$. Then, define the random vector $\mathbf{a}_{J^{\ast
}}\in \mathbb{R}^{J^{\ast }}$ whose coordinates are the random coefficients $%
a_{m}$ for $m\in J^{\ast }$. Let $\mathbf{\Psi }_{J^{\ast }}=\mathbb{E}%
\left( \mathbf{a}_{J^{\ast }}\mathbf{a}_{J^{\ast }}^{\top }\right) $ and
denote by $\mathbf{G}_{J^{\ast }}$ the $n\times |J^{\ast }|$ matrix obtained
by removing the columns of $\mathbf{G}$ whose indices are not in $J^{\ast }$%
. Note that $\mathbf{\Sigma }=\mathbf{G}_{J^{\ast }}\mathbf{\Psi }_{J^{\ast
}}\mathbf{G}_{J^{\ast }}^{\top }$.

Assuming that $\mathbf{G}_{J^{\ast }}^{\top }\mathbf{G}_{J^{\ast }}$ is
invertible, define the matrix 
\begin{equation}
\mathbf{\Sigma }_{J^{\ast }}= \mathbf{\Sigma } + \mathbf{G}_{J^{\ast }} (%
\mathbf{G}_{J^{\ast }}^{\top }\mathbf{G}_{J^{\ast}})^{-1} \mathbf{G}%
_{J^{\ast }}^{\top } \mathbf{\Sigma }_{noise} \mathbf{G}_{J^{\ast }}\left( 
\mathbf{G}_{J^{\ast }}^{\top }\mathbf{G}_{J^{\ast }}\right) ^{-1} \mathbf{G}%
_{J^{\ast }}^{\top }.  \label{eq:SigmaJast}
\end{equation}%
Then, the following theorem gives a control of deviation between $\widehat{%
\mathbf{\Sigma }}_{\hat{J}}$ and $\mathbf{\Sigma }_{J^{\ast }}$ in operator
norm.

\begin{theorem}
\label{theo:opnorm} Suppose that the observations are i.i.d random variables
from model (\ref{eq:sparseX}) and that the conditions of Theorem \ref%
{theo:oracle1} are satisfied with $1\leq s=s_{\ast }\leq \min (n,M)$.
Suppose that $\mathbf{G}_{J^{\ast }}^{\top }\mathbf{G}_{J^{\ast }}$ is an
invertible matrix, and that 
\begin{equation*}
\underset{k\in J^{\ast }}{\min }\frac{\delta _{k}}{\sqrt{n}}\left\Vert 
\mathbf{\mathbf{\Psi }}_{k}^{\ast }\right\Vert _{\ell _{2}}>2C_{1}\left(
n,M,N,s_{\ast },\mathbf{S,}\mathbf{\Psi }^{\ast },\mathbf{G},\mathbf{\Sigma }%
_{noise}\right) ,
\end{equation*}%
where $C_{1}\left( n,M,N,s_{\ast },\mathbf{S,}\mathbf{\Psi }^{\ast },\mathbf{%
G},\mathbf{\Sigma }_{noise}\right) $ is the constant defined in (\ref{eq:C1}%
). Let $\mathbf{Y}=\left( \mathbf{G}_{J^{\ast }}^{\top }\mathbf{G}_{J^{\ast
}}\right) ^{-1}\mathbf{G}_{J^{\ast }}^{\top }\widetilde{\mathbf{X}}$ and $%
\tilde{Z}=\Vert \mathbf{Y}\Vert _{\ell _{2}}$ . Let $\rho \left( \mathbf{%
\Sigma }_{noise}\right) =\left( \sup_{\beta \in \mathbb{R}^{n},\Vert \beta
\Vert _{\ell _{2}}=1}\mathbb{E}|\mathcal{E}^{\top }\beta |^{4}\right) ^{1/4}$
where $\mathbf{\mathcal{E}}=\left( \mathcal{E}\left( t_{1}\right) ,...,%
\mathcal{E}\left( t_{n}\right) \right) ^{\top }$. Then, with probability at
least $1-M^{1-\delta }-M^{-\left( \frac{\delta _{\star }}{\delta _{\ast }}%
\right) ^{\frac{\alpha }{2+\alpha }}}$, with $\delta >1$ and $\delta _{\star
}>\delta _{\ast }$ one has that 
\begin{equation}
\left\Vert \widehat{\mathbf{\Sigma }}_{\hat{J}}-\mathbf{\Sigma }_{J^{\ast
}}\right\Vert _{2}\leq \rho _{\max }\left( \mathbf{G}_{J^{\ast }}^{\top }%
\mathbf{G}_{J^{\ast }}\right) \tilde{\tau}_{N,s_{\ast }}\delta _{\star
}\left( \log (M)\right) ^{\frac{2+\alpha }{\alpha }},  \label{eq:opnorm}
\end{equation}%
where $\tilde{\tau}_{N,s_{\ast }}=\max (\tilde{A}_{N,s_{\ast }}^{2},\tilde{B}%
_{N,s_{\ast }})$, with $\tilde{A}_{N,s_{\ast }}=\Vert \tilde{Z}\Vert _{\psi
_{\alpha }}\frac{\sqrt{\log d^{\ast }}(\log N)^{1/\alpha }}{\sqrt{N}}$, $%
\tilde{B}_{N,s_{\ast }}=\frac{\tilde{\rho}^{2}(\mathbf{\Sigma },\mathbf{%
\Sigma }_{noise})\rho _{\min }^{-1}\left( \mathbf{G}_{J^{\ast }}^{\top }%
\mathbf{G}_{J^{\ast }}\right) }{\sqrt{N}}+\left( \left\Vert \mathbf{\Psi }%
_{J^{\ast }}\right\Vert _{2}+\rho _{\min }^{-1}\left( \mathbf{G}_{J^{\ast
}}^{\top }\mathbf{G}_{J^{\ast }}\right) \left\Vert \mathbf{\Sigma }%
_{noise}\right\Vert _{2}\right) ^{1/2}\tilde{A}_{N,s_{\ast }}$, where $%
d^{\ast }=\min (N,s_{\ast })$ and $\tilde{\rho}(\mathbf{\Sigma },\mathbf{%
\Sigma }_{noise})=8^{1/4}\left( \rho ^{4}\left( \mathbf{\Sigma }\right)
+\rho ^{4}\left( \mathbf{\Sigma }_{noise}\right) \right) ^{1/4}$.
\end{theorem}

First note that the above theorem gives a deviation in operator norm from $%
\widehat{\mathbf{\Sigma }}_{\hat{J}}$ to the matrix $\mathbf{\Sigma }%
_{J^{\ast }}$ \eqref{eq:SigmaJast} which is not equal to the true covariance 
$\mathbf{\Sigma }$ of $X$ at the design points. Indeed, even if we know the
true sparsity set $J^{\ast }$, the additive noise in the measurements in
model \eqref{eq:modnoise} complicates the estimation of $\mathbf{\Sigma }$
in operator norm. However, although $\mathbf{\Sigma }_{J^{\ast }}\neq 
\mathbf{\Sigma }$, they can have the same eigenvectors if the structure of
the additive noise matrix term in \eqref{eq:SigmaJast} is not too complex.
As an example, consider the case of an additive white noise, for which $%
\mathbf{\Sigma }_{noise}=\sigma ^{2}\mathbf{I}_{n}$ where $\sigma $ is the
level of noise and $\mathbf{I}_{n}$ the $n\times n$ identity matrix. Under
such an assumption, if we further suppose for simplicity that $(\mathbf{G}%
_{J^{\ast }}^{\top }\mathbf{G}_{J^{\ast }})^{-1}=\mathbf{I}_{s_{\ast }}$,
then $\mathbf{\Sigma }_{J^{\ast }}=\mathbf{\Sigma }+\sigma ^{2}\mathbf{G}%
_{J^{\ast }}(\mathbf{G}_{J^{\ast }}^{\top }\mathbf{G}_{J^{\ast }})^{-1}%
\mathbf{G}_{J^{\ast }}^{\top }=\mathbf{\Sigma }+\sigma ^{2}\mathbf{I}_{n}$
and clearly $\mathbf{\Sigma }_{J^{\ast }}$ and $\mathbf{\Sigma }$ have the
same eigenvectors. Therefore, the eigenvectors of $\widehat{\mathbf{\Sigma }}%
_{\hat{J}}$ can be used as estimators of the eigenvectors of $\mathbf{\Sigma 
}$ which is suitable for the sparse PCA application described in the next
section on numerical experiments.

Let us illustrate the implications of Theorem \ref{theo:opnorm} on a simple
example. If $X$ is Gaussian, the random vector $\mathbf{Y}=\left( \mathbf{G}%
_{J^{\ast }}^{\top }\mathbf{G}_{J^{\ast }}\right) ^{-1}\mathbf{G}_{J^{\ast
}}^{\top }\left( \mathbf{X}+\mathbf{\mathcal{E}}\right) $ is also Gaussian
and Proposition \ref{prop:orlicz} can be used to prove that 
\begin{eqnarray*}
\Vert \tilde{Z}\Vert _{\psi _{2}} &\leq &\sqrt{8/3}\sqrt{tr\left( \left( 
\mathbf{G}_{J^{\ast }}^{\top }\mathbf{G}_{J^{\ast }}\right) ^{-1}\mathbf{G}%
_{J^{\ast }}^{\top }\left( \mathbf{\Sigma }+\mathbf{\Sigma }_{noise}\right) 
\mathbf{G}_{J^{\ast }}\left( \mathbf{G}_{J^{\ast }}^{\top }\mathbf{G}%
_{J^{\ast }}\right) ^{-1}\right) } \\
&\leq &\sqrt{8/3}\Vert \mathbf{\Sigma +\Sigma }_{noise}\Vert _{2}^{1/2}\rho
_{\min }^{-1/2}\left( \mathbf{G}_{J^{\ast }}^{\top }\mathbf{G}_{J^{\ast
}}\right) \sqrt{s_{\ast }}.
\end{eqnarray*}%
Then Theorem \ref{theo:opnorm} implies that with high probability 
\begin{equation*}
\left\Vert \widehat{\mathbf{\Sigma }}_{\hat{J}}-\mathbf{\Sigma }_{J^{\ast
}}\right\Vert _{2}\leq \rho _{\max }\left( \mathbf{G}_{J^{\ast }}^{\top }%
\mathbf{G}_{J^{\ast }}\right) \tilde{\tau}_{N,s_{\ast },1}\delta \left( \log
(M)\right) ^{\frac{2+\alpha }{\alpha }},
\end{equation*}%
where $\tilde{\tau}_{N,s_{\ast },1}=\max (\tilde{A}_{N,s_{\ast },1}^{2},%
\tilde{B}_{N,s_{\ast },1})$, with 
\begin{equation*}
\tilde{A}_{N,s_{\ast },1}=\sqrt{8/3}\Vert \mathbf{\Sigma +\Sigma }%
_{noise}\Vert _{2}^{1/2}\rho _{\min }^{-1/2}\left( \mathbf{G}_{J^{\ast
}}^{\top }\mathbf{G}_{J^{\ast }}\right) \sqrt{\log d^{\ast }}(\log
N)^{1/\alpha }\sqrt{\frac{s_{\ast }}{N}}
\end{equation*}%
and 
\begin{equation*}
\tilde{B}_{N,s_{\ast },1}=\frac{\tilde{\rho}^{2}(\mathbf{\Sigma },\mathbf{%
\Sigma }_{noise})\rho _{\min }^{-1}\left( \mathbf{G}_{J^{\ast }}^{\top }%
\mathbf{G}_{J^{\ast }}\right) }{\sqrt{N}}+\left( \left\Vert \mathbf{\Psi }%
_{J^{\ast }}\right\Vert _{2}+\rho _{\min }^{-1}\left( \mathbf{G}_{J^{\ast
}}^{\top }\mathbf{G}_{J^{\ast }}\right) \left\Vert \mathbf{\Sigma }%
_{noise}\right\Vert _{2}\right) ^{1/2}\tilde{A}_{N,s_{\ast },1}.
\end{equation*}%
Therefore, in the Gaussian case (but also under other assumptions for $X$
such as those in Proposition \ref{prop:orlicz}) the above equations show
that the operator norm $\left\Vert \widehat{\mathbf{\Sigma }}_{\hat{J}}-%
\mathbf{\Sigma }_{J^{\ast }}\right\Vert _{2}^{2}$ depends on the ratio $%
\frac{s_{\ast }}{N}$. Recall that $\left\Vert \mathbf{S}-\mathbf{\Sigma }%
\right\Vert _{2}^{2}$ depends on the ratio $\frac{n}{N}$. Thus, using $%
\widehat{\mathbf{\Sigma }}_{\hat{J}}$ clearly yields significant
improvements if $s_{\ast }$ is small compared to $n$.

To summarize our results let us finally consider the case of an orthogonal
design. Combining Theorems \ref{theo:oracle1}, \ref{theo:activeset} and \ref%
{theo:opnorm} one arrives at the following corrolary:

\begin{corollary}
\label{cor:ortho} Suppose that the observations are i.i.d random variables
from model (\ref{eq:sparseX}). Suppose that $M=n$ and that $\mathbf{G}^{\top
}\mathbf{G}=\mathbf{I}_{n}$ (orthogonal design) and that $X^{0}$ satisfies
Assumption \ref{ass:X}. Let $\epsilon >0$ and $1\leq s_{\ast }\leq \min
(n,M) $. Consider the group Lasso estimator $\widehat{\mathbf{\Sigma }}%
_{\lambda }$ defined by (\ref{GroupLassoEstimator}) with the choices 
\begin{equation*}
\gamma _{k}=2,k=1,\ldots ,n\mbox{ and }\lambda =\Vert \mathbf{\Sigma }%
_{noise}\Vert _{2}\left( 1+\sqrt{\frac{n}{N}}+\sqrt{\frac{2\delta \log M}{N}}%
\right) ^{2}\mbox{ for some constant }\delta >1.
\end{equation*}%
Suppose that 
\begin{equation}
\underset{k\in J^{\ast }}{\min }\left\Vert \mathbf{\mathbf{\Psi }}_{k}^{\ast
}\right\Vert _{\ell _{2}}>2n^{1/2}\tilde{C}_{1}\left( \sigma ,n,s_{\ast
},N,\delta \right) ,  \label{Hyp11}
\end{equation}%
where $\tilde{C}_{1}\left( \sigma ,n,s,N,\delta \right) =\frac{4\left(
1+\epsilon \right) \sqrt{s_{\ast }}}{\epsilon }\sqrt{\tilde{C}_{0}\left(
\sigma ,n,s_{\ast },N,\delta \right) }$ and 
\begin{equation*}
\tilde{C}_{0}\left( \sigma ,n,s_{\ast },N,\delta \right) =(1+\epsilon
)\left( \frac{8}{n}\left\Vert \mathbf{S}-\mathbf{G}\mathbf{\Psi }^{\ast }%
\mathbf{G}^{\top }\right\Vert _{F}^{2}+C(\epsilon )\Vert \mathbf{\Sigma }%
_{noise}\Vert _{2}^{2}\left( 1+\sqrt{\frac{n}{N}}+\sqrt{\frac{2\delta \log M%
}{N}}\right) ^{4}\frac{s_{\ast }}{n}\right) .
\end{equation*}%
Take $\hat{J}:=\left\{ k:\left\Vert \widehat{\mathbf{\mathbf{\Psi }}}%
_{k}\right\Vert _{\ell _{2}}>n^{1/2}\tilde{C}_{1}\left( \sigma ,n,s,N,\delta
\right) \right\} .$ Let $\mathbf{Y}=\mathbf{G}_{J^{\ast }}^{\top }\widetilde{%
\mathbf{X}}$ and $\tilde{Z}=\Vert \mathbf{Y}\Vert _{\ell _{2}}$ . Then, with
probability at least $1-M^{1-\delta }-M^{-\left( \frac{\delta _{\star }}{%
\delta _{\ast }}\right) ^{\frac{\alpha }{2+\alpha }}}$, with $\delta >1$ and 
$\delta _{\star }>\delta _{\ast }$ one has that 
\begin{equation}
\left\Vert \widehat{\mathbf{\Sigma }}_{\hat{J}}-\mathbf{\Sigma }_{J^{\ast
}}\right\Vert _{2}\leq \tilde{\tau}_{N,s_{\ast }}\delta _{\star }\left( \log
(M)\right) ^{\frac{2+\alpha }{\alpha }},  \label{eq:opnorm1}
\end{equation}%
where $\tilde{\tau}_{N,s_{\ast }}=\max (\tilde{A}_{N,s_{\ast }}^{2},\tilde{B}%
_{N,s_{\ast }})$, with $\tilde{A}_{N,s_{\ast }}=\Vert \tilde{Z}\Vert _{\psi
_{\alpha }}\frac{\sqrt{\log d^{\ast }}(\log N)^{1/\alpha }}{\sqrt{N}}$ and $%
\tilde{B}_{N,s_{\ast }}=\frac{\tilde{\rho}^{2}(\mathbf{\Sigma },\mathbf{%
\Sigma }_{noise})}{\sqrt{N}}+\left( \left\Vert \mathbf{\Psi }_{J^{\ast
}}\right\Vert _{2}+\left\Vert \mathbf{\Sigma }_{noise}\right\Vert
_{2}\right) ^{1/2}\tilde{A}_{N,s_{\ast }}.$
\end{corollary}


\subsection{Comparison with the standard Lasso}

In this work, we chose a Group Lasso estimation procedure rather than a standard Lasso. As a matter of fact, for covariance estimation in our setting, the group structure enables to impose a constraint on the number of non zero columns of the matrix $\Psi$ and not on the single entries of the matrix $\Psi$. This corresponds to the natural assumption of obtaining a sparse representation of the process $X(t)$ in the basis given by the functions $g_m$'s and replacing its dimension by its sparsity. Alternatively, the standard Lasso in our setting would be the estimator defined by 
$$
\widehat{\mathbf{\Psi }}_L= \underset{\mathbf{\Psi \in }\mathcal{S}_{M}}{\argmin} \left\{ \left\Vert \widetilde{\mathbf{S}}-\mathbf{G\Psi G}^{\top }\right\Vert _{F}^{2} +2 \lambda \sum_{k=1}^M \sum_{m=1}^M \gamma_{mk} |\Psi _{mk}|\right\},
$$
where $\lambda \geq 0$ is a regularization parameters and the $ \gamma_{mk}$'s are positive weights. This procedure leads to the following Lasso estimator of the covariance matrix $\mathbf{\Sigma }$
\begin{equation}
\widehat{\mathbf{\Sigma }}_{L }=\mathbf{G}\widehat{\mathbf{\Psi }}%
_{L }\mathbf{G}^{\top }\in \mathbb{R}^{n\times n}.
\label{LassoEstimator}
\end{equation}%
In the orthogonal case (i.e. $M=n$ and $\mathbf{G}^{\top }\mathbf{G}=\mathbf{I}_{n}$), this gives rise to the estimator $\widehat{\mathbf{\Psi }}_L$ obtained by  soft thresholding individually each entry $Y_{mk}$ of the matrix   $\mathbf{Y}=\mathbf{G}^{\top }\widetilde{\mathbf{S}}\mathbf{G}$ with the thresholds $\lambda \gamma_{mk}$. Proposition \ref{prop:Lasso} (see below) allows a simple comparison of the statistical performances of the group Lasso estimator $\widehat{\mathbf{\Sigma }}_{L }$ with those of the standard Lasso estimator $\widehat{\mathbf{\Sigma }}_{\lambda}$ in terms of upper bounds for the Frobenius norm. To simplify the discussion, we only consider the orthogonal case and the simple model
\begin{equation}
\widetilde{X}(t_{j})=X^{0}(t_{j})+\mathcal{E}\left( t_{j}\right)
,\;j=1,\ldots ,n,  \label{eq:idealmodel}
\end{equation}%
where the process $X^{0}$ is defined in (\ref{eq:X0}). The statement of the result for the group Lasso is an immediate consequence of Theorem \ref{theo:oracle1}, while the proof to obtain the upper bound for the standard Lasso is an immediate adaptation of the arguments in the proof  of Theorem \ref{theo:oracle1}.

\begin{proposition}
\label{prop:Lasso} Assume that $X$ satisfies model \eqref{eq:idealmodel} and that the covariance matrix  $\mathbf{\Sigma }_{noise}=\mathbb{E}\left( \mathbf{W}_{1}\right) $ of the noise is positive-definite. Consider the group Lasso
estimator $\widehat{\mathbf{\Sigma }}_{\lambda }$   and the standard Lasso estimator $\widehat{\mathbf{\Sigma }}_{L }$ with the choices 
\begin{equation*}
\gamma _{k}=2, \; \gamma _{mk} = 2, \; \lambda =\Vert \mathbf{\Sigma }_{noise}\Vert _{2}\left( 2 +%
\sqrt{\frac{2\delta \log M}{N}}\right) ^{2}\mbox{ for some constant }\delta
>1.
\end{equation*}%
Then, there exist two positive constants $C_{1},C_{2}$ not depending on $n,N,s_{\ast}$ such that with probability at least $1-M^{1-\delta }$ one has that 
\begin{equation*}
\frac{1}{n}\left\Vert \widehat{\mathbf{\Sigma }}_{\lambda }-\mathbf{\Sigma }%
\right\Vert _{F}^{2} \leq    \frac{C_{1}}{n}\left\Vert \mathbf{S}-\mathbf{%
\Sigma }\right\Vert _{F}^{2}  +C_{2}  
\Vert \mathbf{\Sigma }_{noise}\Vert _{2}^{2}\left( 2 +%
\sqrt{\frac{2\delta \log n}{N}}\right) ^{4}
\frac{s_{\ast}}{n} ,  \notag
\end{equation*}%
and
\begin{equation*}
\frac{1}{n}\left\Vert \widehat{\mathbf{\Sigma }}_{L }-\mathbf{\Sigma }%
\right\Vert _{F}^{2} \leq    \frac{C_{1}}{n}\left\Vert \mathbf{S}-\mathbf{%
\Sigma }\right\Vert _{F}^{2}  +C_{2}  
\Vert \mathbf{\Sigma }_{noise}\Vert _{2}^{2}\left( 2 +%
\sqrt{\frac{2\delta \log n}{N}}\right) ^{4}
\frac{s_{\ast}^2}{n}.  \notag
\end{equation*}%
\end{proposition}
Proposition \ref{prop:Lasso} illustrates the advantages of the Group Lasso over the standard Lasso. Indeed, the second term in the upper bound for the group Lasso is much smaller (of the order $\frac{s_{\ast}}{n}$) than the second term in the upper bound for the standard Lasso  (of the order $\frac{s_{\ast}^2}{n}$). This comes from the fact that the sparsity prior of the Group Lasso  is on the number of vanishing columns of the matrix $\Psi$, while the sparsity prior of the standard Lasso only controls the number of non-zero entries of $\Psi$. However, to really demonstrate the benefits  of our method when compared to the performances of the standard Lasso, it is required to also derive lower bounds. This issue is a difficult task which has been considered in  few papers and that is beyond the scope of this paper. For recent work in this direction, we refer to \cite{MR2676881} for regression models or \cite{2010arXiv1007.1771L} and \cite{LPTG} for linear regression and multi-task learning. 

However, the analysis in \cite{MR2676881,2010arXiv1007.1771L} of Group Lasso regularization is carried out the setting of multiple regression models where the parameters to estimate are vectors and with error terms that are centered. Therefore, the results in \cite{MR2676881,2010arXiv1007.1771L} cannot be  applied to the matrix regression model \eqref{eq:regmatrix} since, in our setting, the parameter to estimate is the matrix $\mathbf{\Sigma }$ and the error terms $\mathbf{U}_{i}+\mathbf{W}_{i}$ in \eqref{eq:regmatrix} are not centered.

\section{Numerical experiments and an application to sparse PCA}

\label{sec:simus}

In this section we present some simulated examples to illustrate the
practical behaviour of the covariance matrix estimator by group Lasso
regularization proposed in this paper. In particular, we show its
performances with an application to sparse Principal Components Analysis
(PCA). In the numerical experiments, we use the explicit estimator described
in Proposition \ref{prop:ortho} in the case $M=n$ and an orthogonal design
matrix $\mathbf{G}$, and also the estimator proposed in the more general
situation when $n<M$. The programs for our simulations were implemented
using the MATLAB programming environment.

\subsection{Description of the estimating procedure and the data}

We consider a noisy stochastic processes $\widetilde{X}$ on $\mathbb{T}%
=[0,1] $ with values in $\mathbb{R}$ observed at fixed location points $%
t_{1},...,t_{n}$ in $[0,1]$, generated according to%
\begin{equation}
\widetilde{X}(t_{j})=X^{0}(t_{j})+\sigma \epsilon _{j},\;j=1,\ldots ,n,
\label{modelsimul}
\end{equation}%
where $\sigma >0$ is the level of noise, $\epsilon _{1},\ldots ,\epsilon
_{n} $ are i.i.d.\ standard Gaussian variables, and $X^{0}$ is a random
process independent of the $\epsilon _{j}$'s. For the process $X^{0}$ we
consider two simple models. The first one is given by%
\begin{equation}
X^{0}(t)=af(t),  \label{eq:modsimu1}
\end{equation}%
where $a$ is a Gaussian random coefficient such that $\mathbb{E}a=0$, $%
\mathbb{E}a^{2}=\gamma ^{2}$, and $f:[0,1]\rightarrow \mathbb{R}$ is an
unknown function. The second model for $X^{0}$ is 
\begin{equation}
X^{0}(t)=a_{1}f_{1}(t) + a_{2}f_{2}(t),  \label{eq:modsimu2}
\end{equation}%
where $a_1$ and $a_2$ are independent Gaussian variables such that $\mathbb{E%
}a_1 = \mathbb{E}a_2 =0$, $\mathbb{E}a_1^{2}=\gamma_1 ^{2}$, $\mathbb{E}%
a_2^{2}=\gamma_2 ^{2}$ (with $\gamma_1 > \gamma_2$), and $%
f_1,f_2:[0,1]\rightarrow \mathbb{R}$ are unknown functions. The simulated
data consists in a sample of $N$ independent observations of the process $%
\widetilde{X}$ at the points $t_{1},...,t_{n}$, which are generated
according to (\ref{modelsimul}). 
Therefore, throughout the numerical experiments, one has that
$$
\mathbf{\Sigma }_{noise}=\sigma ^{2}\mathbf{I}_{n}.
$$

In model \eqref{eq:modsimu1}, the covariance matrix $\mathbf{\Sigma }$ of
the process $X^{0}$ at the locations points is given by $\mathbf{\Sigma }%
=\gamma ^{2}\mathbf{F}\mathbf{F}^{\top }$, where by definition $$\mathbf{F}=\left( f\left(
t_{1}\right) ,...,f\left( t_{1}\right) \right) ^{\top }\in \mathbb{R}^{n}.$$
Note that the largest eigenvalue of $\mathbf{\Sigma }$ is $\gamma ^{2}\Vert 
\mathbf{F}\Vert _{\ell _{2}}^{2}$ with corresponding eigenvector $\mathbf{F}$%
. We suppose that the signal $f$ has some sparse representation in a large
dictionary of basis functions of size $M$, given by $\left\{
g_{m},\,m=1,\ldots ,M\right\} $, meaning that $f\left( t\right)
=\sum_{m=1}^{M}\beta _{m}g_{m}\left( t\right) ,$ with $J^{\ast }=\{m,\beta
_{m}\neq 0\}$ of small cardinality $s_{\ast }$. Then, the process $X^{0}$
can be written as $X^{0}(t)=\sum_{m=1}^{M}a\beta _{m}g_{m}\left( t\right) ,$
and thus $\mathbf{\Sigma }=\gamma ^{2}\mathbf{G\Psi }_{J^{\ast }}\mathbf{G}%
^{\top }$, where $\mathbf{\Psi }_{J^{\ast }}$ is an $M\times M$ matrix with
entries equal to $\beta _{m}\beta _{m^{\prime }}$ for $1\leq m,m^{\prime
}\leq M$. 

Similarly, in model \eqref{eq:modsimu2}, the covariance matrix $\mathbf{%
\Sigma }$ of the process $X^{0}$ at the locations points is given by $%
\mathbf{\Sigma }=\gamma_1^{2}\mathbf{F}_1\mathbf{F}_1^{\top } + \gamma_2^{2}%
\mathbf{F}_2\mathbf{F}_2^{\top } $, where by definition
$$\mathbf{F}_1=\left( f_1\left(
t_{1}\right) ,...,f\left( t_{1}\right) \right) ^{\top }\in \mathbb{R}^{n} \mbox{ and }
 \mathbf{F}_2=\left( f_2\left( t_{1}\right) ,...,f\left(t_{1}\right)
\right) ^{\top }\in \mathbb{R}^{n}.$$ In the following simulations, the
functions $f_1$ and $f_2$ are chosen such that $\mathbf{F}_1$ and $\mathbf{F}%
_2$ are orthogonal vectors in $\mathbb{R}^{n}$ with $\Vert \mathbf{F}_1\Vert
_{\ell_{2}} = 1$ and $\Vert \mathbf{F}_2\Vert _{\ell_{2}} = 1$. Under such
an assumption and since $\gamma_1 > \gamma_2$, the largest eigenvalue of $%
\mathbf{\Sigma }$ is $\gamma_1^{2}$ with corresponding eigenvector $\mathbf{F%
}_1$, and the second largest eigenvalue of $\mathbf{\Sigma }$ is $%
\gamma_2^{2}$ with corresponding eigenvector $\mathbf{F}_2$. We suppose that
the signals $f_1$ and $f_2$ have some sparse representations in a large
dictionary of basis functions of size $M$, given by $f_1\left( t\right)
=\sum_{m=1}^{M}\beta _{m}^1 g_{m}\left( t \right),$ and $f_2\left( t\right)
=\sum_{m=1}^{M}\beta _{m}^2 g_{m}\left( t \right)$. Then, the process $X^{0}$
can be written as $X^{0}(t)=\sum_{m=1}^{M}( a_1\beta_{m}^1 + a_2\beta_{m}^2
)g_{m}\left( t\right)$ and thus $\mathbf{\Sigma }=\mathbf{G}( \gamma_1^{2} 
\mathbf{\Psi }^{1} + \gamma_2^{2}\mathbf{\Psi }^{2}) \mathbf{G}^{\top }$,
where $\mathbf{\Psi }^{1},\mathbf{\Psi }^{2} $ are $M\times M$ matrix with
entries equal to $\beta _{m}^1(\beta _{m}^{1})^{\prime }$ and $\beta
_{m}^2(\beta _{m}^{2})^{\prime }$ for $1\leq m,m^{\prime }\leq M$
respectively.

In models \eqref{eq:modsimu1} and \eqref{eq:modsimu2}, we aim at estimating
either $\mathbf{F}$ or $\mathbf{F}_1,\mathbf{F}_2$ by the eigenvectors
corresponding to the largest eigenvalues of the matrix $\widehat{\mathbf{%
\Sigma }}_{\hat{J}} $ defined in (\ref{SigmaHatJhat}), in a high-dimensional
setting with $n > N$ and by using different type of dictionaries. {The idea
behind this is that }$\widehat{\mathbf{\Sigma }}_{\hat{J}}$ is a consistent
estimator of $\mathbf{\Sigma }_{J^{\ast }}$ (see its definition in \ref%
{eq:SigmaJast}) in operator norm. 
Although the matrices $\mathbf{\Sigma }_{J^{\ast }}$ and $\mathbf{\Sigma }$
may have different eigenvectors (depending on the design points and chosen
dictionary), the examples below show the eigenvectors of $\widehat{\mathbf{%
\Sigma }}_{\hat{J}}$ can be used as estimators of the eigenvectors of $%
\mathbf{\Sigma }$.

The estimator $\widehat{\mathbf{\Sigma }}_{\hat{J}}$ of the covariance
matrix $\mathbf{\Sigma }$ is computed as follows. Once the dictionary has
been chosen, we compute the covariance group Lasso (CGL) estimator $\widehat{%
\mathbf{\Sigma }}_{\widehat{\lambda }}=\mathbf{G}\widehat{\mathbf{\Psi }}_{%
\widehat{\lambda }}\mathbf{G}^{\top }$, where $\widehat{\mathbf{\Psi }}_{%
\widehat{\lambda }}$ is {defined in (\ref{MatrixGroupLassoEstimator}). We }%
use a completely data-driven choice for the regularizarion parameter $%
\lambda $, given by $\widehat{\lambda }=\Vert \widehat{\mathbf{\Sigma }%
_{noise}}\Vert _{2}\left( 1+\sqrt{\frac{n}{N}}+\sqrt{\frac{2\delta \log M}{N}%
}\right) ^{2}$, where $\Vert \widehat{\mathbf{\Sigma }_{noise}}\Vert _{2}=$ $%
\widehat{\sigma }^{2}$ is the median absolute deviation (MAD) estimator of $%
\sigma ^{2}$ used in standard wavelet denoising (see e.g.\ \cite{ABS}) and $%
\delta =1.1$. Hence, the method to compute $\widehat{\mathbf{\Sigma }}_{%
\widehat{\lambda }}$ is fully data-driven. Furthermore, we will show in the
examples below that replacing $\lambda $ by $\widehat{\lambda }$ into the
penalized criterion yields a very good practical performance of the
covariance estimation procedure.

As a final step, one needs to compute the estimator $\widehat{\mathbf{\Sigma 
}}_{\hat{J}}$ of $\mathbf{\Sigma }$, as in (\ref{SigmaHatJhat}). For this, {%
we need to have an idea of the true sparsity $s_{\ast }$, since }$\hat{J}$
defined in (\ref{eq:hatJ}) depends on $s_{\ast }$ and also on unknown upper
bounds on the level of noise in the matrix regression model %
\eqref{eq:regmodel} . A similar problem arises in the selection of a sparse
set of active variables when using regularization by a Lasso penalty in the
standard linear model. As an example, recall that in \cite{MR2386087}, a
second thresholding step is aso used to estimate the true sparsity pattern.
However, the suggested thresholding procedure in \cite{MR2386087} also
depends on a priori unknown quantities (such as the amplitude of the
coefficients to estimate). To overcome this drawback in our case, we can
define the final covariance group Lasso (FCGL) estimator as the matrix%
\begin{equation}
{\widehat{\mathbf{\Sigma }}_{\hat{J}}=\mathbf{G}_{\hat{J}}\widehat{\mathbf{%
\Psi }}_{\hat{J}}\mathbf{G}_{\hat{J}}^{\top },}  \label{FCGLE}
\end{equation}%
{with $\hat{J}=\hat{J}_{\epsilon }=\left\{ k:\left\Vert \widehat{\mathbf{%
\mathbf{\Psi }}}_{k}\right\Vert _{\ell _{2}}>\varepsilon \right\} $, where }$%
\varepsilon \ $is a positive constant. 
To select an appropriate value of $\epsilon $, one can plot the cardinality
of $\hat{J}_{\epsilon }$ as a function of $\epsilon $, and then use an
L-curve criterion to only keep in $\hat{J}$ the indices of the columns of $%
\widehat{\mathbf{\Psi }}_{\widehat{\lambda }}$ with a significant value in $%
\ell _{2}$-norm. This choice for $\hat{J}$ is sufficient for numerical
purposes.

In the simulations, to measure the accuracy of the estimation procedure, we
also use the empirical average of the Frobenius and operator norm of the
estimators $\widehat{\mathbf{\Sigma }}_{\hat{\lambda}}$ and $\widehat{%
\mathbf{\Sigma }}_{\hat{J}}$ with respect to the true covariance matrix $%
\mathbf{\Sigma }$ defined by $EAFN=\frac{1}{P}\sum\limits_{p=1}^{P}\left%
\Vert \widehat{\mathbf{\Sigma }}_{\hat{\lambda}}^{p}-\mathbf{\Sigma }%
\right\Vert _{F}$ and $EAON=\frac{1}{P}\sum\limits_{p=1}^{P}\left\Vert 
\widehat{\mathbf{\Sigma }}_{\hat{J}}^{p}-\mathbf{\Sigma }\right\Vert _{2}$%
respectively, over a number $P$ of iterations, where $\widehat{\mathbf{%
\Sigma }}_{\hat{\lambda}}^{p}$ and $\widehat{\mathbf{\Sigma }}_{\hat{J}}^{p}$
are the CGL and FCGL estimators of $\mathbf{\Sigma }$, respectively,
obtained at the $p$-th iteration. We also compute the empirical average of
the operator norm of the estimator $\widehat{\mathbf{\Sigma }}_{\hat{J}}$
with respect to the matrix $\mathbf{\Sigma }_{J^{\ast }}$, defined by $%
EAON^{\ast }=\frac{1}{P}\sum\limits_{p=1}^{P}\left\Vert \widehat{\mathbf{%
\Sigma }}_{\hat{J}}^{p}-\mathbf{\Sigma }_{J^{\ast }}\right\Vert _{2}$.

\subsection{Model \eqref{eq:modsimu1} - case of an orthonormal design (with $%
n=M$)}

First, the size of the dictionary $M$ as well as the basis functions $%
\left\{ g_{m},m=1,...,M\right\} $ have to be specified. In model %
\eqref{eq:modsimu1}, we will use for the test function $f$ the signals
HeaviSine and Blocks (see e.g.\ \cite{ABS} for a definition), and the
Symmlet 8 and Haar wavelet basis for the HeaviSine and Blocks signals
respectively, which are implemented in the Matlab's open-source library
WaveLab (see e.g.\ \cite{ABS} for further references on wavelet methods in
nonparametric statistics). Then, we took $n=M$ and the location points $%
t_{1},...,t_{n}$ are given by the equidistant grid of points $t_{j}=\frac{j}{%
M}$, $j=1,\ldots ,M$ such that the design matrix $\mathbf{G}$ (using either
the Symmlet 8 or the Haar basis) is orthogonal. 

Figures 1, 2, and 3 present the results obtained for a particular simulated
sample of size $N=25$ according to (\ref{modelsimul}), with $n=M=256$, $%
\sigma =0.015$, $\gamma =0.5$ and with $f$ being either the function
HeaviSine or the function Blocks. It can be observed in Figures 1(a) and
1(b) that, as expected in this high dimensional setting ($N<n$), the
empirical eigenvector of $\widetilde{\mathbf{S}}$ associated to its largest
empirical eigenvalue does not lead to a consistent estimator of $\mathbf{F}$.

The CGL estimator $\widehat{\mathbf{\Sigma }}_{\widehat{\lambda }}$ is
computed directly from Proposition \ref{prop:ortho}. In Figures 2(a) and
2(b), we display the eigenvector associated to the largest eigenvalue of $%
\widehat{\mathbf{\Sigma }}_{\widehat{\lambda }}$ as an estimator of $\mathbf{%
F}$. Note that this estimator behaves poorly. The estimation considerably
improves by taking the FCGL estimator $\widehat{\mathbf{\Sigma }}_{\hat{J}}$
defined in (\ref{FCGLE}).{\ Figures 3(a) and 3(b) illustrate the very good
performance of the eigenvector associated to the largest eigenvalue of the
matrix $\widehat{\mathbf{\Sigma }}_{\hat{J}}$ as an estimator of }$\mathbf{F}
$.

It is clear that the estimators $\widehat{\mathbf{\Sigma }}_{\widehat{\lambda }}$ and {$\widehat{\mathbf{\Sigma }}_{\hat{J}}$ are }random matrices
that depend on the observed sample. Tables 1(a) and 1(b) show the values of $%
EAFN$, $EAON$ and $EAON^{\ast }$ corresponding to $P=100$ simulated samples
of different sizes $N$ and different values of the level of noise $\sigma $.
It can be observed that for both signals the empirical averages $EAFN$, $%
EAON $ and $EAON^{\ast }$ behaves similarly, being the values of $EAON$
smaller than its corresponding values of $EAFN$ as expected. Observing each
table separately we can remark that, for $N$ fixed, when the level of noise $%
\sigma $ increases then the values of $EAFN$, $EAON$ and $EAON^{\ast }$ also
increase. By simple inspection of the values of $EAFN$, $EAON$ and $%
EAON^{\ast }$ in the same position at Tables 1(a) and 1(b) we can check
that, for $\sigma $ fixed, when the number of replicates $N$ increases then
the values of $EAFN$, $EAON$ and $EAON^{\ast }$ decrease in all cases. We
can also observe {how the difference between }$EAON$ and $EAON^{\ast }${\ is
bigger as the level of noise increases. }

\begin{center}
$\overset{%
\begin{array}{c}
\text{Table 1(a). Values of }EAFN\text{, }EAON\text{ and }EAON^{\ast }\text{
corresponding to signals } \\ 
\text{HeaviSine and Blocks for }M=n=256\text{, }N=25\text{.}%
\end{array}%
}{%
\begin{tabular}{|c|c|c|c|c|c|c|c|}
\hline
Signal & $\sigma $ & $0.005$ & $0.01$ & $0.05$ & $0.1$ & $0.5$ & $1$ \\ 
\hline
HeaviSine & $EAFN$ & 0.0634 & 0.0634 & 0.2199 & 0.2500 & 0.2500 & 0.2500 \\ 
\hline
HeaviSine & $EAON$ & 0.0619 & 0.0569 & 0.1932 & 0.2500 & 0.2500 & 0.2500 \\ 
\hline
HeaviSine & $EAON^{\ast }$ & 0.0619 & 0.0569 & 0.1943 & 0.2600 & 0.5000 & 
1.2500 \\ \hline
Blocks & $EAFN$ & 0.0553 & 0.0681 & 0.2247 & 0.2500 & 0.2500 & 0.2500 \\ 
\hline
Blocks & $EAON$ & 0.0531 & 0.0541 & 0.2083 & 0.2500 & 0.2500 & 0.2500 \\ 
\hline
Blocks & $EAON^{\ast }$ & 0.0531 & 0.0541 & 0.2107 & 0.2600 & 0.5000 & 1.2500
\\ \hline
\end{tabular}%
}$

$\overset{%
\begin{array}{c}
\text{Table 1(b). Values of }EAFN\text{, }EAON\text{ and }EAON^{\ast }\text{
corresponding to signals } \\ 
\text{HeaviSine and Blocks for }M=n=256\text{, }N=40\text{.}%
\end{array}%
}{%
\begin{tabular}{|c|c|c|c|c|c|c|c|}
\hline
Signal & $\sigma $ & $0.005$ & $0.01$ & $0.05$ & $0.1$ & $0.5$ & $1$ \\ 
\hline
HeaviSine & $EAFN$ & 0.0501 & 0.0524 & 0.1849 & 0.2499 & 0.2500 & 0.2500 \\ 
\hline
HeaviSine & $EAON$ & 0.0496 & 0.0480 & 0.1354 & 0.2496 & 0.2500 & 0.2500 \\ 
\hline
HeaviSine & $EAON^{\ast }$ & 0.0496 & 0.0480 & 0.1366 & 0.2596 & 0.5000 & 
1.2500 \\ \hline
Blocks & $EAFN$ & 0.0485 & 0.0494 & 0.2014 & 0.2500 & 0.2500 & 0.2500 \\ 
\hline
Blocks & $EAON$ & 0.0483 & 0.0429 & 0.1871 & 0.2500 & 0.2500 & 0.2500 \\ 
\hline
Blocks & $EAON^{\ast }$ & 0.0483 & 0.0429 & 0.1893 & 0.2600 & 0.5000 & 1.2500
\\ \hline
\end{tabular}%
}$
\end{center}

\subsection{Model \eqref{eq:modsimu2} - the case $M=2n$ by mixing two
orthonormal basis}

Consider now the setting of model \eqref{eq:modsimu2} with $\gamma_1=0.5$, $%
\gamma_2=0.2$, $\sigma = 0.045$, $N=25$ and an equidistant grid of design points $%
t_{1},...,t_{n}$ given by $t_{j}=\frac{j}{n}$, $j=1,\ldots ,n$ with $n=128$.
For the signals $f_1$ and $f_2$ we took the test functions displayed in
Figure 4(a) and 4(b). Obviously, the signal $f_1$ has a sparse
representation in a Haar basis while the signal $f_2$ has a sparse
representation in a Fourier basis. Thus, this suggests to construct a
dictionary by mixing two orthonormal basis. More precisely, we construct a $%
n \times n$ orthogonal matrix $\mathbf{G}^1$ using the Haar basis and a $n
\times n$ orthogonal matrix $\mathbf{G}^2$ using a Fourier basis (cosine and
sine at various frequencies) at the design points. Then, we form the $n
\times M$ design matrix $\mathbf{G} = [\mathbf{G}^1 \; \mathbf{G}^2]$ with $%
M=2n$. The CGL estimator $\widehat{\mathbf{\Sigma }}_{\widehat{\lambda }}$
is computed by the minimization procedure {(\ref{MatrixGroupLassoEstimator})
using the }Matlab package \textit{minConf} of \cite{Schmidt}.

In Figures 5(a) and 5(b), we display the eigenvector associated to the
largest eigenvalue of $\widehat{\mathbf{\Sigma }}_{\widehat{\lambda }}$ as
an estimator of $\mathbf{F}_1$, and the eigenvector associated to the second
largest eigenvalue of $\widehat{\mathbf{\Sigma }}_{\widehat{\lambda }}$ as
an estimator of $\mathbf{F}_2$. Note that these estimators behaves poorly.
The estimation considerably improves by taking the FCGL estimator $\widehat{%
\mathbf{\Sigma }}_{\hat{J}}$ defined in (\ref{FCGLE}).{\ Figures 6(a) and
6(b) illustrate the very good performance of the eigenvectors associated to
the largest eigenvalue and second largest eigenvalue of the matrix $\widehat{%
\mathbf{\Sigma }}_{\hat{J}}$ as estimators of }$\mathbf{F}_1$ and $\mathbf{F}%
_2$.

Finally, to illustrate the benefits of mixing two orthonormal basis, we also
display in Figures 7 and 8 the estimation of $\mathbf{F}_1$ and $\mathbf{F}%
_2 $ when computing the matrix $\widehat{\mathbf{\Sigma }}_{\hat{J}}$ by
using either only the Haar basis (i.e.\ $\mathbf{G} = \mathbf{G}^1$ and $M=n$%
) or only the Fourier basis (i.e.\ $\mathbf{G} = \mathbf{G}^1$ and $M=n$).
The results are clearly much worse and not satisfactory.

\subsection{Model \eqref{eq:modsimu1} - case of non equispaced design points
such that $n<M$}

Let us now return to the setting of model \eqref{eq:modsimu1}. The test
functions $f$ are either the signal HeaviSine and or the signal Blocks. We
also use the Symmlet 8 and Haar wavelet basis for the HeaviSine and Blocks
functions respectively. However, we now choose to take a setting where the
number of design points $n$ is smaller than the size $M$ of the dictionary.
Taking $n<M$, the location points are given by a subset $\left\{
t_{1},...,t_{n}\right\} \subset \{\frac{k}{M}:k=1,...,M\}$ of size $n$, such
that the design matrix $\mathbf{G}$ is an $n\times M$ matrix (using either
the Symmlet 8 and Haar basis). For {a fixed value of $n$, the }subset $%
\left\{ t_{1},...,t_{n}\right\} ${\ is chosen by taking the first }$n$
points obtained from {a }random permutation of the elements of the set $\{%
\frac{1}{M},\frac{2}{M},...,1\}$. Figures 9 and 10 present the results
obtained for a particular simulated sample of size $N=25$ according to (\ref%
{modelsimul}), with $n=90$, $M=128$, $\sigma =0.02$, $\gamma =0.5$ and with $%
f$ being either the function HeaviSine or the function Blocks.
It can be observed in Figures 9(a) and 9(b) that, as expected in this high
dimensional setting ($N<n$), the empirical eigenvector of $\widetilde{%
\mathbf{S}}$ associated to its largest empirical eigenvalue are noisy
versions of $\mathbf{F}$. As explained previously, the CGL estimator $%
\widehat{\mathbf{\Sigma }}_{\widehat{\lambda }}$ is computed by the
minimization procedure {(\ref{MatrixGroupLassoEstimator}) using the }Matlab
package \textit{minConf} of \cite{Schmidt}. In Figures 10(a) and 10(b) is
shown the eigenvector associated to the largest eigenvalue of $\widehat{%
\mathbf{\Sigma }}_{\widehat{\lambda }}$ as an estimator of $\mathbf{F}$.
Note that this estimator is quite noisy. Again, the {eigenvector associated
to the largest eigenvalue }of{\ the matrix $\widehat{\mathbf{\Sigma }}_{%
\widehat{J}}{\ }$defined in (\ref{FCGLE}) is much a better estimator of }$%
\mathbf{F}${.} This is {illustrated in Figures 11(a) and 11(b)}. To compare the accuracy of the estimators for different simulated samples,
we compute the values of $EAFN$, $EAON$ and $EAON^{\ast }$ with fixed
values of $\sigma =0.05$, $M=128$, $N=40$, $P=50$ {for different values of 
} the number of design points $n${. For all the values of $n$ considered,
the design points }$t_{1},...,t_{n}$ are selected as {the first }$n$ points
obtained from {the same }random permutation of the elements of the set $\{%
\frac{1}{M},\frac{2}{M},...,1\}$. The chosen subset $\left\{
t_{1},...,t_{n}\right\} $ is used for all the $P$ iterations needed in the
computation of the empirical averages (fixed design over the iterations).
Figure 12 shows the values of $EAFN$, $EAON$ and $EAON^{\ast }$ obtained for
each value of $n$ for both signals HeaviSine and Blocks. It can be observed
that the values of the empirical averages $EAON$ and $EAON^{\ast }$ are much
smaller than its corresponding values of $EAFN$ as expected. We can remark
that, 
when $n$ increases, the values of $EAFN$, $EAON$ and $EAON^{\ast }$ first
increase and then decrease, and the change of monotony occurs when $n>N$.
Note that the case $n=M=128$ is included in these results.


\parbox{7.5cm}{\
\begin{center}
\textbf{Orthonormal case - Model \eqref{eq:modsimu1}}
 \pgfuseimage{Fig1H}
 \textbf{Figure 1(a). Signal HeaviSine and Eigenvector associated to the largest eigenvalue of $\widetilde{\mathbf{S}}$}\vspace{.1cm}\\
 \pgfuseimage{Fig2H}
 \textbf{Figure 2(a). Signal HeaviSine and Eigenvector associated to the largest eigenvalue of $\widehat{\mathbf{\Sigma}}_{\widehat{\lambda }}$}\vspace{.1cm}\\
 \pgfuseimage{Fig3H}
 \textbf{Figure 3(a). Signal HeaviSine and Eigenvector associated to the largest eigenvalue of $\widehat{\mathbf{\Sigma}}_{\widehat{J}}$}\vspace{.1cm}\\
  \end{center}
}\hfill 
\parbox{7.5cm}{\vspace{.5cm}
   \begin{center}
 \pgfuseimage{Fig1B}
 \textbf{Figure 1(b). Signal Blocks and Eigenvector associated to the largest eigenvalue of $\widetilde{\mathbf{S}}$}\vspace{.1cm}\\
 \pgfuseimage{Fig2B}
 \textbf{Figure 2(b). Signal Blocks and Eigenvector associated to the largest eigenvalue of $\widehat{\mathbf{\Sigma}}_{\widehat{\lambda }}$}\vspace{.1cm}\\
 \pgfuseimage{Fig3B}
 \textbf{Figure 3(b). Signal Blocks and Eigenvector associated to the largest eigenvalue of $\widehat{\mathbf{\Sigma}}_{\widehat{J}}$}\vspace{.1cm}\\
  \end{center}}

\parbox{7.5cm}{\
  \textbf{Case $M=2n$ (Haar + Fourier basis)}
\begin{center}
 \pgfuseimage{f1}
 \textbf{Figure 4(a). Signal $F_1$}\vspace{.1cm}\\
 \pgfuseimage{hatf1_naive}
 \textbf{Figure 5(a). Signal $F_1$ and Eigenvector associated to the largest eigenvalue of $\widehat{\mathbf{\Sigma }}_{\widehat{\lambda }}$}\vspace{.1cm}\\
   \pgfuseimage{hatf1_GLF}
 \textbf{Figure 6(b). Signal $F_1$ and Eigenvector associated to the largest eigenvalue of $\widehat{\mathbf{\Sigma}}_{\widehat{J}}$ with  $\mathbf{G} = [\mathbf{G}^1 \; \mathbf{G}^2]$}\vspace{.1cm}\\
  \end{center}
}\hfill 
\parbox{7.5cm}{\vspace{.5cm}
   \begin{center}
 \pgfuseimage{f2}
 \textbf{Figure 4(b). Signal $F_2$}\vspace{.1cm}\\
 \pgfuseimage{hatf2_naive}
 \textbf{Figure 5(b). Signal $F_2$ and Eigenvector associated to the second largest eigenvalue of $\widehat{\mathbf{\Sigma }}_{\widehat{\lambda }}$}\vspace{.1cm}\\
  \pgfuseimage{hatf2_GLF}
 \textbf{Figure 6(b). Signal $F_2$ and Eigenvector associated to the second largest eigenvalue of $\widehat{\mathbf{\Sigma}}_{\widehat{J}}$ with  $\mathbf{G} = [\mathbf{G}^1 \; \mathbf{G}^2]$}\vspace{.1cm}\\
  \end{center}}

\parbox{7.5cm}{\
    \textbf{Orthonormal case $M=n$ (Haar)}
\begin{center}
 \pgfuseimage{hatf1_Haar}
 \textbf{Figure 7(a). Signal $F_1$ and Eigenvector associated to the largest eigenvalue of $\widehat{\mathbf{\Sigma}}_{\widehat{J}}$ with $\mathbf{G} = \mathbf{G}^1$}\vspace{.1cm}\\
  \end{center}
}\hfill 
\parbox{7.5cm}{\vspace{.5cm}
   \begin{center}
 \pgfuseimage{hatf2_Haar}
 \textbf{Figure 7(b). Signal $F_2$ and Eigenvector associated to the second largest eigenvalue of $\widehat{\mathbf{\Sigma}}_{\widehat{J}}$  with $\mathbf{G} = \mathbf{G}^1$}\vspace{.1cm}\\
  \end{center}}

\parbox{7.5cm}{\
     \textbf{Orthonormal case $M=n$ (Fourier)}
\begin{center}

 \pgfuseimage{hatf1_Fourier}
 \textbf{Figure 8(a). Signal $F_1$ and Eigenvector associated to the largest eigenvalue of $\widehat{\mathbf{\Sigma}}_{\widehat{J}}$  with $\mathbf{G} = \mathbf{G}^2$}\vspace{.1cm}\\
  \end{center}
}\hfill 
\parbox{7.5cm}{\vspace{.5cm}
   \begin{center}
 \pgfuseimage{hatf2_Fourier}
 \textbf{Figure 8(b). Signal $F_2$ and Eigenvector associated to the second largest eigenvalue of $\widehat{\mathbf{\Sigma}}_{\widehat{J}}$  with $\mathbf{G} = \mathbf{G}^2$}\vspace{.1cm}\\
  \end{center}}

\parbox{7.5cm}{\
\begin{center}
\textbf{Non equi-spaced points with $n<M$}
 \pgfuseimage{Fig4Hnn}
 \textbf{Figure 9(a). Signal HeaviSine and Eigenvector associated to the largest eigenvalue of $\widetilde{\mathbf{S}}$}\vspace{.1cm}\\
 \pgfuseimage{Fig5Hnn}
 \textbf{Figure 10(a). Signal HeaviSine and Eigenvector associated to the largest eigenvalue of $\widehat{\mathbf{\Sigma}}_{\widehat{\lambda }}$}\vspace{.1cm}\\
 \pgfuseimage{Fig6Hnn}
 \textbf{Figure 11(a). Signal HeaviSine and Eigenvector associated to the largest eigenvalue of $\widehat{\mathbf{\Sigma}}_{\widehat{J}}$}\vspace{.1cm}\\
  \end{center}
}\hfill 
\parbox{7.5cm}{\vspace{.6cm}
   \begin{center}
 \pgfuseimage{Fig4Bn}
 \textbf{Figure 9(b). Signal Blocks and Eigenvector associated to the largest eigenvalue of $\widetilde{\mathbf{S}}$}\vspace{.1cm}\\
 \pgfuseimage{Fig5Bn}
 \textbf{Figure 10(b). Signal Blocks and Eigenvector associated to the largest eigenvalue of $\widehat{\mathbf{\Sigma}}_{\widehat{\lambda }}$}\vspace{.1cm}\\
 \pgfuseimage{Fig6Bn}
 \textbf{Figure 11(b). Signal Blocks and Eigenvector associated to the largest eigenvalue of $\widehat{\mathbf{\Sigma}}_{\widehat{J}}$}\vspace{.1cm}\\
  \end{center}}

\parbox{7.5cm}{\
\begin{center}
 \pgfuseimage{FigHN40}
 \textbf{Figure 12(a). Values of $EAFN$, $EAON$ and $EAON^{\ast}$ for Signal HeaviSine as a function of $n$}\vspace{.1cm}\\
 \end{center}
}\hfill 
\parbox{7.5cm}{
   \begin{center}
 \pgfuseimage{FigBN40}
 \textbf{Figure 12(b). Values of $EAFN$, $EAON$ and $EAON^{\ast}$ for Signal Blocks as a function of $n$}\vspace{.1cm}\\
\end{center}}

\renewcommand*\thesection{A}



\appendix


\section{}

\subsection{Notations}

First let us introduce some notations and properties that will be used
throughout this Appendix. The vectorization of a $p\times q$ matrix $\mathbf{%
A}=(a_{ij})_{1\leq i\leq p,1\leq j\leq q}$ is the $pq\times 1$ column vector
denoted by $vec\left( \mathbf{A}\right) $, obtain by stacking the columns of
the matrix $\mathbf{A}$ on top of one another. That is $vec(\mathbf{A}%
)=[a_{11},...,a_{p1},a_{12},...,a_{p2},...,a_{1q},...,a_{pq}]^{\top }$. If $%
\mathbf{A}=(a_{ij})_{1\leq i\leq k,1\leq j\leq n}$ is a $k\times n$ matrix
and $\mathbf{B}=(b_{ij})_{1\leq i\leq p,1\leq j\leq q}$ is a $p\times q$
matrix, then the Kronecker product of the two matrices, denoted by $\mathbf{A%
}\otimes \mathbf{B}$, is the $kp\times nq$ block matrix%
\begin{equation*}
\mathbf{A}\otimes \mathbf{B=}%
\begin{bmatrix}
a_{11}\mathbf{B} & . & . & . & a_{1n}\mathbf{B} \\ 
. & . &  &  & . \\ 
. &  & . &  & . \\ 
. &  &  & . & . \\ 
a_{k1}\mathbf{B} & . & . & . & a_{kn}\mathbf{B}%
\end{bmatrix}%
.
\end{equation*}%
In what follows, we repeatedly use the fact that the Frobenius norm is
invariant by the $vec$ operation meaning that 
\begin{equation}
\left\Vert \mathbf{A}\right\Vert _{F}^{2}=\Vert vec\left( \mathbf{A}\right)
\Vert _{\ell _{2}}^{2},  \label{eq:P1}
\end{equation}%
and the properties that 
\begin{equation}
vec\left( \mathbf{ABC}\right) =\left( \mathbf{C}^{\top }\mathbf{\otimes A}%
\right) vec\left( \mathbf{B}\right) ,  \label{eq:P2}
\end{equation}%
and 
\begin{equation}
(\mathbf{A}\otimes \mathbf{B})(\mathbf{C}\otimes \mathbf{D})=\mathbf{A}%
\mathbf{C}\otimes \mathbf{B}\mathbf{D},  \label{eq:P3}
\end{equation}%
provided the above matrix products are compatible.

\subsection{Proof of Proposition \protect\ref{prop:ortho}}

\begin{lemma}
\label{lem:KKT} Let $\widehat{\mathbf{\Psi }}=\widehat{\mathbf{\Psi }}%
_{\lambda }$ denotes the solution of \eqref{MatrixGroupLassoEstimator}.
Then, for $k=1,\ldots ,M$ 
\begin{eqnarray*}
\left[ (\mathbf{G}\otimes \mathbf{G})^{\top }\left( vec(\widetilde{\mathbf{S}%
})-(\mathbf{G}\otimes \mathbf{G})vec(\widehat{\mathbf{\Psi }})\right) \right]
^{k} &=&\lambda \gamma _{k}\frac{\widehat{\mathbf{\Psi }}_{k}}{\Vert 
\widehat{\mathbf{\Psi }}_{k}\Vert _{\ell _{2}}}\quad \mbox{ if }\quad 
\mathbf{\Psi }_{k}\neq 0 \\
\left\Vert \left[ (\mathbf{G}\otimes \mathbf{G})^{\top }\left( vec(%
\widetilde{\mathbf{S}})-(\mathbf{G}\otimes \mathbf{G})vec(\widehat{\mathbf{%
\Psi }})\right) \right] ^{k}\right\Vert _{\ell _{2}} &\leq &\lambda \gamma
_{k}\qquad \qquad \;\,\mbox{ if }\quad \widehat{\mathbf{\Psi }}_{k}=0
\end{eqnarray*}%
where $\widehat{\mathbf{\Psi }}_{k}$ denotes the $k$-th column of the matrix 
$\widehat{\mathbf{\Psi }}$ and the notation $\left[ \beta \right] ^{k}$
denotes the vector $(\beta _{k,m})_{m=1,\ldots ,M}$ in $\mathbb{R}^{M}$ for
a vector $\beta =(\beta _{k,m})_{k,m=1,\ldots ,M}\in \mathbb{R}^{M^{2}}$.
\end{lemma}

\noindent \textbf{Proof of Lemma \ref{lem:KKT}} For $\mathbf{\Psi } \in \mathbb{R}^{M
\times M}$ define 
\begin{equation*}
L(\mathbf{\Psi } ) = \left\Vert \widetilde{\mathbf{S}}-\mathbf{G\Psi G}%
^{\top }\right\Vert _{F}^{2} = \left\| vec({\widetilde{\mathbf{S}}}) - (%
\mathbf{G}\otimes \mathbf{G}) vec( \mathbf{\Psi } ) \right\|_{\ell_{2}}^2,
\end{equation*}
and remark that $\widehat{\mathbf{\Psi }}$ is the solution of the convex
optimization problem 
\begin{equation*}
\widehat{\mathbf{\Psi }} =\underset{\mathbf{\Psi \in }\mathcal{S}_{M}}{%
\argmin}\left\{ L(\mathbf{\Psi } ) +2\lambda \sum_{k=1}^{M}\gamma _{k}\sqrt{%
\sum_{m=1}^{M}\Psi _{mk}^{2}}\right\}.
\end{equation*}
It follows from standard arguments in convex analysis (see e.g.\ \cite%
{MR2061575}), that $\widehat{\mathbf{\Psi }}$ is a solution of the above
minimization problem if and only if 
\begin{equation*}
- \nabla L(\widehat{\mathbf{\Psi }}) \in 2 \lambda \partial \left(
\sum_{k=1}^{M}\gamma _{k}\sqrt{\sum_{m=1}^{M} \hat{\Psi}_{mk}^{2}} \right)
\end{equation*}
where $\nabla L(\widehat{\mathbf{\Psi }}) $ denotes the gradient of $L$ at $%
\widehat{\mathbf{\Psi }}$ and $\partial$ denotes the subdifferential given
by 
\begin{equation*}
\partial \left( \sum_{k=1}^{M}\gamma _{k}\sqrt{\sum_{m=1}^{M}\Psi _{mk}^{2}}
\right) = \left\{ \mathbf{\Theta} \in \mathbb{R}^{M \times M} : \mathbf{%
\Theta}_{k} = \gamma_{k} \frac{ \mathbf{\Psi }_{k} }{\| \mathbf{\Psi }_{k}
\|_{\ell_{2}}} \mbox{ if } \mathbf{\Psi }_{k} \neq 0, \left\| \mathbf{\Theta}%
_{k} \right\|_{\ell_{2}} \leq \gamma_{k}\mbox{ if } \mathbf{\Psi }_{k} = 0
\right\}
\end{equation*}
where $\mathbf{\Theta}_{k}$ denotes the $k$-th column of $\mathbf{\Theta}
\in \mathbb{R}^{M \times M}$ which completes the proof. \hfill $\Box $ 
\newline

Now, let $\mathbf{\Psi }\in \mathcal{S}_{M}$ with $M=n$ and suppose that $%
\mathbf{G}^{\top} \mathbf{G} = \mathbf{I}_{n}$. Let $\mathbf{Y}= (\mathbf{Y}%
_{mk})_{1\leq m,k\leq M} = \mathbf{G}^{\top }\widetilde{\mathbf{S}}\mathbf{G}
$ and remark that $vec( \mathbf{Y} ) = \left( \mathbf{G} \otimes \mathbf{G}
\right) ^{\top } vec(\widetilde{\mathbf{S}})$. Then, by using Lemma \ref%
{lem:KKT} and the fact that $\mathbf{G}^{\top }\mathbf{G}=\mathbf{I}_{n}$
implies that $(\mathbf{G}\otimes \mathbf{G})^{\top }\left( \mathbf{G}\otimes 
\mathbf{G}\right) =\mathbf{I}_{n^{2}}$, it follows that $\widehat{\mathbf{%
\Psi }}=\widehat{\mathbf{\Psi }}_{\lambda }$ satisfies for $k=1,\ldots ,M$
the following equations 
\begin{equation*}
\widehat{\mathbf{\Psi }}_{k}\left( 1+\frac{\lambda \gamma _{k}}{\sqrt{%
\sum_{m=1}^{M}\widehat{\Psi }_{mk}^{2}}}\right) =\mathbf{Y}_{k}%
\mbox{
for all }\widehat{\mathbf{\Psi }}_{k}\neq 0,
\end{equation*}%
and 
\begin{equation*}
\sqrt{\sum_{m=1}^{M}\mathbf{Y}_{mk}^{2}}\leq \lambda \gamma _{k}%
\mbox{ for
all }\widehat{\mathbf{\Psi }}_{k}=0.
\end{equation*}%
where $\widehat{\mathbf{\Psi }}_{k}=(\widehat{\Psi }_{mk})_{1\leq m\leq
M}\in \mathbb{R}^{M}$ and $\mathbf{Y}_{k}=(\mathbf{Y}_{mk})_{1\leq m\leq
M}\in \mathbb{R}^{M}$, which implies that the solution is given by 
\begin{equation*}
\widehat{\Psi }_{mk}=\left\{ 
\begin{array}{ccc}
0 & \mbox{ if } & \sqrt{\sum_{m=1}^{M}\mathbf{Y}_{mk}^{2}}\leq \lambda
\gamma _{k} \\ 
Y_{mk}\left( 1-\frac{\lambda \gamma _{k}}{\sqrt{\sum_{j=1}^{M}\mathbf{Y}%
_{jk}^{2}}}\right) & \mbox{ if } & \sqrt{\sum_{m=1}^{M}\mathbf{Y}_{mk}^{2}}%
>\lambda \gamma _{k}%
\end{array}%
\right.
\end{equation*}%
which completes the proof of Proposition \ref{prop:ortho}. \hfill $\Box $

\subsection{Proof of Proposition \protect\ref{prop:orlicz}}

First suppose that $X$ is Gaussian. Then, remark that for $Z=\Vert \mathbf{X}%
\Vert _{\ell _{2}}$, one has that $\Vert Z\Vert _{\psi _{2}}<+\infty $ which
implies that $\Vert Z\Vert _{\psi _{2}}=\Vert Z^{2}\Vert _{\psi _{1}}^{1/2}$%
. Since $Z^{2}=\sum_{i=1}^{n}|X(t_{i})|^{2}$ it follows that 
\begin{equation*}
\Vert Z^{2}\Vert _{\psi _{1}}\leq \sum_{i=1}^{n}\Vert Z_{i}^{2}\Vert _{\psi
_{1}}=\sum_{i=1}^{n}\Vert Z_{i}\Vert _{\psi _{2}}^{2}=\sum_{i=1}^{n}\mathbf{%
\Sigma }_{ii}\Vert \mathbf{\Sigma }_{ii}^{-1/2}Z_{i}\Vert _{\psi _{2}}^{2},
\end{equation*}%
where $Z_{i}=X(t_{i}),i=1,\ldots ,n$ and $\mathbf{\Sigma }_{ii}$ denotes the 
$i$th diagonal element of $\mathbf{\Sigma }$. Then, the result follows by
noticing that $\Vert Y\Vert _{\psi _{2}}\leq \sqrt{8/3}$ if $Y\sim N(0,1)$.
The proof for the case where $X$ is such that $\Vert Z\Vert _{\psi
_{2}}<+\infty $ and there exists a constant $C_{1}$ such that $\Vert \mathbf{%
\Sigma }_{ii}^{-1/2}Z_{i}\Vert _{\psi _{2}}\leq C_{1}$ for all $i=1,\ldots
,n $ follows from the same arguments.

Now, consider the case where $X$ is a bounded process. Since there exists a
constant $R>0$ such that for all $t\in \mathbb{T}$, $|X(t)|\leq R$, it
follows that for $Z=\Vert \mathbf{X}\Vert _{\ell _{2}}$ then $Z\leq \sqrt{n}%
R $ which implies that for any $\alpha \geq 1$, $\Vert Z\Vert _{\psi
_{\alpha }}\leq \sqrt{n}R(\log 2)^{-1/\alpha }$, (by definition of the norm $%
\Vert Z\Vert _{\psi _{\alpha }}$) which completes the proof of Proposition %
\ref{prop:orlicz}. \hfill $\Box $



\subsection{Proof of Proposition \protect\ref{prop:MPsparse}}

Under the assumption that $X=X^{0}$, it follows that $\mathbf{\Sigma }=%
\mathbf{G}\mathbf{\Psi }^{\ast }\mathbf{G}^{\top }$ with $\mathbf{\Psi }%
^{\ast }=\mathbb{E}\left( \mathbf{a}\mathbf{a}^{\top }\right) $, where $%
\mathbf{a}$ is the random vector of $\mathbb{R}^{M}$ with $\mathbf{a}%
_{m}=a_{m}$ for $m\in J^{\ast }$ and $\mathbf{a}_{m}=0$ for $m\notin J^{\ast
}$. Then, define the random vector $\mathbf{a}_{J^{\ast }}\in \mathbb{R}%
^{J^{\ast }}$ whose coordinates are the random coefficients $a_{m}$ for $%
m\in J^{\ast }$. Let $\mathbf{\Psi }_{J^{\ast }}=\mathbb{E}\left( \mathbf{a}%
_{J^{\ast }}\mathbf{a}_{J^{\ast }}^{\top }\right) $. Note that $\mathbf{%
\Sigma }=\mathbf{G}_{J^{\ast }}\mathbf{\Psi }_{J^{\ast }}\mathbf{G}_{J^{\ast
}}^{\top }$ and $\mathbf{S}=\mathbf{G}_{J^{\ast }}\widehat{\mathbf{\Psi }}%
_{J^{\ast }}\mathbf{G}_{J^{\ast }}^{\top }$, with $\widehat{\mathbf{\Psi }}%
_{J^{\ast }}=\frac{1}{N}\sum_{i=1}^{N}\mathbf{a}_{J^{\ast }}^{i}(\mathbf{a}%
_{J^{\ast }}^{i})^{\top }$, where $\mathbf{a}_{J^{\ast }}^{i}\in \mathbb{R}%
^{J^{\ast }}$ denotes the random vector whose coordinates are the random
coefficients $a_{m}^{i}$ for $m\in J^{\ast }$ such that $X_{i}(t)=\sum_{m\in
J^{\ast }}a_{m}^{i}g_{m}(t),\;t\in \mathbb{T}$.

Therefore, $\widehat{\mathbf{\Psi }}_{J^{\ast }}$ is a sample covariance
matrix of size $s_{\ast }\times s_{\ast }$ and we can control its deviation
in operator norm from $\widehat{\mathbf{\Psi }}_{J^{\ast }}$ by using
Proposition \ref{prop:mp}. For this we simply have to verify conditions
similar to \textbf{(A1)} and \textbf{(A2)} in Assumption \ref{ass:X} for the
random vector $\mathbf{a}_{J^{\ast }}=(\mathbf{G}_{J^{\ast }}^{\top }\mathbf{%
G}_{J^{\ast }})^{-1}\mathbf{G}_{J^{\ast }}^{\top }\mathbf{X}\in \mathbb{R}%
^{s_{\ast }}$. First, let $\beta \in \mathbb{R}^{s_{\ast }}$ with $\Vert
\beta \Vert _{\ell _{2}}=1$. Then, remark that $\mathbf{a}_{J^{\ast }}^{\top
}\beta =\mathbf{X}^{\top }\tilde{\beta}$ with $\tilde{\beta}=\mathbf{G}%
_{J^{\ast }}\left( \mathbf{G}_{J^{\ast }}^{\top }\mathbf{G}_{J^{\ast
}}\right) ^{-1}\beta $. Since $\Vert \tilde{\beta}\Vert _{\ell _{2}}\leq
\left( \rho _{\min }\left( \mathbf{G}_{J^{\ast }}^{\top }\mathbf{G}_{J^{\ast
}}\right) \right) ^{-1/2}$ and using that $X$ satisfies Assumption \ref%
{ass:X} it follows that 
\begin{equation}
\left( \mathbb{E}|\mathbf{a}_{J^{\ast }}^{\top }\beta |^{4}\right)
^{1/4}\leq \rho \left( \mathbf{\Sigma }\right) \rho _{\min }^{-1/2}\left( 
\mathbf{G}_{J^{\ast }}^{\top }\mathbf{G}_{J^{\ast }}\right) .  \label{eq:YA1}
\end{equation}%
Now let $\tilde{Z}=\Vert \mathbf{a}_{J^{\ast }}\Vert _{\ell _{2}}\leq \rho
_{\min }^{-1/2}\left( \mathbf{G}_{J^{\ast }}^{\top }\mathbf{G}_{J^{\ast
}}\right) \Vert \mathbf{X}\Vert _{\ell _{2}}$. Given our assumptions on $X$
it follows that there exists $\alpha \geq 1$ such that 
\begin{equation}
\Vert \tilde{Z}\Vert _{\psi _{\alpha }}\leq \rho _{\min }^{-1/2}\left( 
\mathbf{G}_{J^{\ast }}^{\top }\mathbf{G}_{J^{\ast }}\right) \Vert Z\Vert
_{\psi _{\alpha }}<+\infty ,  \label{eq:YA2}
\end{equation}%
where $Z=\Vert \mathbf{X}\Vert _{\ell _{2}}$. Hence, using the relations (%
\ref{eq:YA1}) and (\ref{eq:YA2}), and Proposition \ref{prop:mp} (with $%
\mathbf{a}_{J^{\ast }}$ instead of $\mathbf{X}$), it follows that there
exists a universal constant $\delta _{\ast }>0$ such that for all $x>0$, 
\begin{equation*}
\mathbb{P}\left( \left\Vert \widehat{\mathbf{\Psi }}_{J^{\ast }}-\mathbf{%
\Psi }_{J^{\ast }}\right\Vert _{2}\geqslant \tilde{\tau}_{d^{\ast
},N,s_{\ast },1}x\right) \leqslant \exp \left( -(\delta _{\ast }^{-1}x)^{%
\frac{\alpha }{2+\alpha }}\right) ,
\end{equation*}%
where $\tilde{\tau}_{d^{\ast },N,s_{\ast },1}=\max (\tilde{A}_{d^{\ast
},N,s_{\ast },1}^{2},\tilde{B}_{d^{\ast },N,s_{\ast },1})$, with $\tilde{A}%
_{d^{\ast },N,s_{\ast },1}=\Vert \tilde{Z}\Vert _{\psi _{\alpha }}\frac{%
\sqrt{\log d^{\ast }}(\log N)^{1/\alpha }}{\sqrt{N}}$, $\tilde{B}_{d^{\ast
},N,s_{\ast },1}=\frac{\rho ^{2}\left( \mathbf{\Sigma }\right) \rho _{\min
}^{-1}\left( \mathbf{G}_{J^{\ast }}^{\top }\mathbf{G}_{J^{\ast }}\right) }{%
\sqrt{N}}+\left\Vert \mathbf{\Psi }_{J^{\ast }}\right\Vert _{2}^{1/2}\tilde{A%
}_{d^{\ast },N,s_{\ast },1}$ and $d^{\ast }=\min (N,s_{\ast })$. Then, using
the inequality $\Vert \mathbf{S}-\mathbf{\Sigma }\Vert _{2}\leq \rho _{\max
}\left( \mathbf{G}_{J^{\ast }}^{\top }\mathbf{G}_{J^{\ast }}\right) \Vert 
\widehat{\mathbf{\Psi }}_{J^{\ast }}-\mathbf{\Psi }_{J^{\ast }}\Vert _{2}$,
it follows that%
\begin{eqnarray*}
&&\mathbb{P}\left( \Vert \mathbf{S}-\mathbf{\Sigma }\Vert _{2}\geq \rho
_{\max }\left( \mathbf{G}_{J^{\ast }}^{\top }\mathbf{G}_{J^{\ast }}\right) 
\tilde{\tau}_{d^{\ast },N,s_{\ast },1}x\right)  \\
&\leq &\mathbb{P}\left( \rho _{\max }\left( \mathbf{G}_{J^{\ast }}^{\top }%
\mathbf{G}_{J^{\ast }}\right) \left\Vert \widehat{\mathbf{\Psi }}_{J^{\ast
}}-\mathbf{\Psi }_{J^{\ast }}\right\Vert _{2}\geqslant \rho _{\max }\left( 
\mathbf{G}_{J^{\ast }}^{\top }\mathbf{G}_{J^{\ast }}\right) \tilde{\tau}%
_{d^{\ast },N,s_{\ast },1}x\right)  \\
&=&\mathbb{P}\left( \left\Vert \widehat{\mathbf{\Psi }}_{J^{\ast }}-\mathbf{%
\Psi }_{J^{\ast }}\right\Vert _{2}\geqslant \tilde{\tau}_{d^{\ast
},N,s_{\ast },1}x\right)  \\
&\leqslant &\exp \left( -(\delta _{\ast }^{-1}x)^{\frac{\alpha }{2+\alpha }%
}\right) .
\end{eqnarray*}%
Hence, the result follows with 
\begin{eqnarray*}
\tilde{\tau}_{N,s_{\ast }} &=&\rho _{\max }\left( \mathbf{G}%
_{J^{\ast }}^{\top }\mathbf{G}_{J^{\ast }}\right) \tilde{\tau}_{d^{\ast
},N,s_{\ast },1} \\
&=&\max (\rho _{\max }\left( \mathbf{G}_{J^{\ast }}^{\top }\mathbf{G}%
_{J^{\ast }}\right) \tilde{A}_{d^{\ast },N,s_{\ast },1}^{2},\rho _{\max
}\left( \mathbf{G}_{J^{\ast }}^{\top }\mathbf{G}_{J^{\ast }}\right) \tilde{B}%
_{d^{\ast },N,s_{\ast },1}) \\
&=&\max (\tilde{A}_{d^{\ast },N,s_{\ast }}^{2},\tilde{B}_{d^{\ast
},N,s_{\ast }}),
\end{eqnarray*}%
where $\tilde{A}_{d^{\ast },N,s_{\ast }}=$ $\rho _{\max }^{1/2}\left( 
\mathbf{G}_{J^{\ast }}^{\top }\mathbf{G}_{J^{\ast }}\right) \Vert \tilde{Z}%
\Vert _{\psi _{\alpha }}\frac{\sqrt{\log d^{\ast }}(\log N)^{1/\alpha }}{%
\sqrt{N}}$ and, using the inequality 
\begin{equation*}
\left\Vert \mathbf{\Psi }_{J^{\ast }}\right\Vert _{2}=\left\Vert \left( 
\mathbf{G}_{J^{\ast }}^{\top }\mathbf{G}_{J^{\ast }}\right) ^{-1}\mathbf{G}%
_{J^{\ast }}^{\top }\mathbf{\Sigma }\mathbf{G}_{J^{\ast }}\left( \mathbf{G}%
_{J^{\ast }}^{\top }\mathbf{G}_{J^{\ast }}\right) ^{-1}\right\Vert _{2}\leq
\rho _{\min }^{-1}\left( \mathbf{G}_{J^{\ast }}^{\top }\mathbf{G}_{J^{\ast
}}\right) \left\Vert \mathbf{\Sigma }\right\Vert _{2},
\end{equation*}%
$\tilde{B}_{d^{\ast },N,s_{\ast }}=$ $\left( \frac{\rho _{\max }\left( 
\mathbf{G}_{J^{\ast }}^{\top }\mathbf{G}_{J^{\ast }}\right) }{\rho _{\min
}\left( \mathbf{G}_{J^{\ast }}^{\top }\mathbf{G}_{J^{\ast }}\right) }\right) 
\frac{\rho ^{2}\left( \mathbf{\Sigma }\right) }{\sqrt{N}}+\left( \frac{\rho
_{\max }\left( \mathbf{G}_{J^{\ast }}^{\top }\mathbf{G}_{J^{\ast }}\right) }{%
\rho _{\min }\left( \mathbf{G}_{J^{\ast }}^{\top }\mathbf{G}_{J^{\ast
}}\right) }\right) ^{1/2}\left\Vert \mathbf{\Sigma }\right\Vert _{2}^{1/2}%
\tilde{A}_{d^{\ast },N,s_{\ast }}$.


\subsection{ Proof of Theorem \protect\ref{theo:oracle1}}

Let us first prove the following lemmas.

\begin{lemma}
\label{lem:concentration1} Let $\mathcal{E}_{1},...,\mathcal{E}_{N}$ be
independent copies of a second order Gaussian process $\mathcal{E}$ with
zero mean. Let $\mathbf{W}=\frac{1}{N}\sum\limits_{i=1}^{N}\mathbf{W}_{i}$
with 
\begin{equation*}
\mathbf{W}_{i}=\mathcal{E}_{i}\mathcal{E}_{i}^{\top }\in \mathbb{R}^{n\times
n}\mbox{ and }\mathcal{E}_{i}=\left( \mathcal{E}_{i}\left( t_{1}\right) ,...,%
\mathcal{E}_{i}\left( t_{n}\right) \right) ^{\top },\text{ }i=1,\ldots ,N.
\end{equation*}%
Suppose that $\mathbf{\Sigma }_{noise}=\mathbb{E}\left( \mathbf{W}_{1}\right) $ is positive-definite. For $%
1\leq k\leq M$, let $\mathbf{\eta }_{k}$ be the $k$-th column of the matrix $%
\mathbf{G}^{\top }\mathbf{W}\mathbf{G}$. Then, for any $x>0$, 
\begin{equation*}
\mathbb{P}\left( \Vert \mathbf{\eta }_{k}\Vert _{\ell _{2}}\geq \Vert 
\mathbf{G}_{k}\Vert _{\ell _{2}}\sqrt{\rho _{\max }(\mathbf{G}\mathbf{G}%
^{\top })}\Vert \mathbf{\Sigma }_{noise}\Vert _{2}\left( 1+\sqrt{\frac{n}{N}}%
+\sqrt{\frac{2x}{N}}\right) ^{2}\right) \leq \exp (-x).
\end{equation*}
\end{lemma}

\noindent \textbf{Proof of Lemma \ref{lem:concentration1}:} by definition
one has that $\Vert \mathbf{\eta }_{k}\Vert _{\ell _{2}}^{2}=\mathbf{G}%
_{k}^{\top }\mathbf{W}\mathbf{G}\mathbf{G}^{\top }\mathbf{W}\mathbf{G}_{k}$
where $\mathbf{G}_{k}$ denotes the $k$-th column of $\mathbf{G}$. Hence 
\begin{equation}
\Vert \mathbf{\eta }_{k}\Vert _{\ell _{2}}^{2}\leq \Vert \mathbf{G}_{k}\Vert
_{\ell _{2}}^{2}\rho _{\max }(\mathbf{G}\mathbf{G}^{\top })\Vert \mathbf{W}%
\Vert _{2}^{2}. \label{eq:eta}
\end{equation}%
Using the assumption that  $\mathbf{\Sigma }_{noise}$ is positive-definite define the random vectors $Z_{i} = \mathbf{\Sigma }_{noise}^{-1/2} \mathcal{E}_{i},i=1,\ldots,n$. Note that the $Z_{i}$'s are  i.i.d.\  Gaussian vectors in $\R^{n}$ with zero mean and covariance matrix the identity. Then, define  the $N \times n$ matrix
$$
\Gamma = \frac{1}{\sqrt{N}} \left(\begin{array}{c} Z_{1}^{\top}  \\ \vdots  \\  Z_{N}^{\top} \end{array}\right).
$$
Since $\Gamma$ is a matrix with i.i.d.\ entries following a Gaussian distribution with zero mean and variance $1/N$, it follows from  the arguments in the proof of Theorem II.13 in \cite{MR1863696} that for any $x>0$
\begin{equation}
\mathbb{P}\left( \Vert \mathbf{\Gamma^{\top} \Gamma }\Vert _{2}\geq  \left( 1+\sqrt{\frac{n}{N}}+\sqrt{\frac{2x}{N}}\right)
^{2}\right) \leq \exp (-x). \label{eq:deveig}
\end{equation}%
Now, since $\mathbf{W} = \mathbf{\Sigma }_{noise}^{1/2}  \Gamma^{\top} \Gamma \mathbf{\Sigma }_{noise}^{1/2}$ it follows that $ \Vert \mathbf{W} \Vert _{2} \leq \Vert \mathbf{\Sigma }_{noise} \Vert _{2}  \Vert  \Gamma^{\top} \Gamma  \Vert _{2}$. Hence, inequality \eqref{eq:deveig} implies that  for any $x>0$
\begin{equation*}
\mathbb{P}\left( \Vert \mathbf{W}\Vert _{2}\geq \Vert \mathbf{\Sigma }%
_{noise}\Vert _{2}\left( 1+\sqrt{\frac{n}{N}}+\sqrt{\frac{2x}{N}}\right)
^{2}\right) \leq \exp (-x),
\end{equation*}%
and  the result finally follows from inequality \eqref{eq:eta}.
\hfill $\Box $

\begin{lemma}
\label{lem:eig} Let $1\leq s\leq \min (n,M)$ and suppose that Assumption \ref%
{ass:eig} holds for some $c_{0}>0$. Let $J\subset \{1,\ldots ,M\}$ be a
subset of indices of cardinality $|J|\leq s$. Let $\mathbf{\Delta }\in 
\mathcal{S}_{M}$ and suppose that 
\begin{equation*}
\sum_{k\in J^{c}}\Vert \mathbf{\Delta }_{k}\Vert _{\ell _{2}}\leq
c_{0}\sum_{k\in J}\Vert \mathbf{\Delta }_{k}\Vert _{\ell _{2}},
\end{equation*}%
where $\mathbf{\Delta }_{k}$ denotes the $k$-th column of $\mathbf{\Delta }$%
. Let 
\begin{equation*}
\kappa _{s,c_{0}}=\left( \rho _{\min }(s)^{2}-c_{0}\theta (\mathbf{G})\rho
_{\max }(\mathbf{G}^{\top }\mathbf{G})s\right) ^{1/2}.
\end{equation*}%
Then, 
\begin{equation*}
\left\Vert \mathbf{G\mathbf{\Delta }G}^{\top }\right\Vert _{F}^{2}\geq
\kappa _{s,c_{0}}^{2}\left\Vert \mathbf{\Delta }_{J}\right\Vert _{F}^{2},
\end{equation*}%
where $\mathbf{\Delta }_{J}$ denotes the $M\times M$ matrix obtained by
setting to zero the rows and columns of $\mathbf{\Delta }$ whose indices are
not in $J$.
\end{lemma}

\noindent \textbf{Proof of Lemma \ref{lem:eig}:} first let us introduce some
notations. For $\mathbf{\Delta }\in \mathcal{S}_{M}$ and $J\subset
\{1,\ldots ,M\}$, then $\mathbf{\Delta }_{J^{c}}$ denotes the $M\times M$
matrix obtained by setting to zero the rows and columns of $\mathbf{\Delta }$
whose indices are not in the complementary $J^{c}$ of $J$. 
Now, remark that 
\begin{eqnarray}
\left\Vert \mathbf{G\mathbf{\Delta }G}^{\top }\right\Vert _{F}^{2}
&=&\left\Vert \mathbf{G\Delta }_{J}\mathbf{G}^{\top }\right\Vert
_{F}^{2}+\left\Vert \mathbf{G\Delta }_{J^{c}}\mathbf{G}^{\top }\right\Vert
_{F}^{2}+2tr\left( \mathbf{G\Delta }_{J}\mathbf{G}^{\top }\mathbf{G\Delta }%
_{J^{c}}\mathbf{G}^{\top }\right)  \notag \\
&\geq &\left\Vert \mathbf{G\Delta }_{J}\mathbf{G}^{\top }\right\Vert
_{F}^{2}+2tr\left( \mathbf{G\Delta }_{J}\mathbf{G}^{\top }\mathbf{G\Delta }%
_{J^{c}}\mathbf{G}^{\top }\right) .  \label{eq:F}
\end{eqnarray}%
Let $\mathbf{A}=\mathbf{G\Delta }_{J}\mathbf{G}^{\top }$ and $\mathbf{B}=%
\mathbf{G\Delta }_{J^{c}}\mathbf{G}^{\top }$. Using that $tr\left( \mathbf{A}%
^{\top }\mathbf{B}\right) =vec(\mathbf{A})^{\top }vec(\mathbf{B})$ and the
properties (\ref{eq:P1}) and (\ref{eq:P3}) it follows that 
\begin{equation}
tr\left( \mathbf{G\Delta }_{J}\mathbf{G}^{\top }\mathbf{G\Delta }_{J^{c}}%
\mathbf{G}^{\top }\right) =vec(\mathbf{\Delta }_{J})^{\top }\left( \mathbf{G}^{\top }%
\mathbf{G}\otimes \mathbf{G}^{\top }\mathbf{G}\right) vec(\mathbf{\Delta }%
_{J^{c}}).  \label{eq:tr1}
\end{equation}%
Let $\mathbf{C}=\mathbf{G}^{\top }\mathbf{G}\otimes \mathbf{G}^{\top }%
\mathbf{G}$ and note that $\mathbf{C}$ is a $M^{2}\times M^{2}$ matrix whose
elements can be written in the form of $M\times M$ block matrices given by 
\begin{equation*}
\mathbf{C}_{ij}=(\mathbf{G}^{\top }\mathbf{G})_{ij}\mathbf{G}^{\top }\mathbf{%
G},\mbox{ for }1\leq i,j\leq M.
\end{equation*}%
Now, write the $M^{2}\times 1$ vectors $vec(\mathbf{\Delta }_{J})$ and $vec(%
\mathbf{\Delta }_{J^{c}})$ in the form of block vectors as $vec(\mathbf{%
\Delta }_{J})=[(\mathbf{\Delta }_{J})_{i}^{\top }]_{1\leq i\leq M}^{\top }$
and $vec(\mathbf{\Delta }_{J^{c}})=[(\mathbf{\Delta }_{J^{c}})_{j}^{\top
}]_{1\leq j\leq M}^{\top }$, where $(\mathbf{\Delta }_{J})_{i}\in \mathbb{R}%
^{M}$ $(\mathbf{\Delta }_{J^{c}})_{j}\in \mathbb{R}^{M}$ for $1\leq i,j\leq
M $. Using (\ref{eq:tr1}) it follows that 
\begin{eqnarray}
tr\left( \mathbf{G\Delta }_{J}\mathbf{G}^{\top }\mathbf{G\Delta }_{J^{c}}%
\mathbf{G}^{\top }\right) &=&\sum_{1\leq i,j\leq M}(\mathbf{\Delta }%
_{J})_{i}^{\top }\mathbf{C}_{ij}(\mathbf{\Delta }_{J^{c}})_{j}  \notag \\
&=&\sum_{i\in J}\sum_{j\in J^{c}} (\mathbf{G}^{\top }\mathbf{G})_{ij}  (\mathbf{%
\Delta }_{J})_{i}^{\top }\mathbf{G}^{\top }\mathbf{G}(\mathbf{\Delta }%
_{J^{c}})_{j}.  \notag
\end{eqnarray}%
Now, using that $\left\vert (\mathbf{G}^{\top }\mathbf{G})_{ij}\right\vert
\leq \theta (\mathbf{G})$ for $i\neq j$ and that 
\begin{equation*}
\left\vert (\mathbf{\Delta }_{J})_{i}^{\top }\mathbf{G}^{\top }\mathbf{G}(%
\mathbf{\Delta }_{J^{c}})_{j}\right\vert \leq \Vert \mathbf{G}(\mathbf{%
\Delta }_{J})_{i}\Vert _{\ell _{2}}\Vert \mathbf{G}(\mathbf{\Delta }%
_{J^{c}})_{j}\Vert _{\ell _{2}}\leq \rho _{\max }(\mathbf{G}^{\top }\mathbf{G%
})\Vert (\mathbf{\Delta }_{J})_{i}\Vert _{\ell _{2}}\Vert (\mathbf{\Delta }%
_{J^{c}})_{j}\Vert _{\ell _{2}},
\end{equation*}%
it follows that 
\begin{equation*}
tr\left( \mathbf{G\mathbf{\Delta }_{J}G}^{\top }\mathbf{G\mathbf{\Delta }%
_{J^{c}}G}^{\top }\right) \geq -\theta (\mathbf{G})\rho _{\max }(\mathbf{G}%
^{\top }\mathbf{G})\left( \sum_{i\in J}\Vert (\mathbf{\Delta }_{J})_{i}\Vert
_{\ell _{2}}\right) \left( \sum_{j\in J^{c}}\Vert (\mathbf{\Delta }%
_{J^{c}})_{j}\Vert _{\ell _{2}}\right) .
\end{equation*}%
%
%
%
%
%
%
%
%
%
%
%
%
%
%
%
%
%
%
%
%
%
%
Now, using the assumption that $\sum_{k\in J^{c}}\Vert \mathbf{\Delta }%
_{k}\Vert _{\ell _{2}}\leq c_{0}\sum_{k\in J}\Vert \mathbf{\Delta }_{k}\Vert
_{\ell _{2}}$ it follows that 
\begin{eqnarray}
tr\left( \mathbf{G\mathbf{\Delta }_{J}G}^{\top }\mathbf{G\mathbf{\Delta }%
_{J^{c}}G}^{\top }\right) &\geq &-c_{0}\theta (\mathbf{G})\rho _{\max }(%
\mathbf{G}^{\top }\mathbf{G})\left( \sum_{i\in J}\Vert (\mathbf{\Delta }%
_{J})_{i}\Vert _{\ell _{2}}\right) ^{2}  \notag \\
&\geq &-c_{0}\theta (\mathbf{G})\rho _{\max }(\mathbf{G}^{\top }\mathbf{G}%
)s\left\Vert \mathbf{\Delta }_{J}\right\Vert _{F}^{2},  \label{eq:tr2}
\end{eqnarray}%
where, for the inequality, we have used the properties that for the positive
reals $c_{i}=\Vert (\mathbf{\Delta }_{J})_{i}\Vert _{\ell _{2}},\;i\in J$
then $\left( \sum_{i\in J}c_{i}\right) ^{2}\leq |J|\sum_{i\in
J}c_{i}^{2}\leq s\sum_{i\in J}c_{i}^{2}$ and that $\sum_{i\in J}\Vert (%
\mathbf{\Delta }_{J})_{i}\Vert _{\ell _{2}}^{2}=\left\Vert \mathbf{\Delta }%
_{J}\right\Vert _{F}^{2}$. \newline

Using the properties (\ref{eq:P1}) and (\ref{eq:P2}) remark that 
\begin{eqnarray}
\left\Vert \mathbf{G\mathbf{\Delta }_{J}G}^{\top }\right\Vert _{F}^{2}
&=&\Vert \mathbf{G}_{J}\otimes \mathbf{G}_{J}\;vec(\tilde{\mathbf{\Delta }}%
_{J})\Vert _{\ell _{2}}^{2}  \notag \\
&\geq &\rho _{\min }\left( \mathbf{G}_{J}\otimes \mathbf{G}_{J}\right) \Vert
vec(\tilde{\mathbf{\Delta }}_{J})\Vert _{\ell _{2}}^{2}  \notag \\
&\geq &\rho _{\min }(s)^{2}\left\Vert \mathbf{\Delta }_{J}\right\Vert
_{F}^{2},  \label{eq:rho1}
\end{eqnarray}%
where $vec(\tilde{\mathbf{\Delta }}_{J})=[(\mathbf{\Delta }_{J})_{i}^{\top
}]_{i\in J}^{\top }$. Therefore, combining inequalities (\ref{eq:F}), (\ref%
{eq:tr2}) and (\ref{eq:rho1}) it follows that 
\begin{equation*}
\left\Vert \mathbf{G\mathbf{\Delta }G}^{\top }\right\Vert _{F}^{2}\geq
\left( \rho _{\min }(s)^{2}-c_{0}\theta (\mathbf{G})\rho _{\max }(\mathbf{G}%
^{\top }\mathbf{G})s\right) \left\Vert \mathbf{\Delta }_{J}\right\Vert
_{F}^{2},
\end{equation*}%
which completes the proof of Lemma \ref{lem:eig}. \hfill $\Box $ \newline

Let us now proceed to the proof of Theorem \ref{theo:oracle1}. Part of the
proof is inspired by results in \cite{BRT}. Let $s\leq \min (n,M)$ and $%
\mathbf{\Psi }\in \mathcal{S}_{M}$ with $\mathcal{M}(\mathbf{\Psi })\leq s$.
Let $J=\{k\;;\mathbf{\Psi }_{k}\neq 0\}$. To simplify the notations, write $%
\widehat{\mathbf{\Psi }}=\widehat{\mathbf{\Psi }}_{\lambda }$. By definition
of $\widehat{\mathbf{\Sigma }}_{\lambda }=\mathbf{G}\widehat{\mathbf{\Psi }}%
\mathbf{G}^{\top }$ one has that 
\begin{equation}
\left\Vert \widetilde{\mathbf{S}}-\mathbf{G}\widehat{\mathbf{\Psi }}\mathbf{G%
}^{\top }\right\Vert _{F}^{2}+2\lambda \sum_{k=1}^{M}\gamma _{k}\Vert 
\widehat{\mathbf{\Psi }}_{k}\Vert _{\ell _{2}}\leq \left\Vert \widetilde{%
\mathbf{S}}-\mathbf{G\mathbf{\Psi }G}^{\top }\right\Vert _{F}^{2}+2\lambda
\sum_{k=1}^{M}\gamma _{k}\Vert \mathbf{\Psi }_{k}\Vert _{\ell _{2}}.
\label{Ineq1}
\end{equation}%
Using the scalar product associated to the Frobenius norm $\left\langle
A,B\right\rangle _{F}=tr\left( A^{\top }B\right) $ then 
\begin{eqnarray}
\left\Vert \widetilde{\mathbf{S}}-\mathbf{G}\widehat{\mathbf{\Psi }}\mathbf{G%
}^{\top }\right\Vert _{F}^{2} &=&\left\Vert \mathbf{S}+\mathbf{W}-\mathbf{G}%
\widehat{\mathbf{\Psi }}\mathbf{G}^{\top }\right\Vert _{F}^{2}  \notag \\
&=&\left\Vert \mathbf{W}\right\Vert _{F}^{2}+\left\Vert \mathbf{S}-\mathbf{G}%
\widehat{\mathbf{\Psi }}\mathbf{G}^{\top }\right\Vert _{F}^{2}+2\left\langle 
\mathbf{W},\mathbf{S}-\mathbf{G}\widehat{\mathbf{\Psi }}\mathbf{G}^{\top
}\right\rangle _{F}.  \label{T1}
\end{eqnarray}%
Putting (\ref{T1}) in (\ref{Ineq1}) we get%
\begin{eqnarray*}
\left\Vert \mathbf{S}-\mathbf{G}\widehat{\mathbf{\Psi }}\mathbf{G}^{\top
}\right\Vert _{F}^{2}+2\lambda \sum_{k=1}^{M}\gamma _{k}\Vert \widehat{%
\mathbf{\Psi }}_{k}\Vert _{\ell _{2}} &\leq &\left\Vert \mathbf{S}-\mathbf{G%
\mathbf{\Psi }G}^{\top }\right\Vert _{F}^{2}+2\left\langle \mathbf{W},%
\mathbf{G}\left( \widehat{\mathbf{\Psi }}-\mathbf{\Psi }\right) \mathbf{G}%
^{\top }\right\rangle _{F} \\
&&+2\lambda \sum_{k=1}^{M}\gamma _{k}\Vert \mathbf{\Psi }_{k}\Vert _{\ell
_{2}}.
\end{eqnarray*}

For $k=1,\ldots ,M$ define the $M\times M$ matrix $\mathbf{A}_{k}$ with all
columns equal to zero except the $k$-th which is equal to $\widehat{\mathbf{%
\Psi }}_{k}-\mathbf{\Psi }_{k}$. Then, remark that 
\begin{eqnarray*}
\left\langle \mathbf{W},\mathbf{G}\left( \widehat{\mathbf{\Psi }}-\mathbf{%
\Psi }\right) \mathbf{G}^{\top }\right\rangle _{F}
&=&\sum_{k=1}^{M}\left\langle \mathbf{W},\mathbf{G}\mathbf{A}_{k}\mathbf{G}%
^{\top }\right\rangle _{F}=\sum_{k=1}^{M}\left\langle \mathbf{G}^{\top }%
\mathbf{W}\mathbf{G},\mathbf{A}_{k}\right\rangle _{F}=\sum_{k=1}^{M}\mathbf{%
\eta }_{k}^{\top }(\widehat{\mathbf{\Psi }}_{k}-\mathbf{\Psi }_{k}) \\
&\leq &\sum_{k=1}^{M}\Vert \mathbf{\eta }_{k}\Vert _{\ell _{2}}\Vert 
\widehat{\mathbf{\Psi }}_{k}-\mathbf{\Psi }_{k}\Vert _{\ell _{2}},
\end{eqnarray*}%
where $\mathbf{\eta }_{k}$ is the $k$-th column of the matrix $\mathbf{G}%
^{\top }\mathbf{W}\mathbf{G}$. Define the event 
\begin{equation}
\mathcal{A}=\bigcap_{k=1}^{M}\left\{ 2\Vert \mathbf{\eta }_{k}\Vert _{\ell
_{2}}\leq \lambda \gamma _{k}\right\} .  \label{Event A}
\end{equation}%
Then, the choices 
\begin{equation*}
\gamma _{k}=2\Vert \mathbf{G}_{k}\Vert _{\ell _{2}}\sqrt{\rho _{\max }(%
\mathbf{G}\mathbf{G}^{\top })},\;\lambda =\Vert \mathbf{\Sigma }%
_{noise}\Vert _{2}\left( 1+\sqrt{\frac{n}{N}}+\sqrt{\frac{2\delta \log M}{N}}%
\right) ^{2},
\end{equation*}%
and Lemma \ref{lem:concentration1} imply that the probability of the
complementary event $\mathcal{A}^{c}$ satisfies 
\begin{equation*}
\mathbb{P}\left( \mathcal{A}^{c}\right) \leq \sum_{k=1}^{M}\mathbb{P}\left(
2\Vert \mathbf{\eta }_{k}\Vert _{\ell _{2}}>\lambda \gamma _{k}\right) \leq
M^{1-\delta }.
\end{equation*}%
Then, on the event $\mathcal{A}$ one has that 
\begin{eqnarray*}
\left\Vert \mathbf{S}-\mathbf{G}\widehat{\mathbf{\Psi }}\mathbf{G}^{\top
}\right\Vert _{F}^{2} &\leq &\left\Vert \mathbf{S}-\mathbf{G\mathbf{\Psi }G}%
^{\top }\right\Vert _{F}^{2}+\lambda \sum_{k=1}^{M}\ \gamma _{k}\Vert 
\widehat{\mathbf{\Psi }}_{k}-\mathbf{\Psi }_{k}\Vert _{\ell _{2}} \\
&&+2\lambda \sum_{k=1}^{M}\gamma _{k}\left( \Vert \mathbf{\Psi }_{k}\Vert
_{\ell _{2}}-\Vert \widehat{\mathbf{\Psi }}_{k}\Vert _{\ell _{2}}\right) .
\end{eqnarray*}%
Adding the term $\lambda \sum_{k=1}^{M}\ \gamma _{k}\Vert \widehat{\mathbf{%
\Psi }}_{k}-\mathbf{\Psi }_{k}\Vert _{\ell _{2}}$ to both sides of the above
inequality yields on the event $\mathcal{A}$ 
\begin{eqnarray*}
\left\Vert \mathbf{S}-\mathbf{G}\widehat{\mathbf{\Psi }}\mathbf{G}^{\top
}\right\Vert _{F}^{2}+\lambda \sum_{k=1}^{M}\ \gamma _{k}\Vert \widehat{%
\mathbf{\Psi }}_{k}-\mathbf{\Psi }_{k}\Vert _{\ell _{2}} &\leq &\left\Vert 
\mathbf{S}-\mathbf{G\mathbf{\Psi }G}^{\top }\right\Vert _{F}^{2} \\
&&+2\lambda \sum_{k=1}^{M}\gamma _{k}\left( \Vert \widehat{\mathbf{\Psi }}%
_{k}-\mathbf{\Psi }_{k}\Vert _{\ell _{2}}+\Vert \mathbf{\Psi }_{k}\Vert
_{\ell _{2}}-\Vert \widehat{\mathbf{\Psi }}_{k}\Vert _{\ell _{2}}\right) .
\end{eqnarray*}%
Now, remark that for all $k\notin J$, then $\Vert \widehat{\mathbf{\Psi }}%
_{k}-\mathbf{\Psi }_{k}\Vert _{\ell _{2}}+\Vert \mathbf{\Psi }_{k}\Vert
_{\ell _{2}}-\Vert \widehat{\mathbf{\Psi }}_{k}\Vert _{\ell _{2}}=0$, which
implies that  on the event $\mathcal{A}$ 
\begin{eqnarray}
\left\Vert \mathbf{S}-\mathbf{G}\widehat{\mathbf{\Psi }}\mathbf{G}^{\top
}\right\Vert _{F}^{2}+\lambda \sum_{k=1}^{M}\ \gamma _{k}\Vert \widehat{%
\mathbf{\Psi }}_{k}-\mathbf{\Psi }_{k}\Vert _{\ell _{2}} &\leq &\left\Vert 
\mathbf{S}-\mathbf{G\mathbf{\Psi }G}^{\top }\right\Vert _{F}^{2}   \label{eq:B1}  \\
&&+4\lambda \sum_{k\in J}\gamma _{k}\Vert \widehat{\mathbf{\Psi }}_{k}-%
\mathbf{\Psi }_{k}\Vert _{\ell _{2}}  \notag \\
&\leq &\left\Vert \mathbf{S}-\mathbf{G\mathbf{\Psi }G}^{\top }\right\Vert
_{F}^{2}   \label{eq:B2} \\
&&+4\lambda \sqrt{\mathcal{M}(\mathbf{\Psi })}\sqrt{\sum_{k\in J}\gamma
_{k}^{2}\Vert \widehat{\mathbf{\Psi }}_{k}-\mathbf{\Psi }_{k}\Vert _{\ell
_{2}}^{2}}.  \notag 
\end{eqnarray}%
where for the last inequality we have used the property that for the
positive reals $c_{k}=\gamma _{k}\Vert \widehat{\mathbf{\Psi }}_{k}-\mathbf{%
\Psi }_{k}\Vert _{\ell _{2}},\;k\in J$ then $\left( \sum_{k\in
J}c_{k}\right) ^{2}\leq \mathcal{M}(\mathbf{\Psi })\sum_{k\in J}c_{k}^{2}$.

Let $\epsilon >0$ and define the event
\begin{equation}
\mathcal{A}_{1}=\left\{ 4\lambda
\sum_{k\in J}\gamma _{k}\Vert \widehat{\mathbf{\Psi }}_{k}-\mathbf{\Psi }%
_{k}\Vert _{\ell _{2}}>\epsilon \left\Vert \mathbf{S}-\mathbf{G\mathbf{\Psi }%
G}^{\top }\right\Vert _{F}^{2}\right\} . \label{eq:eventA1}
\end{equation}
Note that on the event $\mathcal{A}%
\cap \mathcal{A}_{1}^{c}$ then the result of the theorem trivially follows
from inequality (\ref{eq:B1}). Now consider the event $\mathcal{A}\cap 
\mathcal{A}_{1}$ (all the following inequalities hold on this event). Using (%
\ref{eq:B1}) one has that 
\begin{equation}
\lambda \sum_{k=1}^{M}\ \gamma _{k}\Vert \widehat{\mathbf{\Psi }}_{k}-%
\mathbf{\Psi }_{k}\Vert _{\ell _{2}}\leq 4(1+1/\epsilon )\lambda \sum_{k\in
J}\gamma _{k}\Vert \widehat{\mathbf{\Psi }}_{k}-\mathbf{\Psi }_{k}\Vert
_{\ell _{2}}.  \label{Ineq5}
\end{equation}%
Therefore, on $\mathcal{A}\cap \mathcal{A}_{1}$ 
\begin{equation*}
\sum_{k\notin J}\ \gamma _{k}\Vert \widehat{\mathbf{\Psi }}_{k}-\mathbf{\Psi 
}_{k}\Vert _{\ell _{2}}\leq (3+4/\epsilon )\sum_{k\in J}\gamma _{k}\Vert 
\widehat{\mathbf{\Psi }}_{k}-\mathbf{\Psi }_{k}\Vert _{\ell _{2}}.
\end{equation*}%
Let $\mathbf{\Delta }$ be the $M\times M$ symmetric matrix with columns
equal to $\mathbf{\Delta }_{k}=\gamma _{k}\left( \widehat{\mathbf{\Psi }}%
_{k}-\mathbf{\Psi }_{k}\right) ,k=1,\ldots ,M$, and $c_{0}=3+4/\epsilon $.
Then, the above inequality means that $\sum_{k\in J^{c}}\Vert \mathbf{\Delta 
}_{k}\Vert _{\ell _{2}}\leq c_{0}\sum_{k\in J}\Vert \mathbf{\Delta }%
_{k}\Vert _{\ell _{2}}$ and thus Assumption \ref{ass:eig} and Lemma \ref%
{lem:eig} imply that 
\begin{equation}
\kappa _{s,c_{0}}^{2}\sum_{k\in J}\gamma _{k}^{2}\Vert \widehat{\mathbf{\Psi 
}}_{k}-\mathbf{\Psi }_{k}\Vert _{\ell _{2}}^{2}\leq \left\Vert \mathbf{G%
\mathbf{\Delta }G}^{\top }\right\Vert _{F}^{2}\leq 4\mathbf{G}_{\max
}^{2}\rho _{\max }(\mathbf{G}^{\top }\mathbf{G})\left\Vert \mathbf{G(%
\widehat{\mathbf{\Psi }}-\mathbf{\Psi })\mathbf{G}^{\top }}\right\Vert
_{F}^{2}.  \label{Ineq4}
\end{equation}%
Let $\gamma _{\max }^{2}=4\mathbf{G}_{\max }^{2}\rho _{\max }(\mathbf{G}%
^{\top }\mathbf{G})$. Combining the above inequality with (\ref{eq:B2})
yields 
\begin{eqnarray*}
\left\Vert \mathbf{S}-\mathbf{G}\widehat{\mathbf{\Psi }}\mathbf{G}^{\top
}\right\Vert _{F}^{2} &\leq &\left\Vert \mathbf{S}-\mathbf{G\mathbf{\Psi }G}%
^{\top }\right\Vert _{F}^{2}+4\lambda \kappa _{s,c_{0}}^{-1}\gamma _{\max }%
\sqrt{\mathcal{M}(\mathbf{\Psi })}\left\Vert \mathbf{G(\widehat{\mathbf{\Psi 
}}-\mathbf{\Psi })\mathbf{G}^{\top }}\right\Vert _{F} \\
&\leq &\left\Vert \mathbf{S}-\mathbf{G\mathbf{\Psi }G}^{\top }\right\Vert
_{F}^{2}+4\lambda \kappa _{s,c_{0}}^{-1}\gamma _{\max }\sqrt{\mathcal{M}(%
\mathbf{\Psi })}\left( \left\Vert \mathbf{G\widehat{\mathbf{\Psi }}\mathbf{G}%
^{\top }}-\mathbf{S}\right\Vert _{F}\right. \\
&&+\left. \left\Vert \mathbf{G\mathbf{\Psi }G}^{\top }-\mathbf{S}\right\Vert
_{F}\right)
\end{eqnarray*}%
Now, arguing as in \cite{BRT}, a decoupling argument using the inequality $%
2xy\leq bx^{2}+b^{-1}y^{2}$ with $b>1$, $x=2\lambda \kappa
_{s,c_{0}}^{-1}\gamma _{\max }\sqrt{\mathcal{M}(\mathbf{\Psi })}$ and $y$
being either $\left\Vert \mathbf{G\widehat{\mathbf{\Psi }}\mathbf{G}^{\top }}%
-\mathbf{S}\right\Vert _{F}$ or $\left\Vert \mathbf{G\mathbf{\Psi }G}^{\top
}-\mathbf{S}\right\Vert _{F}$ yields the inequality 
\begin{equation}
\left\Vert \mathbf{S}-\mathbf{G}\widehat{\mathbf{\Psi }}\mathbf{G}^{\top
}\right\Vert _{F}^{2}\leq \left( \frac{b+1}{b-1}\right) \left\Vert \mathbf{S}%
-\mathbf{G\mathbf{\Psi }G}^{\top }\right\Vert _{F}^{2}+\frac{8b^{2}\gamma
_{\max }^{2}}{(b-1)\kappa _{s,c_{0}}^{2}}\lambda ^{2}\mathcal{M}(\mathbf{%
\Psi }).  \label{Ineq7}
\end{equation}%
Then, taking $b=1+2/\epsilon $ and using the inequalities $\left\Vert 
\mathbf{\Sigma }-\mathbf{G}\widehat{\mathbf{\Psi }}\mathbf{G}^{\top
}\right\Vert _{F}^{2}\leq 2\left\Vert \mathbf{S}-\mathbf{\Sigma }\right\Vert
_{F}^{2}+2\left\Vert \mathbf{S}-\mathbf{G}\widehat{\mathbf{\Psi }}\mathbf{G}%
^{\top }\right\Vert _{F}^{2}$ and $\left\Vert \mathbf{S}-\mathbf{G}{\mathbf{%
\Psi }}\mathbf{G}^{\top }\right\Vert _{F}^{2}\leq 2\left\Vert \mathbf{S}-%
\mathbf{\Sigma }\right\Vert _{F}^{2}+2\left\Vert \mathbf{\Sigma }-\mathbf{G}{%
\mathbf{\Psi }}\mathbf{G}^{\top }\right\Vert _{F}^{2}$ completes the proof
of Theorem \ref{theo:oracle1}. \hfill $\Box $

\subsection{Proof of Theorem \protect\ref{theo:activeset}}

Part of the proof is inspired by the approach followed in \cite{MR2386087}
and \cite{LPTG}. Note first that 
\begin{equation*}
\underset{1\leq k\leq M}{\max }\gamma _{k}\left\Vert \widehat{\mathbf{%
\mathbf{\Psi }}}_{k}-\mathbf{\mathbf{\Psi }}_{k}^{\ast }\right\Vert _{\ell
_{2}}\leq \sum_{k=1}^{M}\gamma _{k}\left\Vert \widehat{\mathbf{\mathbf{\Psi }%
}}_{k}-\mathbf{\mathbf{\Psi }}_{k}^{\ast }\right\Vert _{\ell _{2}}.
\end{equation*}%
Since $\mathbf{\mathbf{\Psi }}^{\ast }\in \left\{ \mathbf{\mathbf{\Psi }}\in 
\mathcal{S}_{M}:M\left( \mathbf{\mathbf{\Psi }}\right) \leq s_{\ast
}\right\} $, we can use some results from the proof of Theorem (\ref%
{theo:oracle1}).
On the event $\mathcal{A} \cap \mathcal{A}_{1}$, with $\mathcal{A}$ defined by (\ref{Event A}) and $\mathcal{A}_{1}$ defined by (\ref{eq:eventA1}), inequality (\ref{Ineq5}) implies that 
\begin{eqnarray*}
\sum_{k=1}^{M}\gamma _{k}\left\Vert \widehat{\mathbf{\mathbf{\Psi }}}_{k}-%
\mathbf{\mathbf{\Psi }}_{k}^{\ast }\right\Vert _{\ell _{2}} &\leq &4\left( 1+%
\frac{1}{\epsilon }\right) \sum_{k\in J^{\ast }}\gamma _{k}\left\Vert 
\widehat{\mathbf{\mathbf{\Psi }}}_{k}-\mathbf{\mathbf{\Psi }}_{k}^{\ast
}\right\Vert _{\ell _{2}} \\
&\leq &4\left( 1+\frac{1}{\epsilon }\right) \sqrt{s_{\ast }}\sqrt{\sum_{k\in
J^{\ast }}\gamma _{k}^{2}\left\Vert \widehat{\mathbf{\mathbf{\Psi }}}_{k}-%
\mathbf{\mathbf{\Psi }}_{k}^{\ast }\right\Vert _{\ell _{2}}^{2}}.
\end{eqnarray*}%
Let $\mathbf{\Delta }^{\ast }$ be the $M\times M$ symmetric matrix with
columns equal to $\mathbf{\Delta }_{k}^{\ast }=\gamma _{k}\left( \widehat{%
\mathbf{\Psi }}_{k}-\mathbf{\Psi }_{k}^{\ast }\right) $, $k=1,\ldots ,M$,
let $\gamma _{\max }=2\mathbf{G}_{\max }\sqrt{\rho _{\max }(\mathbf{G}^{\top
}\mathbf{G})}$ and $c_{0}=3+4/\epsilon $. Then, the above inequality and (%
\ref{Ineq4}) imply that on the event $\mathcal{A}  \cap \mathcal{A}_{1}$

\begin{eqnarray*}
\sum_{k=1}^{M}\gamma _{k}\left\Vert \widehat{\mathbf{\mathbf{\Psi }}}_{k}-%
\mathbf{\mathbf{\Psi }}_{k}^{\ast }\right\Vert _{\ell _{2}} &\leq &\frac{%
4\left( 1+\frac{1}{\epsilon }\right) \sqrt{s_{\ast }}}{\kappa _{s_{\ast
},c_{0}}}\left\Vert \mathbf{G\mathbf{\Delta }}^{\ast }\mathbf{G}^{\top
}\right\Vert _{F}\leq \frac{4\left( 1+\frac{1}{\epsilon }\right) \sqrt{%
s_{\ast }}}{\kappa _{s_{\ast },c_{0}}}\gamma _{\max }\left\Vert \mathbf{G}%
\left( \mathbf{\widehat{\mathbf{\Psi }}-\mathbf{\mathbf{\Psi }}^{\ast }}%
\right) \mathbf{G}^{\top }\right\Vert _{F} \\
&= &\frac{4\left( 1+\epsilon \right) \sqrt{s_{\ast }}}{\epsilon \kappa
_{s_{\ast },c_{0}}}\gamma _{\max }\left\Vert \widehat{\mathbf{\Sigma }}%
_{\lambda }-\mathbf{\Sigma }\right\Vert _{F} \\
&\leq &\frac{4\left( 1+\epsilon \right) \sqrt{s_{\ast }}}{\epsilon \kappa
_{s_{\ast },c_{0}}}\gamma _{\max }\sqrt{n}\sqrt{C_{0}\left( n,M,N,s_{\ast },%
\mathbf{S,}\mathbf{\Psi }^{\ast },\mathbf{G},\mathbf{\Sigma }_{noise}\right) 
},
\end{eqnarray*}%
Then, using (\ref{eq:B1}) one has that on the event $\mathcal{A} \cap \mathcal{A}_{1}^{c}$
$$
\sum_{k=1}^{M}\gamma _{k}\left\Vert \widehat{\mathbf{\mathbf{\Psi }}}_{k}-%
\mathbf{\mathbf{\Psi }}_{k}^{\ast }\right\Vert _{\ell _{2}} \leq \frac{1+\epsilon}{\lambda}  \left\Vert \mathbf{S}-\mathbf{G\mathbf{\Psi }^{\ast }%
G}^{\top }\right\Vert _{F}^{2}.
$$
Therefore, by definition of $C_{1}$, the previous inequalities imply that on the event  $\mathcal{A}$ (of probability $1-M^{1-\delta }$ )
\begin{equation}
\sum_{k=1}^{M}\frac{\Vert \mathbf{G}_{k}\Vert _{\ell _{2}}}{\sqrt{n}\mathbf{G%
}_{\max }}\left\Vert \widehat{\mathbf{\mathbf{\Psi }}}_{k}-\mathbf{\mathbf{%
\Psi }}_{k}^{\ast }\right\Vert _{\ell _{2}}\leq C_{1}\left( n,M,N,s_{\ast },%
\mathbf{S,}\mathbf{\Psi }^{\ast },\mathbf{G},\mathbf{\Sigma }_{noise}\right)
.  \label{Ineq2bis}
\end{equation}%
%
%
%
%
%
%
%
%
%
%
%
%
%
%
%
%
%
%
%
%
%
%
%
%
%
%
Hence $\underset{1\leq k\leq M}{\max }\frac{\delta _{k}}{\sqrt{n}}\left\Vert 
\widehat{\mathbf{\mathbf{\Psi }}}_{k}-\mathbf{\mathbf{\Psi }}_{k}^{\ast
}\right\Vert _{\ell _{2}}\leq C_{1}\left( \sigma ,n,M,N,s_{\ast },\mathbf{G},%
\mathbf{\Sigma }_{noise}\right) $ with probability at least $1-M^{1-\delta }$%
, which proves the first assertion of Theorem \ref{theo:activeset}.

Then, to prove that $\hat{J}=J^{\ast }$ we use that $\frac{\delta _{k}}{%
\sqrt{n}}\left\vert \left\Vert \widehat{\mathbf{\mathbf{\Psi }}}%
_{k}\right\Vert _{\ell _{2}}-\left\Vert \mathbf{\mathbf{\Psi }}_{k}^{\ast
}\right\Vert _{\ell _{2}}\right\vert \leq \frac{\delta _{k}}{\sqrt{n}}%
\left\Vert \widehat{\mathbf{\mathbf{\Psi }}}_{k}-\mathbf{\mathbf{\Psi }}%
_{k}^{\ast }\right\Vert _{\ell _{2}}$ for all $k=1,\ldots ,M$. Then, by (\ref%
{Ineq2bis})%
\begin{equation*}
\left\vert \frac{\delta _{k}}{\sqrt{n}}\left\Vert \widehat{\mathbf{\mathbf{%
\Psi }}}_{k}\right\Vert _{\ell _{2}}-\frac{\delta _{k}}{\sqrt{n}}\left\Vert 
\mathbf{\mathbf{\Psi }}_{k}^{\ast }\right\Vert _{\ell _{2}}\right\vert \leq
C_{1}\left( n,M,N,s_{\ast },\mathbf{S,}\mathbf{\Psi }^{\ast },\mathbf{G},%
\mathbf{\Sigma }_{noise}\right) ,
\end{equation*}%
which is equivalent to%
\begin{equation}
-C_{1}\left( n,M,N,s_{\ast },\mathbf{S,}\mathbf{\Psi }^{\ast },\mathbf{G},%
\mathbf{\Sigma }_{noise}\right) \leq \frac{\delta _{k}}{\sqrt{n}}\left\Vert 
\widehat{\mathbf{\mathbf{\Psi }}}_{k}\right\Vert _{\ell _{2}}-\frac{\delta
_{k}}{\sqrt{n}}\left\Vert \mathbf{\mathbf{\Psi }}_{k}^{\ast }\right\Vert
_{\ell _{2}}\leq C_{1}\left( n,M,N,s_{\ast },\mathbf{S,}\mathbf{\Psi }^{\ast
},\mathbf{G},\mathbf{\Sigma }_{noise}\right) .  \label{Ineq3}
\end{equation}%
If $k\in $ $\hat{J}$ then $\frac{\delta _{k}}{\sqrt{n}}\left\Vert \widehat{%
\mathbf{\mathbf{\Psi }}}_{k}\right\Vert _{\ell _{2}}>C_{1}\left(
n,M,N,s_{\ast },\mathbf{S,}\mathbf{\Psi }^{\ast },\mathbf{G},\mathbf{\Sigma }%
_{noise}\right) $. Inequality $\frac{\delta _{k}}{\sqrt{n}}\left\Vert 
\widehat{\mathbf{\mathbf{\Psi }}}_{k}\right\Vert _{\ell _{2}}-\frac{\delta
_{k}}{\sqrt{n}}\left\Vert \mathbf{\mathbf{\Psi }}_{k}^{\ast }\right\Vert
_{\ell _{2}}\leq C_{1}\left( n,M,N,s_{\ast },\mathbf{S,}\mathbf{\Psi }^{\ast
},\mathbf{G},\mathbf{\Sigma }_{noise}\right) $ from (\ref{Ineq3}) imply that 
$\frac{\delta _{k}}{\sqrt{n}}\left\Vert \mathbf{\mathbf{\Psi }}_{k}^{\ast
}\right\Vert _{\ell _{2}}\geq \frac{\delta _{k}}{\sqrt{n}}\left\Vert 
\widehat{\mathbf{\mathbf{\Psi }}}_{k}\right\Vert _{\ell _{2}}-C_{1}\left(
n,M,N,s_{\ast },\mathbf{S,}\mathbf{\Psi }^{\ast },\mathbf{G},\mathbf{\Sigma }%
_{noise}\right) >0$, where the last inequality is obtained using that $k\in 
\hat{J}$. Hence $\left\Vert \mathbf{\mathbf{\Psi }}_{k}^{\ast }\right\Vert
_{\ell _{2}}>0$ and therefore $k\in J^{\ast }$. If $k\in J^{\ast }$ then $%
\left\Vert \mathbf{\mathbf{\Psi }}_{k}^{\ast }\right\Vert _{\ell _{2}}\neq 0$%
. Inequality $-C_{1}\left( n,M,N,s_{\ast },\mathbf{S,}\mathbf{\Psi }^{\ast },%
\mathbf{G},\mathbf{\Sigma }_{noise}\right) \leq \frac{\delta _{k}}{\sqrt{n}}%
\left\Vert \widehat{\mathbf{\mathbf{\Psi }}}_{k}\right\Vert _{\ell _{2}}-%
\frac{\delta _{k}}{\sqrt{n}}\left\Vert \mathbf{\mathbf{\Psi }}_{k}^{\ast
}\right\Vert _{\ell _{2}}$ from (\ref{Ineq3}) imply that $\frac{\delta _{k}}{%
\sqrt{n}}\left\Vert \widehat{\mathbf{\mathbf{\Psi }}}_{k}\right\Vert _{\ell
_{2}}+C_{1}\left( n,M,N,s_{\ast },\mathbf{S,}\mathbf{\Psi }^{\ast },\mathbf{G%
},\mathbf{\Sigma }_{noise}\right) \geq \frac{\delta _{k}}{\sqrt{n}}%
\left\Vert \mathbf{\mathbf{\Psi }}_{k}^{\ast }\right\Vert _{\ell
_{2}}>2C_{1}\left( n,M,N,s_{\ast },\mathbf{S,}\mathbf{\Psi }^{\ast },\mathbf{%
G},\mathbf{\Sigma }_{noise}\right) $, where the last inequality is obtained
using Assumption (\ref{Hyp1}) on $\frac{\delta _{k}}{\sqrt{n}}\left\Vert 
\mathbf{\mathbf{\Psi }}_{k}^{\ast }\right\Vert _{\ell _{2}}$. Hence $\frac{%
\delta _{k}}{\sqrt{n}}\left\Vert \widehat{\mathbf{\mathbf{\Psi }}}%
_{k}\right\Vert _{\ell _{2}}>2C_{1}\left( n,M,N,s_{\ast },\mathbf{S,}\mathbf{%
\Psi }^{\ast },\mathbf{G},\mathbf{\Sigma }_{noise}\right) -C_{1}\left(
n,M,N,s_{\ast },\mathbf{S,}\mathbf{\Psi }^{\ast },\mathbf{G},\mathbf{\Sigma }%
_{noise}\right) =C_{1}\left( n,M,N,s_{\ast },\mathbf{S,}\mathbf{\Psi }^{\ast
},\mathbf{G},\mathbf{\Sigma }_{noise}\right) $ and therefore $k\in $ $\hat{J}
$. This completes the proof of Theorem \ref{theo:activeset}. \hfill $\Box $

\subsection{Proof of Theorem \protect\ref{theo:opnorm}}

Under the assumptions of Theorem \ref{theo:opnorm}, we have shown in the
proof of Theorem \ref{theo:activeset} that $\hat{J}=J^{\ast }$ on the event $%
\mathcal{A}$ defined by (\ref{Event A}). Therefore, under the assumptions of
Theorem \ref{theo:opnorm} it can be checked that on the event $\mathcal{A}$
(of probability $1-M^{1-\delta }$) 
\begin{equation*}
\widehat{\mathbf{\Sigma }}_{\hat{J}}=\widehat{\mathbf{\Sigma }}_{J^{\ast }}=%
\mathbf{G}_{J^{\ast }}\widehat{\mathbf{\Psi }}_{J^{\ast }}\mathbf{G}%
_{J^{\ast }}^{\top },
\end{equation*}%
with 
\begin{equation*}
\widehat{\mathbf{\Psi }}_{J^{\ast }}=\left( \mathbf{G}_{J^{\ast }}^{\top }%
\mathbf{G}_{J^{\ast }}\right) ^{-1}\mathbf{G}_{J^{\ast }}^{\top }\widetilde{%
\mathbf{S}}\mathbf{G}_{J^{\ast }}\left( \mathbf{G}_{J^{\ast }}^{\top }%
\mathbf{G}_{J^{\ast }}\right) ^{-1}.
\end{equation*}%
Now, from the definition (\ref{eq:SigmaJast}) of $\mathbf{\Sigma }_{J^{\ast
}}$ it follows that on the event $\mathcal{A}$ 
\begin{equation}
\left\Vert \widehat{\mathbf{\Sigma }}_{\hat{J}}-\mathbf{\Sigma }_{J^{\ast
}}\right\Vert _{2}\leq \rho _{\max }\left( \mathbf{G}_{J^{\ast }}^{\top }%
\mathbf{G}_{J^{\ast }}\right) \left\Vert \widehat{\mathbf{\Psi }}_{J^{\ast
}}-\Lambda _{J^{\ast }}\right\Vert _{2}  \label{eq:ineqSigma}
\end{equation}%
where $\Lambda _{J^{\ast }}=\mathbf{\Psi }_{J^{\ast }}+(\mathbf{G}_{J^{\ast
}}^{\top }\mathbf{G}_{J^{\ast }})^{-1}\mathbf{G}_{J^{\ast }}^{\top }\mathbf{%
\Sigma }_{noise}\mathbf{G}_{J^{\ast }}\left( \mathbf{G}_{J^{\ast }}^{\top }%
\mathbf{G}_{J^{\ast }}\right) ^{-1}$. Let $\mathbf{Y}_{i}=\left( \mathbf{G}%
_{J^{\ast }}^{\top }\mathbf{G}_{J^{\ast }}\right) ^{-1}\mathbf{G}_{J^{\ast
}}^{\top }\widetilde{\mathbf{X}}_{i}$ for $i=1,\ldots ,N$ and remark that 
\begin{equation*}
\widehat{\mathbf{\Psi }}_{J^{\ast }}=\frac{1}{N}\sum\limits_{i=1}^{N}\mathbf{%
Y}_{i}\mathbf{Y}_{i}^{\top }\mbox{ with }\mathbb{E}\widehat{\mathbf{\Psi }}%
_{J^{\ast }}=\Lambda _{J^{\ast }}.
\end{equation*}%
Therefore, $\widehat{\mathbf{\Psi }}_{J^{\ast }}$ is a sample covariance
matrix of size $s_{\ast }\times s_{\ast }$ and we can control its deviation
in operator norm from $\Lambda _{J^{\ast }}$ by using Proposition \ref%
{prop:mp}. For this we simply have to verify conditions similar to \textbf{%
(A1)} and \textbf{(A2)} in Assumption \ref{ass:X} for the random vector $%
\mathbf{Y}=\left( \mathbf{G}_{J^{\ast }}^{\top }\mathbf{G}_{J^{\ast
}}\right) ^{-1}\mathbf{G}_{J^{\ast }}^{\top }\widetilde{\mathbf{X}}\in 
\mathbb{R}^{s_{\ast }}$. First, let $\beta \in \mathbb{R}^{s_{\ast }}$ with $%
\Vert \beta \Vert _{\ell _{2}}=1$. Then, remark that $\mathbf{Y}^{\top
}\beta =\widetilde{\mathbf{X}}^{\top }\tilde{\beta}$ with $\tilde{\beta}=%
\mathbf{G}_{J^{\ast }}\left( \mathbf{G}_{J^{\ast }}^{\top }\mathbf{G}%
_{J^{\ast }}\right) ^{-1}\beta $. Since $\Vert \widetilde{\beta }\Vert
_{\ell _{2}}\leq \left( \rho _{\min }\left( \mathbf{G}_{J^{\ast }}^{\top }%
\mathbf{G}_{J^{\ast }}\right) \right) ^{-1/2}$ it follows that 
\begin{equation}
\left( \mathbb{E}|\mathbf{Y}^{\top }\beta |^{4}\right) ^{1/4}\leq \tilde{\rho%
}(\mathbf{\Sigma },\mathbf{\Sigma }_{noise})\rho _{\min }^{-1/2}\left( 
\mathbf{G}_{J^{\ast }}^{\top }\mathbf{G}_{J^{\ast }}\right) ,
\label{eq:YA11}
\end{equation}%
where $\tilde{\rho}(\mathbf{\Sigma },\mathbf{\Sigma }_{noise})=8^{1/4}\left(
\rho ^{4}\left( \mathbf{\Sigma }\right) +\rho ^{4}\left( \mathbf{\Sigma }%
_{noise}\right) \right) ^{1/4}$. Now let $\tilde{Z}=\Vert \mathbf{Y}\Vert
_{\ell _{2}}\leq \rho _{\min }^{-1/2}\left( \mathbf{G}_{J^{\ast }}^{\top }%
\mathbf{G}_{J^{\ast }}\right) \Vert \widetilde{\mathbf{X}}\Vert _{\ell _{2}}$%
. Given our assumptions on the process $\widetilde{X}=X+\mathcal{E}$ it
follows that there exists $\alpha \geq 1$ such that 
\begin{equation}
\Vert \tilde{Z}\Vert _{\psi _{\alpha }}\leq \rho _{\min }^{-1/2}\left( 
\mathbf{G}_{J^{\ast }}^{\top }\mathbf{G}_{J^{\ast }}\right) \left( \Vert
Z\Vert _{\psi _{\alpha }}+\Vert W\Vert _{\psi _{\alpha }}\right) <+\infty ,
\label{eq:YA21}
\end{equation}%
where $Z=\Vert \mathbf{X}\Vert _{\ell _{2}}$ and $W=\Vert \mathbf{\mathcal{E}%
}\Vert _{\ell _{2}}$, with $\mathbf{X}=\left( X\left( t_{1}\right)
,...,X\left( t_{n}\right) \right) ^{\top }$ and $\mathbf{\mathcal{E}}=\left( 
\mathcal{E}\left( t_{1}\right) ,...,\mathcal{E}\left( t_{n}\right) \right)
^{\top }$. Finally, remark that 
\begin{equation}
\left\Vert \Lambda _{J^{\ast }}\right\Vert _{2}\leq \left\Vert \mathbf{\Psi }%
_{J^{\ast }}\right\Vert _{2}+\rho _{\min }^{-1}\left( \mathbf{G}_{J^{\ast
}}^{\top }\mathbf{G}_{J^{\ast }}\right) \left\Vert \mathbf{\Sigma }%
_{noise}\right\Vert _{2}.  \label{eq:eigmax}
\end{equation}%
Hence, using the relations (\ref{eq:YA11}) and (\ref{eq:YA21}), the bound (%
\ref{eq:eigmax}) and Proposition \ref{prop:mp} (with $\mathbf{Y}$ instead of 
$\mathbf{X}$), it follows that there exists a universal constant $\delta
_{\ast }>0$ such that for all $x>0$, 
\begin{equation}
\mathbb{P}\left( \left\Vert \widehat{\mathbf{\Psi }}_{J^{\ast }}-\Lambda
_{J^{\ast }}\right\Vert _{2}\geqslant \tilde{\tau}_{N,s_{\ast }}x\right)
\leqslant \exp \left( -(\delta _{\ast }^{-1}x)^{\frac{\alpha }{2+\alpha }%
}\right) ,  \label{eq:MP3}
\end{equation}%
where $\tilde{\tau}_{N,s_{\ast }}=\max (\tilde{A}_{N,s_{\ast }}^{2},\tilde{B}%
_{N,s_{\ast }})$, with $\tilde{A}_{N,s_{\ast }}=\Vert \tilde{Z}\Vert _{\psi
_{\alpha }}\frac{\sqrt{\log d^{\ast }}(\log N)^{1/\alpha }}{\sqrt{N}}$ and $%
\tilde{B}_{N,s_{\ast }}=\frac{\tilde{\rho}^{2}(\mathbf{\Sigma },\mathbf{%
\Sigma }_{noise})\rho _{\min }^{-1}\left( \mathbf{G}_{J^{\ast }}^{\top }%
\mathbf{G}_{J^{\ast }}\right) }{\sqrt{N}}+\left( \left\Vert \mathbf{\Psi }%
_{J^{\ast }}\right\Vert _{2}+\rho _{\min }^{-1}\left( \mathbf{G}_{J^{\ast
}}^{\top }\mathbf{G}_{J^{\ast }}\right) \left\Vert \mathbf{\Sigma }%
_{noise}\right\Vert _{2}\right) ^{1/2}\tilde{A}_{N,s_{\ast }}$, with $%
d^{\ast }=\min (N,s_{\ast })$. Then, define the event 
\begin{equation*}
\mathcal{B}=\left\Vert \widehat{\mathbf{\Psi }}_{J^{\ast }}-\Lambda
_{J^{\ast }}\right\Vert _{2}\leqslant \tilde{\tau}_{N,s_{\ast }}\delta
_{\star }\left( \log (M)\right) ^{\frac{2+\alpha }{\alpha }},
\end{equation*}%
and note that, for $x=\delta _{\star }\left( \log (M)\right) ^{\frac{%
2+\alpha }{\alpha }}$ with $\delta _{\star }>\delta _{\ast }$, inequality (%
\ref{eq:MP3}) implies that $\mathbb{P}\left( \mathcal{B}\right) \geq
1-M^{-\left( \frac{\delta _{\star }}{\delta _{\ast }}\right) ^{\frac{\alpha 
}{2+\alpha }}}$. Therefore, on the event $\mathcal{A}\cap \mathcal{B}$ (of
probability at least $1-M^{1-\delta }-M^{-\left( \frac{\delta _{\star }}{%
\delta _{\ast }}\right) ^{\frac{\alpha }{2+\alpha }}}$), using inequality (%
\ref{eq:ineqSigma}) and the fact that $\hat{J}=J^{\ast }$ one obtains 
\begin{equation*}
\left\Vert \widehat{\mathbf{\Sigma }}_{\hat{J}}-\mathbf{\Sigma }_{J^{\ast
}}\right\Vert _{2}\leq \rho _{\max }\left( \mathbf{G}_{J^{\ast }}^{\top }%
\mathbf{G}_{J^{\ast }}\right) \tilde{\tau}_{N,s_{\ast }}\delta _{\star
}\left( \log (M)\right) ^{\frac{2+\alpha }{\alpha }},
\end{equation*}%
which completes the proof of Theorem \ref{theo:opnorm}. \hfill $\Box $

\bibliographystyle{apalike}
\bibliography{GroupLassoFinal}

\end{document}